\documentclass[11pt,oneside]{amsart}
\usepackage{amsmath,graphicx,amssymb}
\newtheorem{thm}{Theorem}[section]
\newtheorem{lem}[thm]{Lemma}
\newtheorem{cor}[thm]{Corollary}

\theoremstyle{definition}
\newtheorem{defn}[thm]{Definition}
\newtheorem{conv}[thm]{Convention}
\newtheorem{rem}[thm]{Remark}
\newtheorem{exmp}[thm]{Example}

\newtheorem{question}[thm]{Problem}
\newenvironment{highlight}[2]
  {\vspace{1.5ex}\noindent\textbf{#1}.\textrm{#2}\itshape}
  {\rmfamily\vspace{1.5ex}}

\newcommand{\size}[1]{\ensuremath{\vert #1 \vert}}
\newcommand{\numsides}[1]{\ensuremath{\size{\sides_X{(#1)}}}}
\newcommand{\nummiss}[1]{\ensuremath{\size{\missing_X{(#1)}}}}
\newcommand{\packet}[1]{\ensuremath{\widetilde{#1}}}
\newcommand{\field}[1]{\ensuremath{\mathbb{#1}}}
\newcommand{\script}[1]{\ensuremath{\mathcal{#1}}}
\newcommand{\bfname}[1]{\ensuremath{\mathbf{#1}}}
\newcommand{\scname}[1]{\text{\sf #1}}
\newcommand{\complex}{\script{C}}
\newcommand{\fantypes}{\script{T}}
\newcommand{\missing}{\scname{Missing}}
\newcommand{\edges}{\scname{Edges}}
\newcommand{\cells}{\scname{Cells}}

\newcommand{\sides}{\scname{Sides}}

\newcommand{\area}{\scname{Area}}
\newcommand{\perimeter}{\bfname{P}}
\newcommand{\weight}{\bfname{Wt}}
\newcommand{\relator}{\bfname{R}}
\newcommand{\Z}{\field{Z}}
\newcommand{\R}{\field{R}}
\newcommand{\N}{\field{N}}

\begin{document}
\title[Coherence, Local Quasiconvexity, and Perimeter of $2$-complexes]
      {Coherence, Local Quasiconvexity,\\ and the Perimeter of $2$-complexes}
\author[J.~P.~McCammond]{Jonathan P. McCammond$\ \!{ }^1$}
      \address{Dept. of Math.\\
               University of California, Santa Barbara\\
               Santa Barbara, CA 93106}
      \email{jon.mccammond@math.ucsb.edu}
\author[D.~T.~Wise]{Daniel T. Wise$\ \!{ }^2$}
      \address{Dept. of Math.\\
               McGill University \\
               Montreal, Quebec, Canada H3A 2K6 }
      \email{wise@math.mcgill.ca}
\date{\today}

\begin{abstract}
A group is coherent if all its finitely generated
subgroups are finitely presented.   In this
article we provide a criterion for positively determining the
coherence of a group. This criterion is based upon the
notion of the perimeter of a map between two finite $2$-complexes
which is introduced here.  In the groups to which this theory applies,
a presentation for a finitely generated subgroup can be computed in
quadratic time relative to the sum of the lengths of the generators.
For many of these groups we can show in addition that they are locally
quasiconvex.

As an application of these results we prove that one-relator groups
with sufficient torsion are coherent and locally quasiconvex and we
give an alternative proof of the coherence and local quasiconvexity of
certain $3$-manifold groups.  The main application is to establish the
coherence and local quasiconvexity of many small cancellation groups.
\end{abstract}

\footnotetext[1]{Supported under NSF grant no.\ DMS-99781628}
\footnotetext[2]{Supported as an NSF Postdoctoral Fellow under
grant no.\ DMS-9627506, and partially supported under NSF grant no.\ DMS-9971511.}

\maketitle
\tableofcontents

\section{Introduction}

\subsection{Coherence}
A group is  \emph{coherent} if all its finitely generated
subgroups are finitely presented.  The best known examples of coherent
groups are free groups, surface groups, polycyclic groups, and
$3$-manifold groups.
Outside of these examples, few criteria for determining the coherence
or incoherence of an arbitrary group presentation are known.  Free
groups are easily proven to be coherent by observing that subgroups of
free groups are free and hence finitely presented if they are finitely
generated.  Similarly, surface groups and polycyclic groups are easily
shown to be coherent.  The coherence of the fundamental groups of
$3$-manifolds is a deeper result proved independently by Scott
\cite{Sc73} and Shalen (unpublished).

Theoretical interest in the coherence of various groups
has been prompted in part by a desire to perform calculations.
  Groups in which all of the finitely generated subgroups
have a computable finite presentation are especially amenable to
computer investigation.  The range of possible positive results is
limited by the existence of various counterexamples.  Rips has
produced examples of incoherent word-hyperbolic groups \cite{Ri82},
Bestvina and Brady have produced examples of incoherent right-angled
Artin groups \cite{BeBr97}, and Wise has produced examples of compact
negatively curved $2$-complexes with incoherent fundamental groups
\cite{Wi98}.

Recently, Feighn and Handel proved the remarkable positive result that
the mapping torus of any injective endomorphism of a free group is
coherent \cite{FeHa99}.  Their theorem is related to the coherence of
$3$-manifolds in the following sense.  Many $3$-manifolds arise as
surface bundles over a circle, and their fundamental groups are thus
isomorphic to extensions of surface groups by~$\Z$.  The result of
\cite{FeHa99} shows that extensions of free groups by $\Z$ are also
coherent, thus extending the often successful analogy between free
groups and surface groups.

\subsection{Coherence Results}
The current investigation was primarily motivated by the following
open problem which has remained unresolved for over thirty years with very
little forward progress:

\begin{highlight}{Problem}{(G.Baumslag, \cite{Ba74})}
Is every one-relator group coherent?
\end{highlight}

In this article we describe a criterion
which, if successful, allows one to conclude that the group under
consideration is coherent. The criterion involves a new notion which we
 call ``perimeter''.
Roughly speaking, given a map $Y\rightarrow X$ between 2-complexes,
the perimeter of $Y$ which we denote by $\perimeter(Y)$,
is a measure of how large the ``boundary'' of $Y$ is relative to $X$.
The strategy underlying the results in this paper is that if $Y\rightarrow X$
is unsatisfactory, because for instance, it is not $\pi_1$-injective,
then some $2$-cells can be added to $Y$ which reduce the perimeter,
and after repeating this finitely many times, we obtain a satisfactory map
between $2$-complexes.

We will now state several of our main results.  All
of the undefined terminology, such as perimeter, weighted $2$-complex,
and the various hypotheses will be explained in the course of the
article.  The main coherence result is the following:

\renewcommand{\thethm}{\ref{thm:GCT}}
\begin{thm}[Coherence theorem]
Let $X$ be a weighted $2$-complex which satisfies the perimeter
reduction hypothesis.

A) If $Y$ is a compact connected subcomplex of a cover $\widehat X$ of
$X$, and the inclusion $Y\rightarrow \widehat X$ is not
$\pi_1$-injective, then $Y$ is contained in a compact connected
subcomplex $Y'$ such that $\perimeter(Y') < \perimeter(Y)$.

B) For any compact subcomplex $C\subset \widehat X$, there exists a
compact connected subcomplex $Y$ containing $C$, such that
$\perimeter(Y)$ is minimal among all such compact connected
subcomplexes containing $C$.  Consequently $\pi_1X$ is coherent.
\end{thm}

A geometric consequence of Theorem~\ref{thm:GCT} is the following:

\renewcommand{\thethm}{\ref{thm:compactcore}}
\begin{thm}
Let $X$ be a weighted aspherical $2$-complex which satisfies the
perimeter reduction hypothesis.  If $\widehat X\rightarrow X$ is a
covering space and $\pi_1\widehat X$ is finitely generated, then every
compact subcomplex of $\widehat X$ is contained in a compact core of
$\widehat X$.
\end{thm}
\renewcommand{\thethm}{\thesection.\arabic{thm}}

\subsection{Local Quasiconvexity Results}
A subspace $Y$ of a geodesic metric space $X$ is quasiconvex if
there is an $\epsilon$ neighborhood of $Y$ which contains all of the
geodesics in $X$ which start and end in $Y$.  In group theory, a
subgroup $H$ of a group $G$ generated by $A$ is quasiconvex if
the $0$-cells corresponding to $H$ form a quasiconvex subspace of the
Cayley graph $\Gamma(G,A)$.

The next main result in this article is a criterion which,
 if satisfied, allows one
to conclude that a group is locally quasiconvex, i.e. that all
finitely generated subgroups are quasiconvex. As the reader will observe
from some of the applications described below, most of the groups
which we can show are coherent also satisfy this stronger
criterion, and thus will be locally quasiconvex as well.  The exact
relationship between the two criteria will become clear in the course
of the article.  Our main quasiconvexity theorem is the following:

\renewcommand{\thethm}{\ref{thm:lqc}}
\begin{thm}[Subgroups quasi-isometrically embed]
Let $X$ be a compact weighted $2$-complex.  If $X$ satisfies the
straightening hypothesis, then every finitely generated subgroup of
$\pi_1X$ embeds by a quasi-isometry.  Furthermore, if $\pi_1X$ is
word-hyperbolic then it is locally quasiconvex.
\end{thm}
\renewcommand{\thethm}{\thesection.\arabic{thm}}

\subsection{Some applications}

The statements of the perimeter reduction hypothesis Definition~\ref{def:prh}
and
the straightening hypothesis
(Definition~\ref{def:qrh}) are rather technical, but the strength of the theorems
above can be illustrated easily through some explicit consequences
that we now describe. First of all, as a consequence of Theorem~\ref{thm:GCT} we obtain the following result:

\renewcommand{\thethm}{\ref{thm:equalweights}}
\begin{thm}
Let $W$ be a cyclically reduced word and let $G = \langle a_1,\ldots
\mid W^n \rangle$.  If $n \geq |W|-1$, then  $G$ is coherent.
In particular, for every word $W$, the group $G = \langle a_1,\ldots
\mid W^n \rangle$ is coherent provided that $n$ is sufficiently large.
\end{thm}
\renewcommand{\thethm}{\thesection.\arabic{thm}}

\noindent
With a slightly stronger requirement on the degree of  torsion,
we can obtain the following consequence of Theorem~\ref{thm:lqc}:

\renewcommand{\thethm}{\ref{thm:lqc-one-rel}}
\begin{thm}
Let $G = \langle a_1,\ldots | W^n \rangle$ be a one-relator group
with $n \geq 3|W|$.  Then $G$ is locally quasiconvex.
\end{thm}
\renewcommand{\thethm}{\thesection.\arabic{thm}}

\noindent
A similar result holds for multi-relator groups:

\renewcommand{\thethm}{\ref{thm:powers}}
\begin{thm}[Power theorem]
Let $\langle a_1, \dotsc \mid W_1, \dotsc \rangle$ be a finite
presentation, where each $W_i$ is a cyclically reduced word which is
not a proper power.  If $W_i$ is not freely conjugate to $W_j^{\pm1}$
for $i \neq j$, then there exists a number $N$ such
that for all choices of integers $n_i \geq N$ the group $G = \langle
a_1, \dotsc \mid W_1^{n_1}, \dotsc \rangle$ is coherent.
Specifically, the number
\begin{equation}
  N = 6 \cdot \frac{|W_{\text{max}}|}{|W_{\text{min}}|} \sum |W_i|
\end{equation}

\noindent
has this property, where $W_{\text{max}}$ and $W_{\text{min}}$ denote
longest and shortest words among the $W_i$, respectively.
Moreover, if $n_i >N$ for all $i$, then $G$ is locally quasiconvex.
\end{thm}
\renewcommand{\thethm}{\thesection.\arabic{thm}}

\noindent
A  small-cancellation application with a different flavor is the
following:

\renewcommand{\thethm}{\ref{thm:few-occurrences}}
\begin{thm}
Let $G = \langle a_1,\ldots \mid R_1,\ldots \rangle$ be a small
cancellation presentation which satisfies $C'(1/n)$.  If each $a_i$
occurs at most $n/3$ times among the $R_j$, then $G$ is coherent and
locally quasiconvex.
\end{thm}
\renewcommand{\thethm}{\thesection.\arabic{thm}}

\noindent
More precise applications to additional groups which are important
in geometric group theory can be found in sections~\ref{sec:1relator},
\ref{sec:2-cell-small},~\ref{sec:sc-lqc}, and~\ref{sec:manifold}.
  Many of the individual results
derived in these sections can be summarized by the following
qualitative description:

\begin{highlight}{Qualitative Summary}{}
If a presentation has a large number of generators relative to the sum
of the lengths of the relators, and the relators are relatively long
and sufficiently spread out among the generators, then the group is
coherent and locally quasiconvex.
\end{highlight}

The main triumph of these ideas, is that while we have only partially
solved Baumslag's problem, we have substantially answered the problem
raised by C.T.C.~Wall of whether small-cancellation groups are coherent
\cite{Wall79}.
In a separate paper \cite{McWi-lqc}, we give a much more detailed
application of the strongest results in this paper to
small-cancellation groups. Furthermore, families of examples are
constructed there which show that the applications to
small-cancellation theory are asymptotically sharp.  An
application of our theory towards the local quasiconvexity of
one-relator groups with torsion is given in \cite{HrWi01},
an application towards the coherence of various other one-relator groups is
 given in \cite{McWi-windmills}, and an application towards the subgroup separability of
Coxeter groups is given in \cite{Sch-cox}.

\subsection{Descriptions of the Sections}
We conclude this introduction with a brief section-by-section
description of the article.  The concept of the perimeter of a map is
introduced in section~\ref{sec:perimeter} and it is from this concept
that all of our positive results are derived.  Section~\ref{sec:gch}
shows how this concept leads to the notion of a perimeter reduction
and it also contains our general coherence theorem.
Sections~\ref{sec:attachments}, \ref{sec:2cell}, and \ref{sec:pathCT}
develop the specific case where the perimeter can be reduced through
the addition of a single $2$-cell.  The proofs lead
to procedures which are completely algorithmic.  These algorithmic
approaches are described in section~\ref{sec:algorithms}.  In
particular we show that for the groups included in the $2$-cell
coherence theorem there exists an algorithm to compute an explicit
finite presentation for an arbitrary finitely generated subgroup.
Additionally, we show that the time it takes to produce such a finite
presentation is quadratic in the total length of its set of
generators.  Section~\ref{sec:pathCT} presents a more technical
theorem about coherence using sequences of paths.
In section~\ref{sec:1relator} the theory developed in the first half
of the article is applied to the class of one-relator groups with
torsion.  Similarly, section~\ref{sec:2-cell-small}  presents
some background on small cancellation theory, and gives some applications
of the theory to the coherence of small cancellation groups.

As this work has evolved over the past six years,
two things have become clear:
First of all, a much richer collection of positive results can be obtained
by attaching ``fans'' of 2-cells to reduce perimeter instead of attaching
single $2$-cells.
Secondly, the most significant conclusion of the theory appears to be
 local quasiconvexity rather than coherence, and in fact,
  we know of no word-hyperbolic group which satisfies
our coherence criterion but which is not locally-quasiconvex as well.
The latter part of the paper introduces fans in the context of coherence theorems
that utilize them. Thereafter, fans are employed in the statements and proofs
of the local quasiconvexity theorems.

Fans are defined in section~\ref{sec:fanCT}, where
the general theory is extended  by incorporating fans into the
statements and arguments.
Section~\ref{sec:quasi} contains basic definitions and results about
quasiconvexity which we will need, and
section~\ref{sec:lqc-main} presents our main theorem about local
quasiconvexity.
In section~\ref{sec:sc-lqc} we return to small cancellation groups,
and give additional coherence applications as well as some
local quasiconvexity applications.
Section~\ref{sec:manifold} uses these applications to small
cancellation groups to obtain results about $3$-manifold groups. While
the coherence of $3$-manifold groups has long been known, when
successful, our method gives a different approach towards
understanding the reasons behind this remarkable theorem.
Finally, in section~\ref{sec:other}, we describe theorems and algorithms related
to the finitely generated intersection property and the generalized word problem.

\section{Perimeter}\label{sec:perimeter}
The main goal of this section is to introduce the notion of the
perimeter of a $2$-complex $Y$ relative to a particular map
$\phi\colon Y \rightarrow X$.  We begin with a number of basic
definitions and the definition of the unit perimeter.  In the second
half of the section we broaden the definition to allow for the
introduction of weights.  The weighted perimeter of a $2$-complex is a
measure of the complexity of a map which will be used to prove
coherence and local quasiconvexity theorems throughout the article.

\begin{defn}[Combinatorial maps and complexes]\label{def:combinatorial}
A map $Y\rightarrow X$ between CW complexes is \emph{combinatorial} if
its restriction to each open cell of $Y$ is a homeomorphism onto an
open cell of $X$.  A CW complex $X$ is \emph{combinatorial} provided
that the attaching map of each open cell of $X$ is combinatorial for a
suitable subdivision.
\end{defn}

It will be convenient to be explicit about the cells in a
combinatorial $2$-complex.

\begin{defn}[Polygon]
A \emph{polygon} is a $2$-dimensional disc whose cell structure has
$n$~$0$-cells, $n$~$1$-cells, and one $2$-cell where $n\geq 1$ is a
natural number.  If $X$ is a combinatorial $2$-complex then for each
open $2$-cell $C\hookrightarrow X$ there is a polygon $R$, a
combinatorial map $R\rightarrow X$ and a map $C \rightarrow R$ such
that the diagram
\[\begin{array}{ccc}
C & \hookrightarrow & X\\
\downarrow &\nearrow  & \\
R &          &\\
\end{array}\]
commutes, and the restriction $\partial R \rightarrow X$ is the
attaching map of $C$.  In this article the term \emph{$2$-cell} will
always mean a combinatorial map $R\rightarrow X$ where $R$ is a
polygon. The corresponding \emph{open $2$-cell} is the image of the
interior of $R$.

A similar convention applies to $1$-cells.  Let $e$ denote the graph
with two $0$-cells and one $1$-cell connecting them.  Since
combinatorial maps from $e$ to $X$ are in one-to-one correspondence
with the characteristic maps of $1$-cells of $X$, we will often refer
to a map $e \rightarrow X$ as a {\em $1$-cell} of $X$.
\end{defn}

\begin{defn}[Standard $2$-complex]
In the study of infinite groups, the most commonly considered
combinatorial $2$-complexes correspond to presentations.  Recall that
the \emph{standard $2$-complex} of a presentation is formed by taking
a unique $0$-cell, adding a labeled oriented $1$-cell for each
generator, and then attaching a $2$-cell along the closed
combinatorial path corresponding to each relator.
\end{defn}

\begin{conv}\label{conv:combinatorial}
Unless noted otherwise, all complexes in this article are
combinatorial $2$-complexes, and all maps between complexes are
combinatorial maps.  In addition, we will avoid certain technical
difficulties by always assuming that all of the attaching maps for the
$2$-cells are immersions.  For $2$-complexes with a unique $0$-cell,
this is equivalent to allowing only cyclically reduced relators in the
corresponding presentation.
\end{conv}

\begin{defn}[Basic definitions]\label{def:basic}
A local injection between topological spaces is an \emph{immersion.}
If $\phi\colon Y \rightarrow X$ is an immersion on $Y \setminus
Y^{(0)}$, then $\phi$ is a \emph{near-immersion}.  If $\phi\colon
Y\rightarrow X$ is an immersion on $Y^{(1)}$ then $\phi$ is  a
\emph{$1$-immersion}.  Let $\phi_* \colon \pi_1Y \rightarrow \pi_1X$
be the induced homomorphism between fundamental groups.  The map
$\phi$ is \emph{$\pi_1$-injective} [respectively
\emph{$\pi_1$-surjective}] if $\phi_*$ is injective [surjective].
Finally, if $\phi\colon Y \rightarrow X$ and $\psi:Z \rightarrow
X$ are fixed maps, then a map $\rho:Z \rightarrow Y$ is  a
\emph{lift of $\psi$} or \emph{a lift of $Z$ to $Y$} whenever the
composition $\phi \circ \rho = \psi$.
\end{defn}

\begin{defn}[Path and cycle]
A \emph{path} is a map $P\rightarrow X$ where $P$ is a subdivided
interval or a single $0$-cell.  In the latter case, $P$ is
a \emph{trivial path}.  A \emph{cycle} is a map $C\rightarrow X$
where $C$ is a subdivided circle.  Given two paths $P\rightarrow
X$ and $Q\rightarrow X$ such that the terminal point of $P$ and
the initial point of $Q$ map to the same $0$-cell of $X$, their
concatenation $PQ\rightarrow X$ is the obvious path whose domain
is the union of $P$ and $Q$ along these points.  The path
$P\rightarrow X$ is a \emph{closed path} provided that the
endpoints of $P$ map to the same $0$-cell of $X$.  A path or cycle
is \emph{simple} if the map is injective on $0$-cells.  The
\emph{length} of the path $P$ or cycle $C$ is the number of
$1$-cells in the domain and it is denoted by $\size{P}$ or
$\size{C}$.  The \emph{interior} of a path is the path minus its
endpoints.  In particular, the $0$-cells in the interior of a path
are the $0$-cells other than the endpoints.  A \emph{subpath} $Q$
of a path $P$ [or a cycle $C$] is given by a path $Q \rightarrow P
\rightarrow X$ [$Q \rightarrow C \rightarrow X$] in which distinct
$1$-cells of $Q$ are sent to distinct $1$-cells of $P$ [$C$].
Notice that the length of a subpath is at most that of the path
[cycle] which contains it.  Finally, note that any nontrivial
closed path determines a cycle in the obvious way.  Finally, when
the target space is understood we will often just refer to
$P\rightarrow X$ as the path $P$.
\end{defn}

\begin{conv}
The letters $X$ and $Y$ will always refer to spaces, $R$ will
always refer to a closed $2$-cell, and $P$ will always denote a
path.  We follow the convention that lowercase letters (such as
$x$, $y$, and $r$) refer to specified $1$-cells in the space
denoted by the corresponding uppercase letter.  Thus $r$ is a
$1$-cell in the boundary of the $2$-cell $R$ and $x$ is a $1$-cell
in the space $X$.

We will be very interested in examining the behavior of maps
 and spaces along selected $1$-cells.
Accordingly, the pair $(Y,y)$
will denote a space together with a chosen $1$-cell
in that space, and we will write
$\rho:(R,r) \rightarrow (X,x)$ to denote a
map $\rho:R\rightarrow X$ with the property that $\rho(r)=x$.
\end{conv}

\begin{defn}[Side]\label{def:sides}
Let $X$ be a fixed $2$-complex, and let $R$ be a $2$-cell of $X$.  Let
$r$ be a $1$-cell in $\partial R$ and let $x$ be the image of $r$ in
$X$.  The pair $(R,r)$ will then be called a \emph{side of a $2$-cell of
$X$ which is present at $x$}.  The collection of all sides of $X$
which are present at $x$ will be denoted by $\sides_X(x)$, and the
full collection of sides of $2$-cells of $X$ which are present at
$1$-cells of $X$ will be denoted by $\sides_X$.
Notice that saying a side $(R,r)$ is present at $x$ is equivalent to
saying that the map $R\rightarrow X$ extends to a map $(R,r)
\rightarrow (X,x)$.  Notice also that if $r$ and $r'$ are distinct
$1$-cells of $R$ which are mapped to the same $1$-cell $x$ of $X$,
then $(R,r)\rightarrow (X,x)$ and $(R,r') \rightarrow (X,x)$ are
distinct sides at $x$, even though $r$ and $r'$ come from the same
$2$-cell $R$ and are mapped to the same $1$-cell $x$.  Thus a $2$-cell
$R$ whose boundary has length~$n$ will have exactly~$n$ distinct sides
in $\sides_X$.  Alternatively, the elements of $\sides_X(x)$ can be
viewed as the connected components in $(X - x) \cap B$, where $B$ is a
small open ball around a point in the interior of $x$.

Next, let $\phi\colon Y \rightarrow X$ be a map, let $(R,r)$ be a side
of $X$ which is present at $x$, and let $y$ be a $1$-cell of $Y$ with
$\phi(y)= x$.  We say that the side $(R,r) \rightarrow (X,x)$ is
\emph{present at $y$} if the map $(R,r) \rightarrow (X,x)$
factors through a map $(R,r) \rightarrow (Y,y)$ as indicated in the
following commutative diagram:
\[\begin{array}{ccc}
&& (Y,y)\\
&\nearrow & \downarrow\\
(R,r) &\rightarrow&(X,x)\\
\end{array}\]
Specifically, there must exist a map $\rho: (R,r) \rightarrow (Y,y)$
such that $\phi \circ \rho$ is the map $(R,r) \rightarrow (X,x)$.  If
the map $(R,r)\rightarrow (X,x)$ does not factor through $\phi$ then
$(R,r)$ is said to be \emph{missing at $y$}.  The set of all sides of
$X$ which are present at $y$ will be denoted by $\sides_X(y)$, while
the set of all sides of $X$ which are missing at $y$ will be denoted
$\missing_X(y)$.
\end{defn}

\begin{rem}\label{rem:sides}
It is important to notice that the definitions of the sets
$\sides_X(y)$ and $\sides_X(x)$ both refer to the sides of $2$-cells
of the complex $X$.  In particular, if $\phi(y) = x$ then $\sides_X(y)
\subset \sides_X(x)$, and $\sides_Y(y)$ is not comparable with either
of these since it is a subset of $\sides_Y$.  Moreover, $\sides_X(y)$
can be smaller than $\sides_Y(y)$ if the map $\phi$ is not a near-immersion.
  In fact, $\phi$ is a near-immersion if and only if
$|\sides_X(y)| = |\sides_Y(y)|$ for all $1$-cells $y \in Y$.
\end{rem}

\begin{defn}[Unit perimeter]\label{def:perimeter}
Let $\phi\colon Y \rightarrow X$ be a combinatorial map between
$2$-complexes.  We define the \emph{unit perimeter} of $\phi$ to be
\begin{equation}
 \perimeter(\phi) \ =
 \sum_{y \in \edges(Y)} \nummiss{y} \  =
 \sum_{y \in \edges(Y)} \numsides{\phi(y)} - \numsides{y}
\end{equation}

\noindent For each $1$-cell $y$ of $Y$, we can either count the sides
of $X$ at $x=\phi(y)$ which are missing at $y$, or else we can count
the number of sides of $X$ that are present at $x$ and then subtract
off those which are also present at $y$.  From the first description
it is clear that the perimeter of $\phi$ is nonnegative.
\end{defn}

The following examples will illustrate these distinctions.  In
particular, they will illustrate the significance of the maps
$\phi\colon Y \rightarrow X$ and $R \rightarrow X$, respectively.

\begin{exmp}\label{exmp:cells}
Let $X$ be the complex formed by attaching two squares along a common
$1$-cell $x$, and let $Y$ be another complex which is isomorphic to
$X$ with common $1$-cell $y$. Let $\phi\colon Y\rightarrow X$ be an
isomorphism, and let $\psi:Y\rightarrow X$ be a map which sends $y$ to
$x$ but which folds the two squares of $Y$ to the same square of
$X$. Observe that $\perimeter(\phi\colon Y\rightarrow X) =0$ but
$\perimeter(\psi:Y\rightarrow X)=1$.
\end{exmp}


\begin{exmp}\label{exmp:sides}
Let $X$ be the standard $2$-complex of the presentation $\langle
a,b\mid (aab)^3 \rangle$ and let $R\rightarrow X$ be the unique
$2$-cell of $X$.  Notice that $\partial R$ wraps three times around
the path $aab$ in $X$, and that there are exactly six sides present at
the $1$-cell labeled~$a$ and exactly three sides present at the
$1$-cell labeled~$b$.

Inside the universal cover $\widetilde{X}$ of $X$, one can find three
distinct $2$-cells which share the same boundary cycle.  Let $Y\subset
X$ be the union of two of these three $2$-cells and define $\phi\colon
Y \rightarrow X$ to be the composition $Y\hookrightarrow \widetilde{X}
\rightarrow X$.
Observe that $Y$ is a sphere and $\phi$ is an immersion.  If $y$ is a
$1$-cell labeled $b$ in $Y$ and $x$ is its image under $\phi$, then
$\size{\sides_Y(y)} = 2$, $\size{\sides_X(y)} = 2$, and
$\size{\sides_X(x)} = 3$.  Thus $\size{\missing_X(y)} = 1$.
To see the importance of the map $R\rightarrow X$, let $r$ be a
$1$-cell labeled $b$ in $R$ and let $(R,r) \rightarrow (X,x)$ be the
corresponding map of pairs.  There is a $1$-cell $y$ in $Y$ such that
$(R,r)$ is missing at $y$ even though there are two distinct maps from
$(R,r)$ to $Y$ which send $r$ to $y$ and which when composed with
$\phi$ agree with the map $R\rightarrow X$ on $\partial R$.  The side
$(R,r)$ is missing from $y$ because neither of these maps agree with
$R\rightarrow X$ on the interior of $R$.
\end{exmp}

\begin{exmp}[$\Z \times \Z \times \Z$]\label{exmp:zzz}
Let $G = \langle a,b,c \mid [a,b] = [a,c] = [b,c] = 1 \rangle$ be the
standard presentation of the free abelian group on three generators
and let $X$ be the standard $2$-complex corresponding to this
presentation. The universal cover $\widetilde{X}$ of $X$ is usually
thought of as the points of $\R^3$ with $x \in \Z$ or $y \in \Z$ or $z
\in \Z$.  That is, $\widetilde{X}$ is isomorphic to the union of the
integer translates of the $x$-$y$, $y$-$z$, and $x$-$z$ planes.  The
$1$-cells and $2$-cells of $\widetilde{X}$ are unit intervals and unit
squares.  Moreover, $\widetilde{X}$ is labeled so that $1$-cells
parallel to the $x$-axis are labeled by the generator $a$ and directed
in the positive $x$~direction.  Similarly, $1$-cells parallel to the
$y$-axis and $z$-axis are labeled by $b$ and $c$ directed in the positive
$y$-direction and $z$-direction respectively.

\begin{figure}\centering
\includegraphics[width=.3\textwidth]{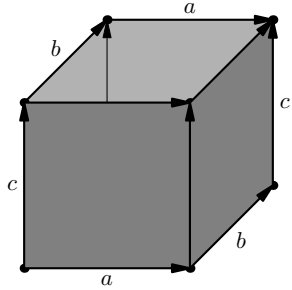}
\caption{The space $Y$ of Example~\ref{exmp:zzz} \label{fig:zzz}}
\end{figure}

If $Y$ is a 1-by-1-by-1 box with four walls, a bottom and no top (see
Figure~\ref{fig:zzz}), and $\phi$ is the obvious embedding of $Y$ into
$\widetilde{X}$, then the perimeter of $Y$ is as follows.  Each
$1$-cell $e$ along the top of the open box contributes a perimeter of
$3$ corresponding to the sides of the three squares which are incident
at $\phi(e)$ in $\widetilde{X}$ but which do not lift to $Y$.  The
vertical $1$-cells along the sides of the box contribute a perimeter
of $2$ each, as do the $1$-cells on the bottom of the box.  The total
perimeter of $Y$ is $28$ sides of squares.

Next, consider a more typical example of a map which is not an
embedding.  Let $Y\rightarrow X$ denote the composition map
$Y\rightarrow \widetilde{X} \rightarrow X$.  The exact same count
shows that the perimeter of $Y \rightarrow X$ is also $28$.
\end{exmp}

We will now define a more flexible notion of perimeter which employs a
weighting on the sides of the $2$-cells in $X$.  The weighted
perimeter of $\phi\colon Y \rightarrow X$ is intuitively just the sum
of the weights of the corresponding missing sides.  We will now make
this notion more precise.

\begin{defn}[Weighted perimeter]\label{def:weighted-perimeter}
A \emph{weight function} on a $2$-complex $X$ is a function of the form
$\weight:\sides(X)\rightarrow\R$.  For most of our applications we
will require that the weight of a side be nonnegative.  Let
$\weight:\sides(X)\rightarrow\R$ be a weight function on $X$, and let
$\phi\colon Y \rightarrow X$ be a combinatorial map of $2$-complexes.
The \emph{weighted perimeter} of $\phi\colon Y \rightarrow X$ is the sum
of the weights of the sides of $X$ which are missing at $1$-cells of
$Y$.  More precisely, the weighted perimeter is defined to be the
following double sum:
\begin{equation}\label{eq:weighted-perimeter}
\perimeter(\phi\colon Y\rightarrow X) \  =
 \sum_{y \in \edges(Y)} \hskip 1mm
 \sum_{(R,r) \in \missing_X(y)}
 \weight \big( (R,r) \big)
\end{equation}

\noindent
Notice that the weighted perimeter is equivalent to the unit perimeter
when each side is assigned a weight of $1$.  This weight function will
be called the \emph{unit weighting}.  Note also that if the assigned
weights are nonnegative, then the weighted perimeter of $\phi$ will be
nonnegative.  When the map $\phi$ is understood we will write
$\perimeter(Y)$ or $\perimeter(Y\rightarrow X)$ for
$\perimeter(\phi\colon Y\rightarrow X)$.
\end{defn}

\begin{defn}[Weighted $2$-complex]
A $2$-complex $X$ is a \emph{weighted $2$-complex} if each of the sides
of $X$ has been assigned a nonnegative integer weight, the perimeter of
each $1$-cell is finite, and the weight of each $2$-cell is positive.
\end{defn}

Although the definition of a weighted $2$-complex adopts the
requirement of an integer weighting, it is sufficient, and often quite
natural, to use a finite set of rationals.  In the theorems which
follow, a successful real-valued weight function can always be
approximated by a successful rational weight function.  After clearing
the denominators, we then would obtain a successful integer-valued
weight function.  Thus there is no real loss of generality in assuming
that the values of the weights are integers.

Perimeters of weighted $2$-complexes satisfy the following useful
property:

\begin{lem}\label{lem:surj}
Let $X$ be a weighted $2$-complex and consider maps $\rho:Z
\rightarrow Y$, $\phi\colon Y \rightarrow X$, and $\psi = \phi \circ
\rho: Z \rightarrow X$.  If $\rho$ is surjective then
$\perimeter(Z\rightarrow X) \geq \perimeter(Y\rightarrow X)$.
\end{lem}

\begin{proof}
First notice that $\missing_X(z) \supset \missing_X(y)$ whenever
$\rho(z) = y$ (or equivalently that $\sides_X(z) \subset
\sides_X(y)$).  Since $\rho$ is surjective, for every $1$-cell $y$ in
$Y$ we can select a $1$-cell $z$ in $Z$ with $\rho(z) = y$.  It is also
clear that the $z$'s chosen for distinct $y$'s are themselves
distinct.  Thus the terms in the sum for $\perimeter(\phi)$ can be
identified with distinct terms in the sum for $\perimeter(\psi)$.
Finally, since the weights are nonnegative it follows that
$\perimeter(\psi) \geq \perimeter(\phi)$.
\end{proof}

\begin{defn}[Induced weights]
Given a weight function on $X$ there is also an induced value assigned
to each of the $1$-cells and $2$-cells of $X$. We define the
\emph{perimeter of a $1$-cell} $x$ in $X$ to be the sum of the weights
assigned to the sides in $\sides_X(x)$.  This agrees with our earlier
definition of perimeter in the sense that it is the weighted perimeter
of the map $\phi\colon x \rightarrow X$ which sends the single
$1$-cell to $x$ in $X$.  In particular, it measures the weights of the
sides which are {\em not} present when the $1$-cell $x$ is considered
in isolation.

We define the \emph{weight of a $2$-cell} $R$ in $X$ to be the sum of
the weights assigned to the sides of the form $(R,r)$ for some $r$ in
$\partial R$.  The sum of the weights of the sides of a $2$-cell, on
the other hand, is called a weight since it is the sum of the weights
of sides which are present in the $2$-cell itself and it ignores the
weights of the other sides which are incident at $1$-cells in its
boundary. Formally, we have the equations
\begin{equation}\label{eq:side}
  \perimeter(x) \ = \sum_{(R,r)\in\sides_X(x)} \weight((R,r))
\end{equation}
\begin{equation}\label{eq:2-cell}
  \weight(R) \ = \sum_{r \in \text{Edges}(\partial R)} \weight((R,r))
\end{equation}
\end{defn}

If $Y$ is compact and the map $\phi\colon Y\rightarrow X$ is a near-immersion,
 then the perimeter of $Y$ can be calculated from the
perimeters of its $1$-cells and the weights of its $2$-cells.
Specifically we have the following result.

\begin{lem}\label{lem:perimeter-rewritten}
If $X$ is a weighted $2$-complex, $Y$ is compact, and the map
$\phi\colon Y\rightarrow X$ is a near-immersion, then
\begin{equation}\label{eq:weight-stuff}
\perimeter(\phi) = \sum_{y \in \edges(Y)} \hskip 1mm
\perimeter(\phi(y)) - \sum_{S \in \cells(Y)}
\hskip 1mm \weight(\phi(S))
\end{equation}
\end{lem}

\begin{proof}
The second summation is the one which requires the immersion
hypothesis.  By Remark~\ref{rem:sides}, the restriction on $\phi$
implies that the sides of $X$ which are present at $y$ are in
one-to-one correspondence with the sides of $Y$ which are present at
$y$.  These sides of $Y$ can then be collected together according to
the $2$-cell in $Y$ to which they belong, and then the sum of the
weights of the sides of a particular $2$-cell $S$ in $Y$ can be
rewritten as the weight of the $2$-cell in $X$ which is the image of
$S$ under $\phi$.
\end{proof}

The following example illustrates these types of calculations.

\begin{exmp}[Weighted $\Z \times \Z \times
 \Z$]\label{exmp:weighted zzz} Let $X$, $Y$, and $\phi\colon Y
\rightarrow X$ be the spaces and maps described in
Example~\ref{exmp:zzz}.  Specifically, let $X$ denote the standard
$2$-complex of the presentation
\[\langle a,b,c \mid aba^{-1}b^{-1},aca^{-1}c^{-1}, bcb^{-1}c^{-1}
\rangle\]
Denote the three $2$-cells of $X$ by $R_1$, $R_2$ and $R_3$, and
observe that each of them has four sides.  Since the four letters of
each defining word are in one-to-one correspondence with the four
sides of the $2$-cell, the weights of the sides of the $2$-cell can be
indicated by a sequence of four numbers.  If we assign weights to the
sides of $2$-cells of $X$ via the sequences $(1,2,3,4), (1,2,0,0)$,
and $(1,3,5,0)$, then the reader can verify that $\weight(R_1) = 10$,
$\weight(R_2) = 3$, $\weight(R_3) = 9$, $\perimeter(a) = 5$,
$\perimeter(b) = 12$, and $\perimeter(c) = 5$.  Since the map $Y
\rightarrow X$ is an immersion, by Lemma~\ref{lem:perimeter-rewritten}
the weighted perimeter of the map $\phi\colon Y\rightarrow X$ can be
calculated as follows:
\[\perimeter(Y) = 4\perimeter(a) + 4\perimeter(b) + 4\perimeter(c) -
\weight(R_1) - 2\weight(R_2) - 2\weight(R_3) = 54\]
\end{exmp}

\section{Coherence theorem}\label{sec:gch}
In this section we describe a general framework for showing that
groups are coherent and then employ the notion of perimeter to
state our main hypothesis and to prove our main coherence theorem.

\begin{defn}[Complexity function]\label{def:complexity-fns}
Let $X$ be a fixed $2$-complex, let $(N,<)$ be a well-ordered set, and
let $\complex$ be a function which assigns an element of $N$ to each
map $\phi\colon Y \rightarrow X$ with a compact domain.  The function
$\complex$ is a \emph{complexity function} for $X$, and the value
$\complex(\phi)$ is the \emph{complexity of the map $\phi$}.  In
practice, $(N,<)$ will either be $\R^+$ with the usual ordering, or
$N$ will be $\R^+ \times \R^+$ and $<$ is defined so that $(a,b) <
(c,d)$ if either $a < c$ or $a = c$ and $b <d$.  This is the usual
lexicographic ordering on ordered pairs.
\end{defn}

\begin{defn}[Reduction method]\label{reduction-method}
Let $X$ be a fixed $2$-complex and let $\complex$ be a complexity
function for $X$.  If for all compact spaces $Y$ and maps $\phi\colon Y
\rightarrow X$ such that $\phi$ is not already $\pi_1$-injective,
there is a ``procedure'' (in any sense of the word) which produces a
compact space $Z$ and a map $\rho:Z\rightarrow X$ such that
$\rho_*(\pi_1Z) = \phi_*(\pi_1Y)$ and such that
$\complex(\rho) < \complex(\phi)$, then this procedure will be
called a \emph{reduction method} for $\complex$.
\end{defn}

\begin{rem}\label{rem:compact}
Notice that for every finitely generated subgroup $H$ in $\pi_1X$
there exists a compact space $Y$ and map $\phi\colon Y\rightarrow X$
such that the image of $\pi_1Y$ under $\phi_*$ is exactly $H$.
One procedure for creating $Y$ and $\phi$ goes as follows:
Suppose that $H$ is generated by $n$~elements of $\pi_1X$ and
represent each of these generators by a closed path in the $1$-skeleton
of $X$ starting at the basepoint.  Next let $Y$ be a bouquet of
$n$~circles, and after subdividing $Y$, define $\phi\colon Y\rightarrow X$
so that the restriction of $\phi$ to the $i$-th subdivided circle is
identical to the $i$-th closed path.

Alternatively, we could let $\widehat X$ be the based covering space of
$X$ corresponding to the subgroup $H$, and let $Y$ be the union of the
based lifts to $\widehat X$ of a finite set of closed based
paths representing the generators of $H$ in $\pi_1X$.
 It is clear that both of these
constructions yield $\pi_1$-surjective maps.  The latter has the
advantage of being an immersion.
\end{rem}

The following theorem is the philosophical basis for the coherence
results in this paper.

\begin{thm}\label{thm:complexity=>coherent}
Let $X$ be a fixed space and let $\complex$ be a complexity function
for $X$.  If there is a reduction method for $\complex$
 then $\pi_1X$ is coherent.
\end{thm}

\begin{proof}
Let $H$ be an arbitrary finitely generated subgroup of $\pi_1X$. By
Remark~\ref{rem:compact}, there is at least one combinatorial map
$\phi\colon Y \rightarrow X$ such that $Y$ is compact and
$\phi_*(\pi_1Y) = H$.  If $\phi$ is not $\pi_1$-injective then there
is another combinatorial map with the same properties which has a
strictly lower complexity.  Since $(N,<)$ is well ordered, there
cannot be an infinite sequences of reductions.  Hence the process of
replacing one combinatorial map with another must terminate at a
$\pi_1$-injective combinatorial map $\rho:Z \rightarrow X$ where $Z$
is compact and $\rho_*(\pi_1Z) = H$. Since $\rho$ is
$\pi_1$-injective, $\pi_1Z$ is itself isomorphic to $H$.  Since $Z$ is
compact, $H$ is finitely presented. In particular, a standard
$2$-complex for a finite presentation of $H$ can be obtained by
contracting a maximal tree in $Z^{(1)}$.
\end{proof}

Note that if the reduction method for $\complex$ is constructive, then
the proof of Theorem~\ref{thm:complexity=>coherent} can be used as an
algorithm to effectively compute finite presentations for finitely
generated subgroups.  Many of the reduction methods we introduce are
in fact constructive and in Section~\ref{sec:algorithms} we explicitly
describe a resulting algorithm.

\begin{rem}\label{coherence=>complexity}
The converse of Theorem~\ref{thm:complexity=>coherent} is also true in
the following sense.  Given a space $X$ with a coherent fundamental
group, we define the complexity of a map $\phi\colon Y \rightarrow X$
where $Y$ is compact to be the minimum number of $2$-cells which must
be added to $Y$ to yield a $\pi_1$-injection.  It is easy to see that
this is indeed a complexity function, that the ``procedure'' of adding
one of the necessary $2$-cells is a method of reducing the complexity,
and that there can be no infinite sequences of reductions.
\end{rem}

We will now specialize to the case where weighted perimeter is used to
measure the complexity of a map.

\begin{defn}[Reduction hypothesis]\label{def:gprh}
Let $X$ be a weighted $2$-complex.  It will satisfy the \emph{perimeter
reduction hypothesis} if for any compact and connected space $Y$ and
for any based $1$-immersion $\phi\colon Y \rightarrow X$ which is not
$\pi_1$-injective, there exists a based map $\phi^+:Y^+\rightarrow X$
and a commutative diagram
\[\begin{array}{ccc}
Y&\rightarrow &X\\ \downarrow&\nearrow&\\ Y^+&&\\
\end{array}\]
such that $Y^+$ is compact and connected, $\perimeter(Y^+) <
\perimeter(Y)$ and $Y^+\rightarrow X$ has the same $\pi_1$-image as
$Y\rightarrow X$.
Typically, the final requirement that
$\phi^+_*(\pi_1Y^+) = \phi_*(\pi_1Y)$ is deduced from a more
stringent requirement that $Y\rightarrow Y^+$ is $\pi_1$-surjective.
\end{defn}

\begin{thm}[Coherence theorem]\label{thm:GCT}
Let $X$ be a weighted $2$-complex which satisfies the perimeter
reduction hypothesis.

A) If $Y$ is a compact connected subcomplex of a cover $\widehat X$ of
$X$, and the inclusion $Y\rightarrow \widehat X$ is not
$\pi_1$-injective, then $Y$ is contained in a compact connected
subcomplex $Y'$ such that $\perimeter(Y') < \perimeter(Y)$.

B) For any compact subcomplex $C\subset \widehat X$, there exists a
compact connected subcomplex $Y$ containing $C$, such that
$\perimeter(Y)$ is minimal among all such compact connected
subcomplexes containing $C$.  Consequently $\pi_1X$ is coherent.
\end{thm}

\begin{proof}
To prove Statement A, suppose that $Y$ is connected and compact but
the inclusion map $Y\rightarrow \widehat X$ is not $\pi_1$-injective.
Then by the perimeter reduction hypothesis, there exists a commutative
diagram
\[\begin{array}{ccc}
Y&\rightarrow &X\\ \downarrow&\nearrow&\\ Y^+&&\\
\end{array}\]
such that $Y^+$ is compact and connected, such that $\perimeter(Y^+) <
\perimeter(Y)$, and such that $Y^+\rightarrow X$ has the same
$\pi_1$-image as $Y\rightarrow X$.  Observe that $Y^+\rightarrow X$
lifts to a map $Y^+\rightarrow \widehat X$ which extends the lift of
$Y$ to $\widehat X$.  Let $Y'$ denote the image of $Y^+$ in $\widehat
X$.  By Lemma~\ref{lem:surj}, $\perimeter(Y')\leq \perimeter(Y^+) <
\perimeter(Y)$.

Statement B follows immediately from the fact that the perimeters of
compact subcomplexes containing $C$ are nonnegative integers. To see
that $\pi_1X$ is coherent, let $\widehat X$ be a based cover of $X$
such that $\pi_1\widehat X$ is finitely generated.  Observe that there
exists a based compact connected subspace $C\subset \widehat X$ whose
inclusion induces a $\pi_1$-surjection.
Let $Y\subset \widehat X$ denote a compact connected subspace
containing $C$ such that $\perimeter(Y)$ is minimal among all such
compact connected subspaces. Then $Y\rightarrow \widehat X$ is
$\pi_1$-surjective because $C\subset Y$, and $Y\rightarrow \widehat X$
is $\pi_1$-injective by Statement A.
\end{proof}

In the remainder of the article we will provide three conditions which
will imply the perimeter reduction hypothesis: the $2$-cell reduction
hypothesis, the path reduction hypothesis, and the fan reduction
hypothesis.  These three hypotheses are more concrete than the perimeter
reduction hypothesis and thus tend to be more useful in establishing
the coherence of specific presentations.

We conclude this section by noting that Theorem~\ref{thm:GCT} very
nearly shows that finitely generated covers of complexes satisfying
the perimeter reduction hypothesis have compact cores.

\begin{defn}[Core]\label{def:core}
A subcomplex $Y$ of the complex $Z$ is a \emph{core} of $Z$ if
the inclusion map $Y \hookrightarrow Z$ is a homotopy equivalence.
Since $Y$ and $Z$ are CW-complexes, $Y$ is a core of $Z$ if and only
if there is a strong deformation retraction from $Z$ to $Y$, which is
true if and only if the map $Y \rightarrow Z$ induces an isomorphism
on all of the homotopy groups (\cite{Wh78}).
Note that when $Z$ is an aspherical $2$-complex, $Y$ will be a core
for $Z$ if and only if the inclusion of $Y$ induces a
$\pi_1$-isomorphism.  Indeed, if $Y\hookrightarrow Z$ is
$\pi_1$-injective, then the based component of the preimage of $Y$ in
the universal cover $\widetilde Z$ is clearly isomorphic to the
universal cover $\widetilde Y$ of $Y$. But then $\pi_2(\widetilde Y) =
H_2(\widetilde Y) \subset H_2(\widetilde Z) = 0$, and so we see that
$\widetilde Y$ and thus $Y$ is aspherical.
\end{defn}

\begin{thm}\label{thm:compactcore}
Let $X$ be a weighted aspherical $2$-complex which satisfies the
perimeter reduction hypothesis.  If $\widehat X\rightarrow X$ is a
covering space and $\pi_1\widehat X$ is finitely generated, then every
compact subcomplex of $\widehat X$ is contained in a compact core of
$\widehat X$.
\end{thm}

\begin{proof}
This follows immediately from Theorem~\ref{thm:GCT} and
Definition~\ref{def:core}.
\end{proof}

The existence of a compact core in a $2$-complex is a nontrivial fact.
For example, there exists a covering space $\widehat X$ of a
$2$-complex $X$ with a single $2$-cell, such that $\pi_1\widehat X$ is
finitely generated, but $\widehat X$ has no compact core.  See
\cite{Wi-nocore} for details.

The restriction in Theorem~\ref{thm:compactcore} that $X$ be
aspherical is not particularly stringent since one of our main sources
of applications will be small cancellation complexes, and small
cancellation complexes in which none of the $2$-cells are attached by
proper powers are known to be aspherical \cite[\S III.11]{LySch77}.
Roughly speaking, a $2$-cell $R$ is attached by a proper power if
$\partial R \rightarrow X$ is obtained by traversing a closed path in
$X$ two or more.  See Definition~\ref{def:exponent}.

\begin{question}[Asphericity]
It appears likely that if $X$ is a compact $2$-complex which satisfies
the perimeter reduction hypothesis then $\pi_1X$ acts properly
discontinuously on a contractible $2$-complex.  We have been unable to
decide whether this is the case.
\end{question}

\section{Attachments}\label{sec:attachments}
By Theorem~\ref{thm:GCT}, the fundamental group of a weighted
$2$-complex $X$ is coherent if there is a method for reducing the
perimeter of the maps $Y \rightarrow X$ which are not
$\pi_1$-injections.  One of the simplest possibilities is where the
perimeter is reduced through the attachment of a single $2$-cell, a
possibility which will be examined in detail in
Section~\ref{sec:2cell}.  In this section, we provide the definitions
and results about paths, $2$-cells and attachments which will be
needed.

We will now describe the two elementary ways of changing a path $P
\rightarrow X$: to remove a backtrack and to push across a $2$-cell.

\begin{defn}[Removing backtracks]\label{def:removing-backtracks}
If $P \rightarrow X$ contains a subpath of the form $ee^{-1}$ where
$e$ is a $1$-cell of $X$, then there is another path $P'\rightarrow X$
obtained by simply removing these two $1$-cells from the path.  Such a
change is called \emph{removing a backtrack}.  Notice that the paths
$P\rightarrow X$ and $P'\rightarrow X$ are homotopic relative to their
endpoints, that a path $P$ is immersed if and only if it has no
backtracks to remove, and that removing a backtrack reduces the length
of the path.
\end{defn}

\begin{defn}[Complement]
Let $R\rightarrow X$ be a $2$-cell, and let $Q$ be a subpath of
$\partial R$.  There exists a unique subpath $S$ of $\partial
R$, called the \emph{complement of $Q$ in $R$}, such that the
concatenation $QS^{-1}$ is a closed path which corresponds to the
boundary cycle $\partial R$.  Note that if $|Q|=|\partial R|$, then
$S$ is a trivial path.
\end{defn}

\begin{defn}[Pushing across a $2$-cell]\label{def:pushing-across}
Let $P \rightarrow X$ be a path, let $R \rightarrow X$ be a $2$-cell,
and let $Q$ be a subpath of both $P$ and $\partial R$, so that we have
the following commutative diagram:
\[\begin{array}{ccc} Q & \rightarrow & P \\ \downarrow &
& \downarrow \\ R & \rightarrow & X \end{array}\]
Let $S$ be the complement of $Q$ in $R$, and observe that since $S$
and $Q$ have the same endpoints in $X$, we can form a new path $P'$ by
substituting $S$ for the subpath $Q$ of $P$.  In particular, if the
path $P\rightarrow X$ is the concatenation of a path $P_1$ followed by
the path $Q$ followed by a path $P_2$, then the modified path
$P'\rightarrow X$ is the concatenation $P_1SP_2$.  The replacement of
$Q \rightarrow X$ by $S \rightarrow X$ is called \emph{pushing across
the $2$-cell $R \rightarrow X$}.  It is clear that if $P'$ is obtained
from $P$ by pushing across a $2$-cell, then $P$ and $P'$ will be
homotopic relative to their endpoints.  Notice also that $\size{P} >
\size{P'}$ whenever $|Q| > |\partial R|/2$.
\end{defn}

\begin{defn}[Exponent of a $2$-cell]\label{def:exponent}
Let $X$ be a $2$-complex, and let $R \rightarrow X$ be one of its
$2$-cells.  Let $n$ be the largest number such that the map $\partial
R \rightarrow X$ can be expressed as a path $W^n$ in $X$, where $W$ is
a closed path in $X$.  This number $n$, which measures the periodicity
of the map of $\partial R\rightarrow X$, is the \emph{exponent} of $R$,
and a path such as $W$ is a \emph{period} for $\partial R$.  Notice
that any other closed path which determines the same cycle as $W$ will
also be a period of $\partial R$.  If the exponent $n$ is greater
than $1$, then the $\partial R\rightarrow X$ is called a \emph{proper
power}.
\end{defn}

\begin{defn}[Packet]\label{def:packets}
Let $R$ be a $2$-cell in $X$ of exponent~$n$ and let $W$ be a period
of $\partial R$.  The attaching map $\partial R\rightarrow X$ can be
expressed as a path $W^n \rightarrow X$.  Consider a circle subdivided
into $|W|$ $1$-cells, and attach a copy of $R$ by wrapping $\partial
R$ around the circle $n$~times.  We call the resulting $2$-complex
$\bar R$.  Note that there is a map $\bar R\rightarrow X$ such that
$R\rightarrow X$ factors as $R \rightarrow \bar R \rightarrow X$.
Observe that $\pi_1\bar R \cong \Z/n\Z$ and that the universal cover
of $\bar R$ has a $1$-skeleton which is identical to that of $R$
together with $n$~distinct copies of $R$ attached by embeddings.  The
universal cover of $\bar R$ is the \emph{packet of $R$} and is denoted
by $\widetilde{R}$.  Technically we should write $\widetilde{\bar R}$
but we will use the notation of $\packet{R}$ since $R$ is its own
universal cover and thus there is no danger of confusion. Notice that
if the exponent of $R$ is~$1$ then the packet $\packet{R}$ is the same
as $R$ itself.  Notice also that the map $\packet{R} \rightarrow X$
can be viewed as an extension of the map $R \rightarrow X$.

 Let $\phi\colon Y \rightarrow X$ be a fixed map.  The map $\phi$ will
be called \emph{packed} if whenever there is a lift of a $2$-cell $R
\rightarrow X$ to a $2$-cell $R \rightarrow Y$, there is also a lift
of $\packet{R} \rightarrow X$ to a map $\packet{R} \rightarrow Y$
which extends the map $R \rightarrow Y$.  Since we will treat the
packets $\packet{R}$ as the basic building blocks of our
$2$-complexes, almost all of the maps under discussion will be packed.
\end{defn}

\begin{defn}[$2$-cell attachment]\label{def:attachments}
Let $\phi\colon  Y \rightarrow X$ be an arbitrary packed map and let $R
\rightarrow X$ be a $2$-cell in $X$.  The pair of paths $R \leftarrow
Q \rightarrow Y$ will be called a \emph{$2$-cell attachment site} if
they satisfy the following conditions:
\begin{enumerate}

  \item the path $Q \rightarrow R$ is a subpath of $\partial R$

  \item the diagram \small$\begin{array}{ccc} Q &
        \rightarrow & Y \\ \downarrow & & \downarrow \\ R &
        \rightarrow & X \end{array}$\normalsize commutes

  \item there does not exist a map $R \rightarrow \packet{R}
      \rightarrow Y$ which is a lift of the map $R \rightarrow
      \packet{R} \rightarrow X$ such that the composition $Q
      \rightarrow R \rightarrow \packet{R} \rightarrow Y$ equals the
      path $Q \rightarrow Y$.
\end{enumerate}
\noindent
\end{defn}

Intuitively, a $2$-cell attachment site is a portion of the boundary
of $R$ which is found in the complex $Y$ at a location where the
packet $\packet{R}$ does not already exist.  In other words, it is a
place at which attaching a copy of $\packet{R}$ will have an effect on
the perimeter of the map.  Notice that when the length of $Q$ is equal
to the length of $\partial R$, the path $Q \rightarrow Y$ may have
distinct endpoints even though the endpoints of the path $Q
\rightarrow \partial R$ are identical.

\begin{defn}[Maximal attachment]\label{def:maximal-attachments}
A $2$-cell attachment site is \emph{maximal} if there does not exist
another pair of maps $R \leftarrow Q' \rightarrow Y$ where $Q
\rightarrow R$ is a proper subpath of $Q' \rightarrow R$ and $Q
\rightarrow Y$ is a proper subpath of $Q' \rightarrow Y$.
Technically, we require that there does not exist a proper inclusion
$Q \rightarrow Q'$ such that $Q \rightarrow Q' \rightarrow R$ is the
map $Q \rightarrow R$ and $Q \rightarrow Q' \rightarrow Y$ is the map
$Q \rightarrow Y$.  This forces $Q$ to appear as a proper subpath of
$Q'$ in the same manner in both cases.  If $|Q| < |\partial R|$ and
the $2$-cell attachment site is maximal we will call it an
\emph{incomplete $2$-cell attachment}.  When $|Q| = |\partial R|$, we
will call this a \emph{complete $2$-cell attachment}.  Notice that
complete attachments are automatically maximal.
\end{defn}

\begin{defn}[$2$-cell reduction]
Let $X$ be a weighted $2$-complex and let $Y \rightarrow X$ be a
packed map.  A $2$-cell attachment $R \leftarrow Q \rightarrow Y$ will
be called a \emph{$2$-cell perimeter reduction} if
$\perimeter(\packet{R}) < \perimeter(Q)$.  If $\perimeter(\packet{R})
\leq \perimeter(Q)$ it will be called a \emph{weak $2$-cell perimeter
reduction}.
\end{defn}

\begin{rem}[The main idea]\label{rem:peri-reduction}
A $2$-cell perimeter reduction $R \leftarrow Q \rightarrow Y$ is so
named because it can be used to reduce the perimeter of the map
$\phi\colon Y \rightarrow X$ (Lemma~\ref{lem:attachments}).  The main
idea is as follows: Simply attach the packet $\packet{R}$ to $Y$ along
the path $Q$.  Technically, the identification space $Y \cup_Q
\packet{R}$ is formed by identifying the image of the $1$-cells of $Q$
in $\packet{R}$ with their image in $Y$.  For simplicity, we write
$Y^+ = Y \cup_Q \packet{R}$ for the resulting complex, and we call the
extended map $\phi^+:Y^+ \rightarrow X$.  Since
$\perimeter(\packet{R}) < \perimeter(Q)$, the cells which are in
$\packet{R}$ and not in $Q$ must make a net negative contribution to
the perimeter and consequently $\perimeter(Y^+) < \perimeter(Y)$.  The
details and the qualifications which are necessary to justify this
calculation are contained in Lemma~\ref{lem:complete} and
Lemma~\ref{lem:incomplete}.  These two technical lemmas will be the
key ingredients in the proof of Theorem~\ref{thm:2cellCT}.
\end{rem}

The following lemma shows how the relationship between the perimeter
of $Q$ and the perimeter of the packet $\packet{R}$ can be
reformulated as a relationship between the exponent of $R$, the weight
of $R$, and the perimeter of the complement of $Q$.  This alternative
form makes it easier to verify that a specific reduction is a
perimeter reduction.  The original form is easier to understand
conceptually.

\begin{lem}\label{lem:S-form}
Let $X$ be a weighted $2$-complex, let $R \rightarrow X$ be a
$2$-cell, and let $Q \rightarrow R$ be a subpath of $\partial R$.  If
$n$ is the exponent of $R$ and $S$ is the complement of $Q$ in $R$,
then $\perimeter(\packet{R}) = \perimeter(Q) + \perimeter(S) - n \cdot
\weight(R)$.  Consequently, $\perimeter(\packet{R}) \leq
\perimeter(Q)$ if and only if $\perimeter(S) \leq n \cdot \weight(R)$,
and the first inequality is strict if and only if the second one is
strict.
\end{lem}

\begin{proof}
Since the map $\packet{R} \rightarrow X$ is an immersion,
Lemma~\ref{lem:perimeter-rewritten} can be used to yield the first
equation.  The inequalities then follow as simple rearrangements of
this basic equation.
\end{proof}

We conclude this section with the notion of a redundant $2$-cell.

\begin{defn}[Redundant $2$-cell]
Let $\phi\colon Y \rightarrow X$ be a fixed map and let $R_1
\rightarrow Y$ and $R_2 \rightarrow Y$ be $2$-cells in $Y$.  We say
that $R_2\rightarrow Y$ is {\em redundant (relative to $R_1$ and the
map $Y\rightarrow X$)} provided that $R_1$ and $R_2$ are distinct
$2$-cells in $Y$ which have the same boundary cycle, but $R_1$ and
$R_2$ project to the same $2$-cell in $X$.  More precisely, their
interiors in $Y$ are disjoint, but there exists a map $R_1\rightarrow
R_2$ which restricts to $\partial R_1\rightarrow \partial R_2$, such
that the following two diagrams commute:
\[
\begin{array}{ccc}
\partial R_1&\rightarrow &\partial R_2\\
&\searrow&\downarrow\\
 &&Y\\
\end{array}
\hspace{1in}
\begin{array}{ccc}
R_1&\rightarrow &R_2\\
&\searrow&\downarrow\\
 &&X\\
\end{array}
\]
\end{defn}

Because of the way that perimeter is calculated, redundant $2$-cells
have no effect on the perimeter of $Y \rightarrow X$. This is made
precise below and will be used in the proofs in
Section~\ref{sec:2cell}.

\begin{lem}\label{lem:redundant}
Let $X$ be a weighted $2$-complex, let $Y \rightarrow X$ be a map, and
let $R_1$ and $R_2$ be redundant $2$-cells of $Y$.  If $Y'$ is $Y$
minus the interior of $R_1$ and $\phi'$ is the restriction of $\phi$
to $Y'$, then $\perimeter(Y') = \perimeter(Y)$.  More generally, if
$Y$ and $Y'$ differ by the addition or removal of redundant $2$-cells,
then $\perimeter(Y) = \perimeter(Y')$.
\end{lem}

\begin{proof}
Since the $1$-skeletons are identical and $Y' \subset Y$, it is clear
that each side that is missing at $y$ in $Y$ is also missing at $y$ in
$Y'$.  To see the reverse implication, let $(R,r) \rightarrow (X,x)$
be a side of $X$ which is present at $y$ in $Y$.  If $(R,r)
\rightarrow (X,x)$ lifts to a side of $y$ which is a side of the
$2$-cell $R_1 \rightarrow Y$, then by the definition of redundant
$2$-cells, it also lifts to a side of the $2$-cell $R_2 \rightarrow Y$
at $y$.  Thus every side at $x$ which is present at $y$ in $Y$ is also
present at $y$ in $Y'$.  The final assertion is now immediate.
\end{proof}

Finally, we relate the lack of redundant $2$-cells to immersions
in the following lemma whose proof is immediate.

\begin{lem}\label{lem:no-redundant->immersion}
If $\phi\colon Y \rightarrow X$ is a $1$-immersion and $Y$ has no
redundant $2$-cells, then $\phi$ is an immersion.
\end{lem}

\section{$2$-cell coherence theorem}\label{sec:2cell}
In this section we show how $2$-cell perimeter reductions can be used
to lower the perimeter of a map $Y \rightarrow X$.  At the end of the
section we use this to prove a $2$-cell version of our coherence
theorem.

\begin{lem}[Complete attachment]\label{lem:complete}
Let $X$ be a weighted $2$-complex, let $\phi\colon Y \rightarrow X$ be
a packed $1$-immersion, and suppose that $\perimeter(Y)$ is finite.
If $R \leftarrow Q \rightarrow Y$ is a complete $2$-cell attachment,
then the perimeter of the induced map $\phi^+: Y^+ \rightarrow X$
satisfies the equation
\begin{equation}\label{eq:complete}
  \perimeter(Y^+) \leq \perimeter(Y) - \weight(R) < \perimeter(Y)
\end{equation}
\end{lem}

\begin{proof}
Since by assumption $|Q| = |\partial R|$, the space $Y^+ = Y \cup_Q
\packet{R}$ can be formed by first identifying the endpoints of $Q$
in $Y$, if they are not already identical, and then attaching the
packet $\packet{R}$ along its boundary.

Next, since the space $Y$, with the two endpoints of $Q$ identified,
is a subcomplex of $Y^+$ with an identical $1$-skeleton, any side of
$X$ which is missing at $y$ in $Y^+$ is also missing at $y$ in $Y$.
This shows that the terms in the sum defining $\perimeter(Y^+)$ are
contained as distinct terms in the sum defining $\perimeter(Y)$.

Let $(R,r)$ be a side of $X$, and let $y$ be the image of this $1$-cell
$r$ under the map $Q \rightarrow Y$.  If the side $(R,r)$ was already
present at $y$, then, using the fact that $\phi\colon Y \rightarrow X$ is
a packed $1$-immersion, we find that there already
existed a lift of $\packet{R} \rightarrow X$ to $Y$ for which the
composition $Q \rightarrow R \rightarrow \packet{R} \rightarrow Y$ is
the given map $Q \rightarrow Y$.  Since this contradicts our
assumption that $R \leftarrow Q \rightarrow Y$ is a $2$-cell
attachment site, we have shown that the side $(R,r)$ was missing at
$y$ in $Y$, even though it is clearly present at $y$ in $Y^+$.  If we
repeat this argument for each of the sides of $R$ we can conclude that
$\perimeter(Y^+) \leq \perimeter(Y) - \weight(R)$, which is less than
$\perimeter(Y)$ since $\weight(R) > 0$.
\end{proof}

A careful argument would show that $\perimeter(Y^+) = \perimeter(Y)
-\frac{n}{d}\weight{(R)}$ where $d$ is the exponent of the $2$-cell $R
\rightarrow Y^+$.

\begin{lem}[Incomplete attachment]\label{lem:incomplete}
Let $X$ be a weighted $2$-complex, let $\phi\colon Y \rightarrow X$ be a
packed $1$-immersion and suppose that $\perimeter(Y)$
is finite.  If $R \leftarrow Q \rightarrow Y$ is an incomplete
$2$-cell attachment then the perimeter of the induced map $\phi^+: Y^+
\rightarrow X$ satisfies the equation
\begin{equation}\label{eq:incomplete} \perimeter(Y^+) = \perimeter(Y)
+ \perimeter(\packet{R}) - \perimeter(Q) \end{equation}
\end{lem}

\begin{proof}
Since the perimeter of $\phi\colon Y \rightarrow X$ is unaffected by the
addition or removal of redundant $2$-cells from $Y$
(Lemma~\ref{lem:redundant}), we might as well assume that $Y$ has no
redundancies.  By Lemma~\ref{lem:no-redundant->immersion} this means
that we may assume that $\phi$ is an immersion.  We will now show that
the map $Y^+ \rightarrow X$ is a near-immersion.

Since the maps $Y \rightarrow X$ and $\packet{R} \rightarrow X$ are
immersions we only need to show that this is true when $y$ lies in the
image of $Q$ under the map $Q \rightarrow Y^+$.  Let $(R,r)$ be a side
of $X$ which is present at $y$ in $Y^+$.  If this side was already
present at $y$ in $Y$, then, using the fact that $\phi\colon Y \rightarrow
X$ is a packed immersion, we find that there already existed a lift of
$\packet{R} \rightarrow X$ to $Y$ for which the composition $Q
\rightarrow R \rightarrow \packet{R} \rightarrow Y$ is the map $Q
\rightarrow Y$.  Since this contradicts our assumption that $R
\leftarrow Q \rightarrow Y$ is a $2$-cell attachment, we conclude that
$(R,r)$ must be missing at the $1$-cell $y$ in $Y$.

Next, suppose that $(R,r)$ is a side of $X$ which is present in
$\packet{R}$ at $r_1$ and present in $\packet{R}$ at $r_2$.  Suppose
further that both $r_1$ and $r_2$ lie in $Q$ and that they are sent to
the same $1$-cell $y$ in $Y^+$.  Since all of the sides of $R$ are
distinct, the only way in which this could happen is if the exponent
of $R$ is nontrivial, these two copies of $R$ in $\packet{R}$ are
distinct, and the $1$-cells $r_1$ and $r_2$ differ by a path which is
a multiple of the period $W$ of $\partial R$.  As a consequence we
find that the path from $r_1$ to $r_2$ in $\packet{R}$ is sent to a
closed path in $Y$ which is a multiple of a period of $\partial R$,
and it is possible to extend the path $Q \rightarrow Y$ to the entire
boundary of $R$, thereby contradicting the maximality assumption on
$Q$.  We thus conclude that distinct sides of $2$-cells in
$\packet{R}$ are sent to distinct sides of $2$-cells in $Y^+$.  Since
we also showed that these sides are disjoint from the sides of $X$
which are present at $y$ in $Y$, we now know that the map from $Y^+$
to $X$ is an immersion in a small neighborhood of a point in the
interior of each $1$-cell.

If we assume for the moment that $Y$ is compact, then we can calculate
the perimeter of $Y^+$ using Equation~(\ref{eq:weight-stuff}) of
Lemma~\ref{lem:perimeter-rewritten}.  According to
Equation~(\ref{eq:weight-stuff}), the perimeter of $Y^+$ equals the
weight of its $1$-cells minus the weight of its $2$-cells.  If we
apply Equation~(\ref{eq:weight-stuff}) to $\packet{R}$ and $Y$
separately then we add the weight of their $1$-skeletons and subtract
the weights of their $2$-cells.  The difference between these counts
is precisely the $1$-cells of $Q$ in $\packet{R}$ which get identified
to $1$-cells in $Y$ in the space $Y^+$.  This proves
Equation~(\ref{eq:incomplete}).

In the general case where we assume that $\perimeter(Y)$ is finite but
not that $Y$ is compact, then we cannot use
Equation~(\ref{eq:weight-stuff}) as we did above.  Instead we argue as
follows: Let $S$ be the complement of $Q$ in $\partial R$.  The change
in perimeter from $Y$ to $Y^+$ can be computed by first adding
$\perimeter(S)$ corresponding to the new $1$-cells in $Y^+$ and then
subtracting $n \cdot \weight(R)$ corresponding to the new sides.  The
resulting change in perimeter is $\perimeter(S) - n \cdot \weight(R)$,
which is equal to $\perimeter(\widetilde R) -\perimeter(Q)$ by
Lemma~\ref{lem:S-form}.
\end{proof}

Combining Lemma~\ref{lem:complete} and Lemma~\ref{lem:incomplete}, we
have the following.

\begin{lem}[$2$-cell attachment]\label{lem:attachments}
Let $X$ be a weighted $2$-complex, let $\phi\colon Y \rightarrow X$ be
a packed $1$-immersion and suppose that $\perimeter(Y)$ is finite.  If
$R \leftarrow Q \rightarrow Y$ is a $2$-cell perimeter reduction then
the perimeter of the induced map $\phi^+: Y^+ \rightarrow X$ is
strictly less than $\perimeter(Y)$.  If it is a weak $2$-cell
perimeter reduction, $\perimeter(Y^+) \leq \perimeter(Y)$.
\end{lem}

\begin{proof}
Without loss of generality we may assume that the $2$-cell perimeter
reduction $R \leftarrow Q \rightarrow Y$ is maximal.  Let $Y^+ = Y
\cup_Q \widetilde{R}$.  That $\perimeter(Y^+) < \perimeter(Y)$
[$\perimeter(Y^+) \leq \perimeter(Y)$] now follows immediately from
either Lemma~\ref{lem:complete} or Lemma~\ref{lem:incomplete}
depending on whether the reduction is complete or incomplete.
\end{proof}

In addition to the process of attaching $2$-cells, we will also need a
second operation called folding.

\begin{defn}[Folding along a path]\label{def:folding}
Let $Y\rightarrow X$ be a map between $2$-complexes and let
$P\rightarrow Y$ be a length~$2$ path whose projection to $X$ is
of the form $ee^{-1}$ (i.e. a backtrack).  If the $1$-cells of $P$ are
distinct in $Y$, then the map $Y \rightarrow X$ can be factored as
$Y\rightarrow Y'\rightarrow X$ where the complex $Y'$ is obtained from
$Y$ by identifying the endpoints of $P$ (if they are not already
identical) and then identifying the $1$-cells in the image of $P\rightarrow Y$
in the obvious way.  The complex $Y'$ is said to be
obtained from $Y$ by \emph{folding along the path $P$}.  If $Y$ can be
folded along some path $P\rightarrow Y$, then $Y\rightarrow X$
\emph{admits a fold}.
\end{defn}

\begin{defn}[$2$-cell reduction hypothesis]\label{def:prh}
A space $X$ is said to satisfy the \emph{$2$-cell reduction
hypothesis} if for any map $\phi\colon Y \rightarrow X$ which is a
packed $1$-immersion which is not a $\pi_1$-injection, there exists a
$2$-cell $R\rightarrow X$ and a $2$-cell perimeter reduction $R
\leftarrow Q \rightarrow Y$. Notice that if $X$ satisfies the $2$-cell
reduction hypothesis and $Y \rightarrow X$ is a packed map which does
not admit a fold or a $2$-cell perimeter reduction, then the induced
map $\pi_1Y \rightarrow \pi_1X$ is injective.
\end{defn}

\begin{thm}[$2$-cell coherence]\label{thm:2cellCT}
If $X$ is a weighted $2$-complex that satisfies the $2$-cell reduction
hypothesis, then it satisfies the perimeter reduction hypothesis, and
thus $\pi_1X$ is coherent.
\end{thm}

\begin{proof}
Let $Y \rightarrow X$ be a $1$-immersion which is not
$\pi_1$-injective.  Since adding the $2$-cells necessary to make
$Y\rightarrow X$ a packed map does not increase perimeter, we may
assume it is packed without loss of generality.  By hypothesis, there
is a $2$-cell perimeter reduction and by Lemma~\ref{lem:attachments}
the perimeter of $Y^+$ will be smaller.  The fact that $Y$ and $Y^+$
have the same $\pi_1$ image in $X$ is obvious.  That $\pi_1X$ is
coherent now follows from Theorem~\ref{thm:GCT}.
\end{proof}

\section{Algorithms}\label{sec:algorithms}
The $2$-cell coherence theorem (Theorem~\ref{thm:2cellCT}) can also be
presented as an algorithm for constructing finite presentations from a
given finite set of generators.  The algorithm may be viewed as a
generalization of Stallings' algorithm for graphs \cite{St83}.

\begin{thm}[Algorithm]\label{thm:algorithm}
If $X$ is a compact weighted $2$-complex which satisfies the $2$-cell
reduction hypothesis, then there is an algorithm which produces a
finite presentation for any subgroup of $\pi_1X$ given by a finite set
of generators.
\end{thm}

\begin{proof}
To help clarify that the algorithm terminates, we will use a
complexity function other than the usual perimeter.  We define the
complexity of a map $\phi\colon Y\rightarrow X$ to be the ordered pair
$\big(\perimeter(Y), \size{Y}\big)$ where $\perimeter(Y)$ is the
perimeter of the map, $\size{Y}$ is the number of $1$-cells in $Y$,
and the ordering is the dictionary ordering.
Let $H$ be a subgroup of $\pi_1X$ generated by a set of $r$~elements
represented by closed based paths.  We let $Y_1$ be a based bouquet of
$r$~circles corresponding to these paths, and we define
$\phi_1:Y_1\rightarrow X$ so that $\phi_1$ takes each circle of $Y_1$
to the closed based path that it corresponds to. We then subdivide
$Y_1$ so that $\phi_1$ is combinatorial.  Clearly the image of
$\pi_1Y_1$ equals $H$.  Observe that since $Y_1$ is compact and
$\perimeter(x)$ is finite for each $1$-cell $x$ of $X$, both
$\perimeter(Y_1)$ and $\size{Y_1}$ are finite.  Finally, note that
$Y_1$ is packed.

Beginning with $Y_1$, the algorithm produces a sequence of maps
$\phi_i: Y_i\rightarrow X$ such that for each $i$, $\pi_1(Y_i)$ is
mapped onto $H$.  For each $i$, $Y_{i+1}$ is obtained from $Y_i$ by
either folding along a path in $Y_i$ or by adding a copy of
$\widetilde {R}$ along a path $Q$ in $Y_i$ such that $R \leftarrow
Q\rightarrow Y_i$ is a $2$-cell perimeter reduction.  We will give a
detailed description of these procedures below.  Each of these
procedures will decrease the complexity and so we know that the
sequence must terminate at a $1$-immersion $\phi_t: Y_t \rightarrow X$
such that $Y_t$ does not admit a
$2$-cell perimeter reduction.  Since $X$ satisfies the $2$-cell
reduction hypothesis, we conclude that $\phi_t$ induces a
$\pi_1$-injection, and therefore maps $\pi_1Y_t$ isomorphically onto
$H$, thus yielding a finite presentation for $H$.
As will be seen from the descriptions given below, each of these
procedures can be implemented algorithmically.  Assume inductively
that $Y_i$ is compact and packed and that $\phi_i:Y_i \rightarrow X$
maps $\pi_1Y_i$ onto $H$.

\textbf{ Folding along a path:}
If $\phi_i$ is not an immersion on $Y_i^{(1)}$, then there exists a
map $\rho_i: Y_i \rightarrow Y_{i+1}$ which is obtained by folding
along a path.  There is also a map $\phi_{i+1}:Y_{i+1} \rightarrow X$
such that $\phi_i$ factors as $Y_i\rightarrow Y_{i+1}\rightarrow X$.
Because $\rho_i: Y_i\rightarrow Y_{i+1}$ is $\pi_1$-surjective, we see
that $\phi_{i+1}$ maps $\pi_1Y_{i+1}$ onto $H$.  Thus by
Lemma~\ref{lem:surj} we have $\perimeter(\phi_{i+1}) \leq
\perimeter(\phi_i)$.  Since we also have $\size{Y_{i+1}} <
\size{Y_i}$, we see that the complexity of $\phi_{i+1}$ is strictly
less than the complexity of $\phi_i$.

We can continue folding along paths until we reach a map
$\phi_j:Y_j\rightarrow X$ where the restriction of $\phi_j$ to
$Y_j^{(1)}$ is an immersion.  At this point we begin looking for a
$2$-cell perimeter reduction.

\textbf{ Adding $2$-cells:}
Suppose that the restriction of $\phi_i$ to $Y_i^{(1)}$ is an
immersion, but that there exists a $2$-cell perimeter reduction, $R
\leftarrow Q \rightarrow Y_i$.  The reduction can be chosen to be
maximal, and the result is a $2$-cell reduction which is either
complete or incomplete.  In both cases we define $Y_{i+1}$ to be the
identification space $Y_i \cup_Q \packet{R}$ obtained by identifying
the $1$-cells of $Q$ in $\packet{R}$ with their images in $Y_i$ under
the map $Q\rightarrow Y_i$.  The map $\phi_{i+1}:Y_{i+1}\rightarrow X$
is well defined since $\phi_i:Y_i\rightarrow X$ and
$\packet{R}\rightarrow X$ agree on the respective images of the
$1$-cells of $Q$ which were identified to form $Y_{i+1}$.  Also it is
easy to see that the natural map $\psi_i:Y_i \rightarrow Y_{i+1}$ is
$\pi_1$-surjective, and therefore since $\phi_i = \phi_{i+1} \circ
\psi_i$ we conclude that the $\pi_1$-image of $\phi_{i+1}$ is $H$.  It
is again clear that $Y_{i+1}$ is packed and compact.

Finally, the complexity of $\phi_{i+1}$ is strictly less than the
complexity of $\phi_i$ since by Lemma~\ref{lem:attachments}
$\perimeter(Y_{i+1}) < \perimeter(Y_i)$.
\end{proof}

The sequence of spaces described above is very similar to the sequence
of spaces which would be constructed by Theorem~\ref{thm:GCT} when $X$
satisfies the $2$-cell reduction hypothesis.  The main difference
between the two is that the sequence of spaces in the proof of
Theorem~\ref{thm:GCT} are subcomplexes of the covering space $\widehat
X$.  In an algorithmic approach, the structure of this covering space
is unavailable and the spaces described above have been constructed
without reference to the space $\widehat X$.  In fact, these spaces
may not embed or even immerse into $\widehat X$ throughout the course
of the proof.

We also note the following features of the algorithm:
\begin{enumerate}
\item The algorithm gives an alternate proof of
Theorem~\ref{thm:2cellCT}.
\item The compactness assumption can be replaced by an appropriate
recursiveness hypothesis and the algorithm is still effective.
\end{enumerate}

The algorithm can be used to prove coherence when the hypothesis is
weakened to allow for weak $2$-cell perimeter reductions.  In this
more general context, however, one will not know when to stop running
the algorithm and the algorithm as stated cannot be used effectively,
even when $X$ is compact.  We note, however, that Oliver Payne
\cite{Pa-alg} has developed a variation of our algorithm which is
effective for weak $2$-cell perimeter reductions so long as all of the
sides have positive weights.

We conclude this section with an estimate of the efficiency of the
algorithm for finding finite presentations.

\begin{cor}\label{cor:quadratic-time}
Let $X$ be a compact weighted $2$-complex which satisfies the
perimeter reduction hypothesis and let $Y_1$ be the based bouquet of
$r$ circles corresponding to a set of generators of a subgroup of
$\pi_1X$.  There exist constants $C_1$ and $C_2$ depending only on $X$
such that the algorithm described above terminates in fewer than
$C_1|Y_1|$ steps and the time it takes to complete each step is
bounded by $C_2|Y_1|$, where $|Y_1|$ denotes the number of $1$-cells in
$Y_1$.  In particular, the algorithm to calculate a finite
presentation for the subgroup with these $r$ generators is $O(|Y_1|^2)$.
\end{cor}

\begin{proof}
Since $X$ is compact, there is a bound $C$ on the perimeter of any
$1$-cell in $X$ and a bound $C'$ on the length of the boundary of a
$2$-cell in $X$.  Next, notice that both folds and perimeter reductions
will decrease the integer $C'\perimeter(Y_i) + |Y_i|$.
This is because a fold will decrease the number of $1$-cells without
increasing the perimeter, while a perimeter reduction will decrease
$\perimeter(Y_i)$ by 1 while the number of $1$-cells is increased by at
most $C'-1$.  Since the perimeter remains nonnegative, the
number of steps will be bounded by
\[
        C' \perimeter(Y_1) + |Y_1|
        \leq C'C |Y_1| + |Y_1|
\]

\noindent Thus we can choose $C_1 = C'C +1$.  Notice that since the
number of steps is $O(|Y_1|)$ and since each step adds at most a
bounded number of $1$-cells, the number of $1$-cells in $Y_i$ is also
$O(|Y_1|)$.  And since $Y_i^{(1)}$ is connected, the number of
$0$-cells in $Y_i$ is also $O(|Y_1|)$.

Next we show that the time it takes to complete each step is
$O(|Y_1|)$.  Let $|X|$ be the number of $1$-cells in $X$ and let the
$1$-skeleton of $Y_i$ be represented as an adjacency list.  To check
for the existence of a fold in $Y_i$ only requires an examination of
the links (the adjacency lists) of each $0$-cell.  In each link we
only need to check $2|X|+1$ $1$-cells before we either find a fold or
exhaust the link.  Thus a single link can be checked in constant
time. Since it is well-known that the time it takes to implement a
breadth-first search of a connected graph represented by adjacency
lists is $O(|E|)$ (see \cite[Section~23.2]{CoLeRi90}), the time it
takes to visit each $0$-cell in $Y_i$ is $O(|Y_1|)$, and thus checking
for a fold in $Y_i$ is $O(|Y_1|)$.

Next, suppose that $Y_i^{(1)}$ is immersed into $X$.  Since $X$ is
compact, there is a finite list of paths $Q \rightarrow R \rightarrow
X$ which can lead to $2$-cell perimeter reductions.  Given one of
these paths and a $0$-cell $v$ in $Y_i$ it takes a finite amount of
time to check whether there is a lift of $Q$ which starts at $v$.  The
constant nature of this search depends on the fact that the
$1$-skeleton of $Y_i$ is immersed into $X$.  This guarantees that the
links of the $0$-cells are bounded in size and that at each point there
is at most one extension of the lift which is a viable candidate.
Thus the search for a $2$-cell perimeter reduction in this type of
complex is also $O(|Y_1|)$.  Since the final complex can easily be
converted into a finite presentation in quadratic time, the proof is
complete.
\end{proof}

\section{Path coherence theorem}\label{sec:pathCT}
In this section we provide a second, more technical application of
Theorem~\ref{thm:GCT}. Our new hypothesis will imply that any
immersion which is not $\pi_1$-injective admits a $2$-cell attachment
which does not increase the perimeter, but with additional
restrictions.  The new hypothesis will involve sequences of closed
paths in the space $X$.  We begin with an example which shows why
these technicalities might be desirable.

\begin{figure}\centering
\includegraphics[width= .11\textwidth]{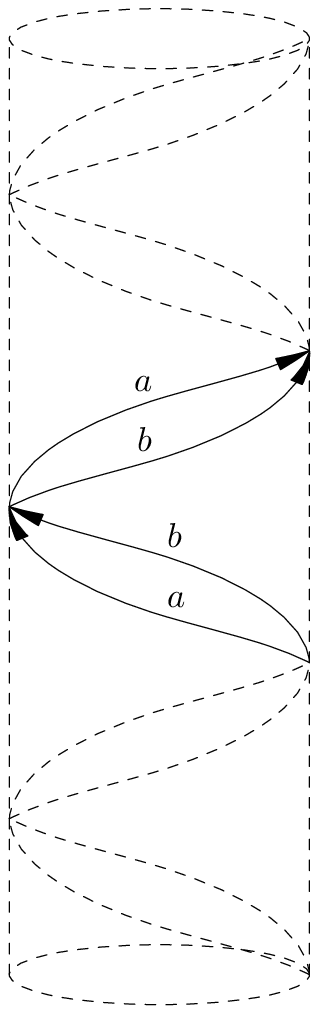}\hspace{1in}
\includegraphics[width= .11\textwidth]{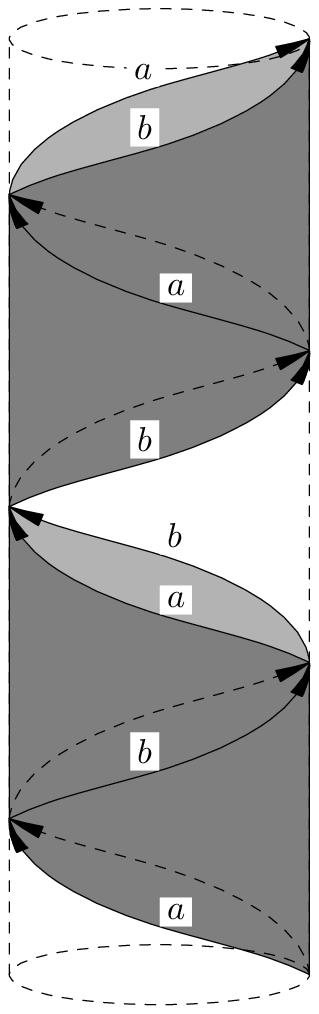}
\caption{The spaces described in Example~\ref{exmp:infinite}
 \label{fig:infinite}}
\end{figure}

\begin{exmp}[Infinite reductions]\label{exmp:infinite}
Let $X$ be the standard $2$-complex of $\langle a,b\mid ab = ba
\rangle \cong \Z^2$, and give $X$ the unit weighting.  Let $H =
\langle ab^{-1} \rangle$, and let $\widehat X$ be the based cover of
$X$ corresponding to $H$.  As illustrated in
Figure~\ref{fig:infinite}, $\widehat X$ is an infinite cylinder.
Observe that every proper $\pi_1$-surjective subcomplex of $\widehat
X$ admits a weak $2$-cell perimeter reduction.  This example will show
that hypothesizing weak $2$-cell perimeter reductions is insufficient
to guarantee that the process of successively attaching $2$-cells will
stop.

Consider the situation where we begin with the subcomplex $Y$ which is
the image of the closed path $a^{-1}b^{-1}ab$ in $\widehat X$.  This
subcomplex is shown on the left side of Figure~\ref{fig:infinite}. The
four $1$-cells determine a length~$4$ closed path which is the
boundary path of a $2$-cell in the cylinder.
Note that the inclusion map $Y \rightarrow \widehat X$ is not
$\pi_1$-injective.  Although there is an obvious complete $2$-cell
attachment which will make the inclusion a $\pi_1$-injection, it is
also possible to apply an infinite sequence of $2$-cell attachments
which are \textit{weak} $2$-cell perimeter reductions, but at each
stage the inclusion map will still fail to be $\pi_1$-injective. These
$2$-cell attachments are formed by adding squares above or below the
square hole bounded by the original closed path.  The right side of
Figure~\ref{fig:infinite} shows the subcomplex obtained by adding two
squares above the original closed path and two squares below it.  The
perimeter is $8$, which is the same as $\perimeter(Y)$.  Clearly, the
operation of adding squares which do not change the perimeter can
continue indefinitely.  We conclude that a weak version of the
$2$-cell reduction hypothesis, in and of itself, is insufficient to
guarantee that a $\pi_1$-injective subcomplex will be obtained after a
finite number of steps.
\end{exmp}

The reason why we never reach a $\pi_1$-injective subcomplex in
Example~\ref{exmp:infinite} is that attached $2$-cells were not
linked in any way to the failure of the $\pi_1$-injectivity.  Our plan
will be that the order in which the $2$-cell attachments are applied
will be tied to the existence of curves which are essential in $Y$ and
null-homotopic in $\widehat X$.  Such precision was not needed for the
$2$-cell reduction hypothesis because the process automatically
stopped after a finite number of steps at a $\pi_1$-injective
subcomplex.

\begin{defn}[Pushing across perimeter-reducing $2$-cell]
Suppose that the path $P'\rightarrow X$ is obtained from $P\rightarrow
X$ by pushing across the $2$-cell $R\rightarrow X$.  According to
Definition~\ref{def:pushing-across}, this means there is a certain
subpath $Q$ of $P$ which is also a subpath of $\partial R$, and $P'$
is obtained from $P$ by replacing $Q$ with its complement $S$ in $R$.
We will now augment this definition with certain perimeter
requirements.  If, in addition, $\perimeter(\packet{R}) <
\perimeter(Q)$, then $P'$ is obtained from $P$ by \emph{pushing across
a perimeter-reducing $2$-cell}.  Similarly, if $\perimeter(\packet{R})
\leq \perimeter(Q)$, then it is obtained by \emph{pushing across a
weakly perimeter-reducing $2$-cell}.
\end{defn}

\begin{defn}[Path reduction hypothesis]\label{def:path-reductions}
We say that a weighted $2$-complex $X$ satisfies the
\emph{path reduction hypothesis} provided the following condition holds:

For every nontrivial closed null-homotopic path $P \rightarrow X$,
there is a sequence of closed paths $\{P_1, \ldots, P_t\}$ which
starts at the path $P = P_1$, ends at trivial path $P_t$, and for each
$i$, $P_{i+1}$ is obtained from $P_i$ by either the removal of a
backtrack or a weakly perimeter-reducing push across a $2$-cell.
\end{defn}

Two elementary conditions which imply the path reduction hypothesis
are a decrease in length and a descrease in area.  In order to make
the second condition precise we recall the defintion of area of a disc
diagram.

\begin{defn}[Disc diagram]\label{def:diagrams}
A {\em disc diagram} $D$ is a compact contractible $2$-complex with a
fixed embedding in the plane.  A {\em boundary cycle} $P$ of $D$ is a
closed path in $\partial D$ which travels entirely around $D$ (in a
manner respecting the planar embedding of $D$).

Let $P\rightarrow X$ be a closed null-homotopic path.  A {\em disc
diagram in $X$ for $P$} is a disc diagram $D$ together with a map
$D\rightarrow X$ such that the closed path $P\rightarrow X$ factors as
$P\rightarrow D\rightarrow X$ where $P\rightarrow D$ is the boundary
cycle of $D$.  The van Kampen's lemma \cite{vK33} essentially states
that every null-homotopic path $P\rightarrow X$ is the boundary cycle
of a disc diagram.  We define $\area(D)$ to be the number of $2$-cells
in $D$.  For a null-homotopic path $P\rightarrow X$, we define
$\area(P)$ to equal the minimal number of $2$-cells in a disc diagram
$D\rightarrow X$ that has boundary cycle $P$.  The disc diagram
$D\rightarrow X$ will then be referred to as a {\em minimal area disc
diagram} for $P$.
\end{defn}

\begin{lem}\label{lem:path-reductions}
Each of the following implies the path reduction hypothesis:
\begin{enumerate}
  \item Every immersed nontrivial null-homotopic path $P\rightarrow X$
  admits a push across a weakly perimeter-reducing $2$-cell which
  yields a strictly shorter path $P' \rightarrow X$.

  \item Every immersed nontrivial null-homotopic path $P\rightarrow X$
  admits a push across a weakly perimeter-reducing $2$-cell which
  yields a path $P' \rightarrow X$ satisfying $\area(P) > \area(P')$.
\end{enumerate}
\end{lem}

\begin{proof}
In either case, there is an obvious procedure for creating the
sequence of paths $P_i \rightarrow X$ which starts at a given closed
null-homotopic path $P \rightarrow X$ and ends at the trivial path.
We first remove backtracks repeatedly until we obtain an immersed
path, then use the condition to find a weakly perimeter-reducing push
across a $2$-cell, and then repeat.  In each case, the process must
terminate at a trivial path after finitely many steps because the
removal of backtracks does not increase either length or area.
\end{proof}

We will now show that the path reduction hypothesis implies coherence.

\begin{thm}[Path coherence]\label{thm:pathCT}
If $X$ is a weighted $2$-complex which satisfies the path reduction
hypothesis, then $X$ satisfies the perimeter reduction hypothesis, and
thus $\pi_1X$ is coherent.
\end{thm}

\begin{proof}
We will assume that the map $Y\rightarrow X$ is packed, for otherwise
we could attach $2$-cells to form a packed map $Y^+ \rightarrow X$
with $\perimeter(Y^+)<\perimeter(Y)$ and with the inclusion $Y
\rightarrow Y^+$ a $\pi_1$-surjection.

If $Y \rightarrow X$ is not $\pi_1$-injective then there is a closed
essential path $P\rightarrow Y$ such that the composition
$P\rightarrow Y\rightarrow X$ is a null-homotopic path in $X$, and by
the path reduction hypothesis there exists a sequence of paths $P_i
\rightarrow X$ for $1\leq i \leq t$ which starts at $P_1 = P
\rightarrow X$, ends at the trivial path $P_t \rightarrow X$, and for
each $i$, $P_{i+1}$ is obtained from $P_i$ by either removing a
backtrack or a weakly perimeter-reducing push across a $2$-cell.  We
will use this sequence of paths to create a sequence of compact spaces
$Y_1 \rightarrow Y_2 \rightarrow \cdots \rightarrow Y_t$ and maps
$Y_i\rightarrow X$ and a sequence of paths $P_i \rightarrow Y_i$ which
are lifts of the paths $P_i \rightarrow X$.

Let $Y_1 = Y$ and let $P_1\rightarrow Y_1$ equal $P\rightarrow Y$ and
assume that $Y_i$ and $P_i \rightarrow Y_i$ have been defined for some
$i$.  The space $Y_{i+1}$ is obtained from $Y_i$ as follows.  If the
operation transforming $P_i$ into $P_{i+1}$ is either the removal of a
backtrack or a weakly perimeter-reducing push across a $2$-cell $R$
where the map $\packet{R} \rightarrow X$ already lifts to $Y$ at the
appropriate point, then $Y_{i+1} = Y_i$.  The exact requirement in the
latter case is that $\packet{R} \rightarrow X$ lift to a map
$\packet{R} \rightarrow Y_i$ such that the composition $Q \rightarrow
R \rightarrow \packet{R} \rightarrow Y_i$ is the map $Q \rightarrow
Y_i$ obtained by restricting the path $P_i \rightarrow Y$.  Since $Y_i
= Y_{i+1}$ it is clear that $\perimeter(Y_i) = \perimeter(Y_{i+1})$.
The path $P_{i+1} \rightarrow Y_{i+1}$ is defined to be the obvious
modification of the path $P_i \rightarrow Y_i$.

\begin{figure}\centering
\includegraphics[width=.33 \textwidth]{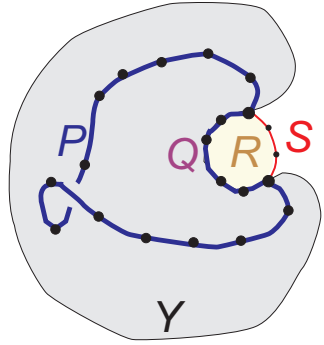}
\caption{A $2$-cell attachment to $Y$ along a portion of the essential
path $P$\label{fig:reductionsite}}
\end{figure}

If the operation is a weakly perimeter-reducing push across a $2$-cell
$R$ and the map $\packet{R} \rightarrow X$ does not lift to $Y$ at the
appropriate point, then $Y_{i+1}$ is defined to be $Y_i \cup_Q
\packet{R}$ and $P_{i+1}$ is again the obvious modification of the
path $P_i \rightarrow Y_i \rightarrow Y_{i+1}$.
Figure~\ref{fig:reductionsite} illustrates a $2$-cell attachment which
arises in this way.  The technical condition is that there does not
exist a lift to a map $\packet{R} \rightarrow Y_i$ such that the
composition $Q \rightarrow R \rightarrow \packet{R} \rightarrow Y_i$
is the map $Q \rightarrow Y_i$ obtained by restricting the path $P_i
\rightarrow Y$.  Notice that in this case $R \leftarrow Q \rightarrow
Y_i$ is a $2$-cell attachment which is a weak $2$-cell perimeter
reduction.  If this attachment is complete, then by
Lemma~\ref{lem:complete}, $\perimeter(Y_{i+1}) < \perimeter(Y_i)$.  If
the attachment is incomplete, then by Lemma~\ref{lem:incomplete},
$\perimeter(Y_{i+1}) = \perimeter(Y_i) + \perimeter(\packet{R}) -
\perimeter(Q) \leq \perimeter(Y_i)$.

In each instance the path $P_{i+1} \rightarrow Y_{i+1}$ is obtained
from $P_i\rightarrow Y_i$ by lifting the operation which occurred in
$X$ to $Y_{i+1}$.  Combining the sequence of perimeter inequalities we
see that $\perimeter(Y) \geq \perimeter(Y_t)$ with a strict inequality
if any of the $2$-cell attachments were complete attachments.  It only
remains to show that at least one of the attachments was complete.
Notice that at each stage the closed path $P_i \rightarrow Y_i \subset
Y_{i+1}$ is homotopic to the closed path $P_{i+1} \rightarrow
Y_{i+1}$.  The crucial observation is that the removal of a backtrack
or an incomplete attachment can never change an essential path into a
null-homotopic one.  Since the final path $P_t$ is null-homotopic, at
least one of the attachments must have been complete.  Thus
$\perimeter(Y) > \perimeter(Y_t)$.

To complete the proof, let $Y^+=Y_t$.  We note $Y_1\rightarrow Y^+$ is
$\pi_1$-surjective because it is the composition $Y=Y_1\rightarrow Y_2
\rightarrow \cdots \rightarrow Y_t$ and for each $i$ the map $Y_i
\rightarrow Y_{i+1}$ is either a homeomorphism, a complete attachment
or an incomplete attachment and thus always a $\pi_1$-surjection.
\end{proof}

\begin{rem}\label{rem:subspace sequence}
Let $X$ be a weighted $2$-complex which satisfies the path reduction
hypothesis.  If $Y$ is a compact connected subcomplex of a cover
$\widehat X$ of $X$, then according to Theorem~\ref{thm:GCT} and
Theorem~\ref{thm:pathCT}, there exists a sequence of subcomplexes $Y =
Y_1 \subset Y_2 \subset \cdots \subset Y_f =Z$ such that for $i\geq
1$, the space $Y_{i+1}$ is the image of $Y_i \cup_Q R$ in $\widehat
X$, where $R \leftarrow Q \rightarrow Y_i$ is a $2$-cell attachment
which is a weak perimeter reduction and $Z$ is a $\pi_1$-isomorphic
subcomplex of $\widehat X$.
\end{rem}

Theorem~\ref{thm:pathCT} and Remark~\ref{rem:subspace sequence} lead
to a pair of interesting corollaries. The first corollary is a bound
on the number of relators needed in the presentation for a finitely
generated subgroup.

\begin{cor}\label{cor:relatorbound}
Let $X$ be a weighted $2$-complex which satisfies the path reduction
hypothesis and assume that each $2$-cell of $X$ is attached along a
simple cycle in $X^{(1)}$.  If $\{ W_i \rightarrow X \}$ is a finite
set of closed based paths in $X$, then the finitely generated subgroup
$H \subset \pi_1X$  generated by the closed based paths
$\{W_i \rightarrow X \}$ has a finite
 presentation with at most $\sum_i \perimeter(W_i)$ relations.
 Similarly, for any $\pi_1$-surjective
subcomplex $Y\subset \widehat X$ whose perimeter is finite and whose
fundamental group is free, there is a finite presentation for $\pi_1
\widehat X$ where the number of relators is at most $\perimeter(Y
\rightarrow X)$.
\end{cor}

\begin{proof}
We first note that the second assertion includes the first assertion
as a special case.  In particular, given closed paths $W_i \rightarrow
X$ we can take the based cover $\widehat{X}$ which corresponds to $H$
and lift the closed paths $W_i \rightarrow X$ to closed paths $W_i
\rightarrow \widehat{X}$ with a common basepoint.  The union of the
images of these closed paths is a finite graph $Y$ in $\widehat{X}$
which satisfies the description in the second half of the corollary.
Thus it suffices to prove the second assertion.

According to Remark~\ref{rem:subspace sequence}, there is a sequence
of subspaces $Y=Y_1\subset Y_2\subset \cdots \subset Y_f=Z$ such that
for $i\geq 1$, the space $Y_{i+1}$ is the image of $Y_i \cup_Q R$ in
$\widehat X$, where $Y_i \leftarrow Q \rightarrow R$ is a $2$-cell
attachment which is a weak perimeter reduction, and such that $Z$ is a
$\pi_1$-isomorphic subcomplex of $\widehat X$.

If all of the $2$-cell attachments are perimeter reductions, then the
argument is easy, because then $f \leq \perimeter(Y)$ and so $Z$ can
be obtained from $Y$ by the addition of fewer than $\perimeter(Y)$
$2$-cells. Consequently, $\pi_1Z$ has a presentation with fewer than
$\perimeter(Y)$ $2$-cells and we are done.

When some of the $2$-cell attachments are weak perimeter reductions,
we argue as follows: Let $\relator(Y_i)$ be the minimum number of
relators which are needed to give a finite presentation of $\pi_1Y_i$.
We will show inductively that $\relator(Y_i) + \perimeter(Y_i) \leq
\perimeter(Y_1)$.  This is true for $i=1$ since the fundamental group
of $Y_1=Y$ is free.  Suppose that $Y_{i+1}$ is obtained from $Y_i$ by
an incomplete attachment. In this case there is a new $1$-cell which
is added to $Y_{i+1}$ and this $1$-cell appears exactly once in the
attaching map of the new $2$-cell.  This is where we use the
additional hypothesis that the attaching map of each $2$-cell embeds
in $X$. Now we can collapse the new $2$-cell across this new $1$-cell
to see that no new relations have been added, although this $2$-cell
attachment may have added new generators.  Since the perimeter has not
increased, the inequality is still true.  If, on the other hand,
$Y_{i+1}$ is obtained from $Y_i$ by a complete attachment, then
$\relator(Y_{i+1}) \leq \relator(Y_i)+1$, but $\perimeter(Y_{i+1})
\leq \perimeter(Y_i)-1$.  Thus the inequality holds in this case.
When the process stops, the perimeter is still nonnegative and thus
$\relator(Y_f)$ is bounded by $\perimeter(Y_1)$, which is the
assertion.
\end{proof}

Note that a similar result (using essentially the same proof) can be
proved under the assumption that each $2$-cell is attached along a
(possibly trivial) power of some simple cycle.  The same type of proof
can also be used to provide an upper bound on the Euler characteristic
of a subgroup.

\begin{cor}\label{cor:Eulerbound}
Let $X$ be a weighted $2$-complex which satisfies the path reduction
hypothesis and assume that no $2$-cell of $X$ is attached by a proper
power.  For any $\pi_1$-surjective compact subcomplex $Y\subset
\widehat X$ there is a compact $\pi_1$-isomorphic subcomplex $Z\subset
\widehat X$ such that $\chi(Z) + \perimeter(Z) \leq \chi(Y) +
\perimeter(Y)$.
\end{cor}

\begin{proof}
According to Remark~\ref{rem:subspace sequence}, there is a sequence
of subspaces $Y=Y_1\subset Y_2\subset \cdots \subset Y_f=Z$ such that
for $i\geq 1$, the space $Y_{i+1}$ is the image of $Y_i \cup_Q R$ in
$\widehat X$, where $Y_i \leftarrow Q \rightarrow R$ is a $2$-cell
attachment which is a weak perimeter reduction, and such that $Z$ is a
$\pi_1$-isomorphic subcomplex of $\widehat X$.

We will deduce that $\chi(Z) + \perimeter(Z) \leq \chi(Y) +
\perimeter(Y)$ by showing that for each $i$ we have $\chi(Y_{i+1}) +
\perimeter(Y_{i+1}) \leq \chi(Y_i) + \perimeter(Y_i)$.  For each $i$,
$Y_{i+1}$ is the union of $Y_i$ and the closure of some $2$-cell.
First suppose that $Y_{i+1}^{(1)} = Y_i^{(1)}$.  In this case,
$Y_{i+1}$ is obtained from $Y_i$ by the addition of a single $2$-cell
and so $\chi(Y_{i+1})= \chi(Y_i)+1$ but $\perimeter(Y_{i+1}) \leq
\perimeter(Y_{i}) -1$, so the inequality holds.  Next suppose that
$Y_{i+1}^{(1)} \neq Y_{i}^{(1)}$, in which case $\chi(Y_{i+1}) \leq
\chi(Y_i)$ because while a $2$-cell has been added, at least one nontrivial arc
of $1$-cells is added to $Y_{i}$ along its endpoints. Since
$\perimeter(Y_{i+1}) \leq \perimeter(Y_i\cup_Q R) \leq
\perimeter(Y_i)$, we see that $\chi(Y_{i+1}) + \perimeter(Y_{i+1})
\leq \chi(Y_i) + \perimeter(Y_i)$.
\end{proof}

We note that a similar statement can be proved in case some of the
$2$-cells are attached by proper powers.  We close the section with
the following problem.

\begin{question}\label{prob:coherence}
Let $X$ be a compact weighted $2$-complex.  Suppose that
$\perimeter(\Delta) < \perimeter(\partial \Delta)$ for every minimal
area disc diagram $\Delta\rightarrow X$.  Does it follow that $\pi_1X$
is coherent?
\end{question}

We conjecture that the answer is yes, but it is not clear how to
proceed.  The problem which arises is that maps $\Delta
\rightarrow X$ which do not send the sides of the boundary of $\Delta$
injectively to the sides of $X$ are, in a fairly strong sense,
unavoidable.

\section{One-relator groups with torsion}\label{sec:1relator}
In this section we present a criterion for the coherence of
one-relator groups with torsion, followed by some applications.
Additional criteria for the coherence of other types of one-relator
groups are developed in \cite{McWi-windmills} and a similar criterion
will be described for small cancellation groups in
Section~\ref{sec:2-cell-small}. The coherence criterion for
one-relator groups is a combination of Theorem~\ref{thm:pathCT} and
the ``spelling theorem'' of B.B.~Newman.  (The original reference is
\cite{Ne68}; see \cite{LySch77} and \cite{HrWi01} for combinatorial
and geometric proofs.)  Here is the theorem as it is usually
formulated.

\begin{thm}[B.B.~Newman]\label{thm:spelling}
Let $G = \langle a_1,\ldots\mid W^n \rangle$ where $W$ is a cyclically
reduced word and $n>1$. Let $U$ and $V$ be words in
$\{a_1^{\pm1}, a_2^{\pm1}, \dotsc\}$ which are equivalent in $G$.  If $U$ is
freely reduced and $V$ omits a generator which occurs in $U$, then $U$
contains a subword $W'$ which is also a subword of $W^n$ and $|W'|>
|W^{n-1}|$.  In particular, if $U$ is a nontrivial word which
represents the identity in $G$, then $U$ contains such a subword $W'$.
\end{thm}

Because of the correspondence between presentations and their standard
$2$-complexes, we will express our main theorem about one-relator
groups in the language of $2$-complexes.  Recall that by
Convention~\ref{conv:combinatorial} the $2$-complexes under
consideration will be those which correspond to presentations whose
defining relators are cyclically reduced.

\begin{thm}[Coherence criterion for one-relator groups]
\label{thm:wt-one-rel}
Let $X$ be a weighted $2$-complex with a unique $2$-cell $R
\rightarrow X$ and a unique $0$-cell.  Let $W\rightarrow X$ be the
period and let $n>1$ be the exponent of $\partial R \rightarrow X$.
If the inequality $\perimeter(S) \leq n \cdot \weight(R)$ holds for
every subpath $S$ of $\partial R$ satisfying $|S| < |W|$, then
$\pi_1X$ is coherent.
\end{thm}

\begin{proof}
Let $P$ be a closed immersed null-homotopic path in $X$ and let $U$ be
the word corresponding to $P$ in the generators of the presentation
corresponding to $X$.  Since $P$ is immersed, $U$ is freely reduced,
and so by Theorem~\ref{thm:spelling} there exists a subpath $Q$ in $P$
such that $Q$ is a subpath of $\partial R$ and $|Q| > (n-1)|W|$.  Note
that we are applying the spelling theorem in the special case where
$V$ is the trivial word.  Since the complement of $Q$ is a path
$S\rightarrow \partial R$ with $|S| < |W|$, we know by assumption that
$\perimeter(S) \leq n \cdot \weight(R)$.  By Lemma~\ref{lem:S-form},
it follows that $\perimeter(\packet{R}) \leq \perimeter(Q)$. Therefore
$P$ can be pushed across a weakly perimeter-reducing $2$-cell.
Moreover, the new path obtained by replacing $Q$ with $S$, is strictly
shorter than $P$ because $n >1$.  Thus by
Lemma~\ref{lem:path-reductions}, $X$ satisfies the path reduction
hypothesis, and so $\pi_1X$ is coherent by Theorem~\ref{thm:pathCT}.
\end{proof}

As an application of Theorem~\ref{thm:wt-one-rel}, we obtain the
following:

\begin{thm}\label{thm:equalweights}
Let $W$ be a cyclically reduced word and let $G = \langle a_1,\ldots
\mid W^n \rangle$.  If $n \geq |W|-1$, then $G$ is coherent.  In
particular, for every word $W$, the group $G = \langle a_1,\ldots \mid
W^n \rangle$ is coherent provided that $n$ is sufficiently large.
\end{thm}

\begin{proof}
Let $X$ be the standard $2$-complex of the presentation with the unit
weighting.  We can assume that $n>1$ since otherwise $G$ is virtually
free, and hence obviously coherent.  Without loss of generality we can
also assume that $W$ is not a proper power since this would only serve
to make the hypothesis more stringent.

Let $R$ denote the unique $2$-cell of $X$, and regard the word $W$ as
a path $W\rightarrow X$.  Then $\partial R \rightarrow X$ has period
$W\rightarrow X$ and has exponent~$n$. Since the perimeter of a
$1$-cell $e$ in $X$ will be the number of times its associated
generator occurs in $\partial R = W^n$ in either orientation, we can
estimate that $\perimeter(e) \leq n \cdot |W|$, and that
$\perimeter(S) \leq n \cdot |W| \cdot (|W| -1)$ for any word $S \in
\partial R$ with $|S| < |W|$.  On the other hand, the weight of the
$2$-cell $R$ is exactly $n \cdot |W|$.  The coherence criterion of
Theorem~\ref{thm:wt-one-rel} will be satisfied so long as
\[
    n \cdot |W| \cdot (|W| -1) \leq n \cdot n \cdot |W|
\]
\noindent
In particular, if $|W|-1\leq n$ then Theorem~\ref{thm:wt-one-rel}
shows that the group $G$ is coherent.
\end{proof}

The next theorem lowers the bound on the exponent by choosing a more
appropriate weight function.

\begin{thm}\label{thm:mingenerator}
Let $W$ be a cyclically reduced word, let $G = \langle a_1,\dotsc, a_r
\mid W^n \rangle$, and let $a_1$ occurs exactly $k$ times $(k > 0)$ in
the word $W$.  If $n \geq k$, then the group $G$ is coherent.  In
particular, if every $a_i$ ($1\leq i\leq r$) occurs in $W$, then $G$
is coherent for all $n \geq \frac{\size{W}}{r}$.
\end{thm}

\begin{proof}
The proof is nearly identical to the previous one, except that the
weight function on $X$ has changed.  Assign a weight of~$1$ to any
side labeled by the generator $a_1$ and assign a weight of~$0$
otherwise.  The perimeter of the $1$-cell labeled by $a_1$ is exactly
$n \cdot k$ (since this is the number of occurrences of $a_1$ in
$W^n$), and the perimeter of any other $1$-cell is $0$.  Since any
word $S \in \partial R$ with $|S| < |R|$ contains at most $k$
$1$-cells labeled by $a_1$, we estimate that $\perimeter(S) \leq n
\cdot k \cdot k$.  On the other hand, $n \cdot \weight(R)$ is exactly
$n \cdot n \cdot k$.  Thus whenever $n \geq k$, the criterion of
Theorem~\ref{thm:wt-one-rel} will be satisfied, and the group $G =
\pi_1X$ will be coherent.
The final assertion is immediate since the word $W$ contains at least
$r$ letters and thus one of them occurs at most
$\frac{\size{W}}{r}$ times.
\end{proof}

\section{Small cancellation I}\label{sec:2-cell-small}
In this section we apply our coherence results to small cancellation
groups.  We begin with a brief review of the basic notions of small cancellation
theory. The reader is referred to \cite{McWi-fanladder} for a rigorous development
of these notions that is consistent with their use here.

\begin{defn}[Piece]\label{def:piece}
Let $X$ be a combinatorial $2$-complex.  Intuitively, a piece of $X$
is a path which is contained in the boundaries of the $2$-cells of $X$
in at least two distinct ways.  More precisely, a nontrivial path
$P\rightarrow X$ is a \emph{piece} of $X$ if there are $2$-cells $R_1$
and $R_2$ such that $P\rightarrow X$ factors as $P \rightarrow R_1
\rightarrow X$ and as $P\rightarrow R_2\rightarrow X$ but there does
not exist a homeomorphism $\partial R_1\rightarrow \partial R_2$ such
that there is a commutative diagram
\begin{equation*}\label{eq:piece}
\begin{array}{ccc}
P             & \rightarrow & \partial R_2\\
\downarrow    & \nearrow    & \downarrow\\
\partial R_1 & \rightarrow & X
\end{array}
\end{equation*}
Excluding commutative diagrams of this form ensures that $P$ occurs in
$\partial R_1$ and $\partial R_2$ in essentially distinct ways.
\end{defn}

\begin{defn}[$C(p)$-$T(q)$-complex]\label{def:sc}
An \emph{arc} in a diagram is a path whose internal vertices
have valence~$2$ and whose initial and terminal vertices have valence~$\geq 3$
The arc is {\em internal} if it its interior lies in the interior of $D$,
and it is a {\em boundary arc} if it lies entirely in $\partial D$.

A $2$-complex $X$ satisfies the $T(q)$ condition if for every
minimal area disc diagram $D\rightarrow X$, each internal $0$-cell of $D$
has valence~$2$ or valence~$\geq q$.
Similarly, $X$ satisfies the $C(p)$ condition if the boundary path of each
$2$-cell in $D$ either contains a nontrivial boundary arc,
 or is the concatenation of at
least~$p$ nontrivial internal arcs.
Finally, for a fixed positive real number $\alpha$,
the complex $X$ satisfies
$C'(\alpha)$ provided that for each $2$-cell
$R\rightarrow X$, and each piece $P\rightarrow X$ which factors as
$P\rightarrow R\rightarrow X$, we have $\size{P}<\alpha\size{\partial
R}$. Note that if $X$ satisfies $C'(\alpha)$ and $n\leq
\frac1\alpha+1$ then $X$ satisfies $C(n)$.

It is a fact that if $D\rightarrow X$ is
minimal area then each nontrivial arc in the interior of $D$ is
a piece in the sense of Definition~\ref{def:piece}.
Although the rough definition given above is
 not quite technically correct (for instance,
it uses minimal area diagrams instead of reduced diagrams), it
should give the reader unfamiliar with small cancellation complexes an
approximate idea of their properties.  We refer the interested reader
to \cite{McWi-fanladder} for precise definitions.
\end{defn}

When $p$ and $q$ are sufficiently large, minimal area diagrams over
$X$ will always contain either spurs or $i$-shells.

\begin{defn}[$i$-shells and spurs]\label{def:i-shells}
Let $D$ be a diagram.  An \emph{$i$-shell} of $D$ is a $2$-cell $R
\hookrightarrow D$ whose boundary cycle $\partial R$ is the concatenation
$P_0P_1\cdots P_i$ where $P_0 \rightarrow D$ is a boundary arc,
the interior of $P_1\cdots P_i$ maps to the interior of $D$,
and  $P_j \rightarrow D$ is a nontrivial interior arc of $D$ for all $j > 0$.
The path $P_0$ is the {\em outer path} of the $i$-shell.

A $1$-cell $e$ in $\partial D$ which is incident with a valence~$1$ $0$-cell $v$
is a \emph{spur}.
In analogy with the outer path of an $i$-shell, we will refer to the
length~$2$ path (either $ee^{-1}$ or $e^{-1}e$)
that passes through $v$ as the {\em outer path} of the spur.
\end{defn}

\begin{figure}\centering
\includegraphics[width=\textwidth]{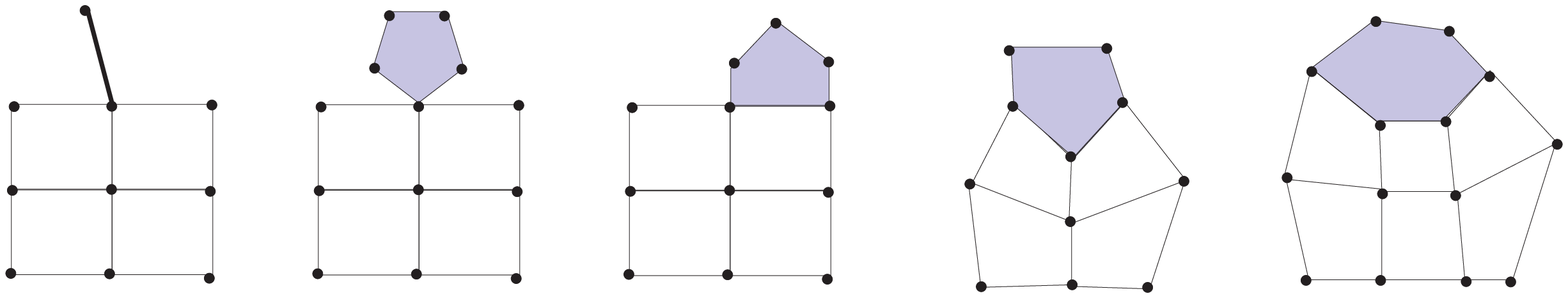}
\caption{Spurs and $i$-shells
\label{fig:outercells}}
\end{figure}

Illustrated from left to right in Figure~\ref{fig:outercells} are disc
diagrams containing a spur, a $0$-shell, a $1$-shell, a $2$-shell, and
a $3$-shell.  In each case, the $2$-cell~$R$ is shaded, and the
boundary arc~$P_0$ is $\partial R \cap \partial D$.

The classical result which forms the basis of small cancellation
theory is called Greendlinger's Lemma (see~\cite[Thm~V.4.5]{LySch77}).
The following strengthening of Greendlinger's Lemma was proven in
\cite[Thm~9.4]{McWi-fanladder}.  While the results of this section
only require Greendlinger's lemma itself, we will require the full
strength of the following theorem in Section~\ref{sec:sc-lqc}.

\begin{thm}\label{thm:fan-classification}
If $D$ is a $C(4)$-$T(4)$ $[C(6)$-$T(3)]$ disc diagram, then one of
the following holds:
\begin{enumerate}
\item $D$ contains at least three spurs and/or $i$-shells with $i\leq
2$ $[i\leq 3]$.
\item $D$ is a ladder of width~$\leq 1$, and hence has a spur, $0$-shell or
$1$-shell at each end.
\item $D$ consists of a single $0$-cell or a single $2$-cell.
\end{enumerate}
Moreover, if $D$ is nontrivial and $v$ is a $0$-cell in $\partial D$,
then $D$ contains a spur or an $i$-shell with $i\leq 2$ $[$$i\leq
3$$]$ which avoids $v$, and if the cut-tree of $D$ has $\ell$ leaves,
then $D$ contains at least $\ell$ separate such spurs and $i$-shells.
\end{thm}

See \cite{McWi-fanladder} for details.  In the present article we will
only need the following immediate corollary.

\begin{cor}\label{cor:fan-classification}
Let $D$ be a $C(4)$-$T(4)$ $[C(6)$-$T(3)]$ disc diagram and let $P$
and $Q$ be immersed paths such that $PQ^{-1}$ is the boundary cycle of
$D$.  If neither path contains the outer path of an $i$-shell in $D$
with $i\leq 2$ $[i\leq 3]$, then every $2$-cell of $D$ contains an
edge in $P$ and an edge in $Q$.
\end{cor}

\begin{thm}[Coherence using $i$-shells]\label{thm:coh-sc-1}
Let $X$ be a weighted $2$-complex which satisfies $C(6)$-$T(3)$
$[C(4)$-$T(4)]$.  Suppose
$\perimeter(S)\leq n\weight(R)$
 for each $2$-cell $R \rightarrow X$
and path $S\rightarrow \partial R$
which is the concatenation
at most three~$[$two$]$ consecutive pieces in the
boundary of $R$.
Then $\pi_1X$ is coherent.
\end{thm}

\begin{proof}
We will prove the $C(6)$-$T(3)$ case; the $C(4)$-$T(4)$ case is
handled similarly.  By Theorem~\ref{thm:pathCT}, it is sufficient to
show that $X$ satisfies the path reduction hypothesis.  Let $P
\rightarrow X$ be a closed immersed nontrivial null-homotopic path.
Let $D\rightarrow X$ be a minimal area disc diagram with boundary
cycle $P$.  According to Theorem~\ref{thm:fan-classification}, there
exists an $i$-shell of $D$ ($i\leq 3$) which avoids the basepoint of
$P$.  By hypothesis, the new boundary path, obtained by removing the
boundary arc and the $2$-cell of this $i$-shell from the diagram, is a
path $P'$ which can be obtained from $P$ by a weakly perimeter-reducing
push across a $2$-cell. Since $P'$ is a path satisfying $\area(P) >
\area(P')$,  Lemma~\ref{lem:path-reductions} shows that $X$ satisfies the path
reduction hypothesis and the proof is complete.
\end{proof}

Theorem~\ref{thm:coh-sc-1} can be improved by using more complicated
weight functions, by using more complicated reductions, or by altering
the presentation substantially before a weight function is applied.
The following example is an illustration of the latter possibility.
Additional examples can be found in Section~\ref{sec:sc-lqc} and in
\cite{McWi-lqc}.

\begin{figure}\centering
\includegraphics[width=.4\textwidth]{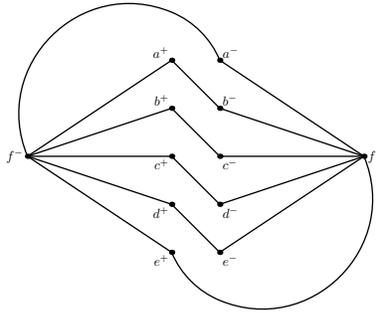}
\caption{The link of $0$-cell of the standard $2$-complex for the
presentation in Example~\ref{exmp:modify} \label{fig:modify}}
\end{figure}

\begin{exmp}\label{exmp:modify}
Consider the following one-relator group.
 \[G = \langle a,b,c,d,e\mid  (abcde) a (abcde) b (abcde) c (abcde) d
(abcde) e \rangle. \]

\noindent Since the relator is not a proper power,
the theorems in Section~\ref{sec:1relator} do not apply.  If we alter
the presentation of $G$ by introducing a new generator $f =
abcde$, then $G = \langle a,b,c,d,e,f\mid abcdef^{-1}, fafbfcfdfe
\rangle$ and the new presentation satisfies certain small
cancellation conditions. This can be seen from the link of the
$0$-cell of the standard $2$-complex $X$ of the modified presentation.
As illustrated in Figure~\ref{fig:modify}, the link is simplicial
and so all pieces are of length~$1$, and since
both relators have length at least~$6$, $X$ is a $C(6)$
presentation. Because the shortest circuit in the link has
length~$4$, $X$ satisfies $T(4)$.  If we assign a weight of $1$ to
each of the sides in the relations which are labeled by $a,b,c,d$
or $e$, and assign a weight of $0$ to the sides labeled $f$, then
the corresponding $1$-cell perimeters and $2$-cell weights are as
follows.  The $1$-cells labeled $a,b,c,d$, and $e$ have a
perimeter of~$2$, and both $2$-cells have a weight of~$5$.  Since the presentation satisfies $C(6)$-$T(4)$ we can use
the coherence criterion for $C(4)$-$T(4)$-complexes
(Theorem~\ref{thm:coh-sc-1}).  The criterion is satisfied since
$\weight(\packet{R}) = 5$ and $\perimeter(Q) \leq 4$ for all
appropriate $R$ and $Q$.  Consequently this group is coherent.
\end{exmp}

The reader may have noticed that although a different weight is
allowed for each side of each $2$-cell of $X$, in all of the
examples we have given so far, we have always chosen the weights
to be equal on all of the sides incident at any particular
$1$-cell in $X$.  Since it is clear that the perimeter of a
$1$-cell $e$ in $X$ is unaffected by the distribution of the
weights among the sides present at $e$ so long as their total is
left invariant, this raises the question of whether the added
flexibility we have allowed will ever be needed.  In our final
example we show that the weights of the sides at $e$ sometimes do
need to be different.

\begin{exmp}
Consider a presentation of the form $\langle a_1,\dotsc \mid U, V
\rangle$. Suppose that for each $i$, the generator $a_i$ appears
exactly the same number of times in $U$ as in $V$, so that in
particular $\size{U} =\size{V}$.
And suppose further that the pieces of $V$ are longer than the
pieces of $U$.  This is the situation in which it makes sense
that a side at $a_i$ in $V$ will need more weight than a side at
$a_i$ in $U$.  The following is a concrete example. Consider the
two-relator presentation:
\begin{displaymath}
\bigg\langle 1,2,3,4,5,6,7,8 \bigg|
\begin{array}{l}
  (1437)(2548)(3651)(4762)(5873)(6184)(7215)(8326), \\
  (1111)(2222)(3333)(4444)(5555)(6666)(7777)(8888)
\end{array}
\bigg\rangle
\end{displaymath}

The parentheses are included for emphasis only.  We will call the
first relator $U$ and the second relator $V$.  Observe that the
presentation is invariant under a cyclic shift of the generators.
Notice also that the presentation satisfies $T(4)$ and $C(16)$,
that every piece in $U$ has length~$1$, and that $V$ has pieces
of length at most~$3$.  Finally it is clear that the subpath
$111222$ in $V$ is a union of two consecutive pieces.

If we assign the sides of $U$ weight~$1$ and we assign the sides of
$V$ weight~$3$ then the perimeter of each $1$-cell is $16$ and the
weights of the $2$-cells corresponding to $U$ and $V$ are $32$ and
$96$ respectively.  Consequently the coherence criterion of
Theorem~\ref{thm:coh-sc-1} is satisfied and so the group is coherent.
On the other hand, if we used the unit weighting, then the perimeter
of each $1$-cell is~$8$.  Observe that the path $111222$ has perimeter
$48$ which is greater than the weight $32$ of the $2$-cell
corresponding to $V$, and so the criterion of
Theorem~\ref{thm:coh-sc-1} fails.

We will now show that more is true.  For this presentation, there
does not exist a way to assign weights to the sides of the
$2$-cells so that (1) all of the sides labeled by a given
generator receive the same weight, and (2) the coherence
criterion of Theorem~\ref{thm:coh-sc-1} is satisfied.  A set of
weights which satisfies the coherence criterion of
Theorem~\ref{thm:coh-sc-1} will be called {\it satisfactory}.
The argument now goes as follows: observe that the sum of any two
sets of weights which are satisfactory will also be satisfactory,
and that a cyclic shift of a set of weights which are
satisfactory will remain satisfactory.  Next suppose that a set
of weights existed which satisfied conditions (1) and (2).  By the
above observations we could add this set of weights to all of its
cyclic shifts to show that a scalar multiple of the unit
perimeter is satisfactory.  But since we know that the unit
perimeter fails the weight criterion, this contradiction shows
that no such set of weights can exist.
\end{exmp}

\section{Fan coherence theorems}\label{sec:fanCT}
In this section we introduce our final coherence hypotheses and our
final coherence theorems which employ fans instead of single
$2$-cells.  Many of the definitions, statements, and proofs will be
analogous to those in previous sections.

\begin{figure}\centering
\includegraphics[width=.4\textwidth]{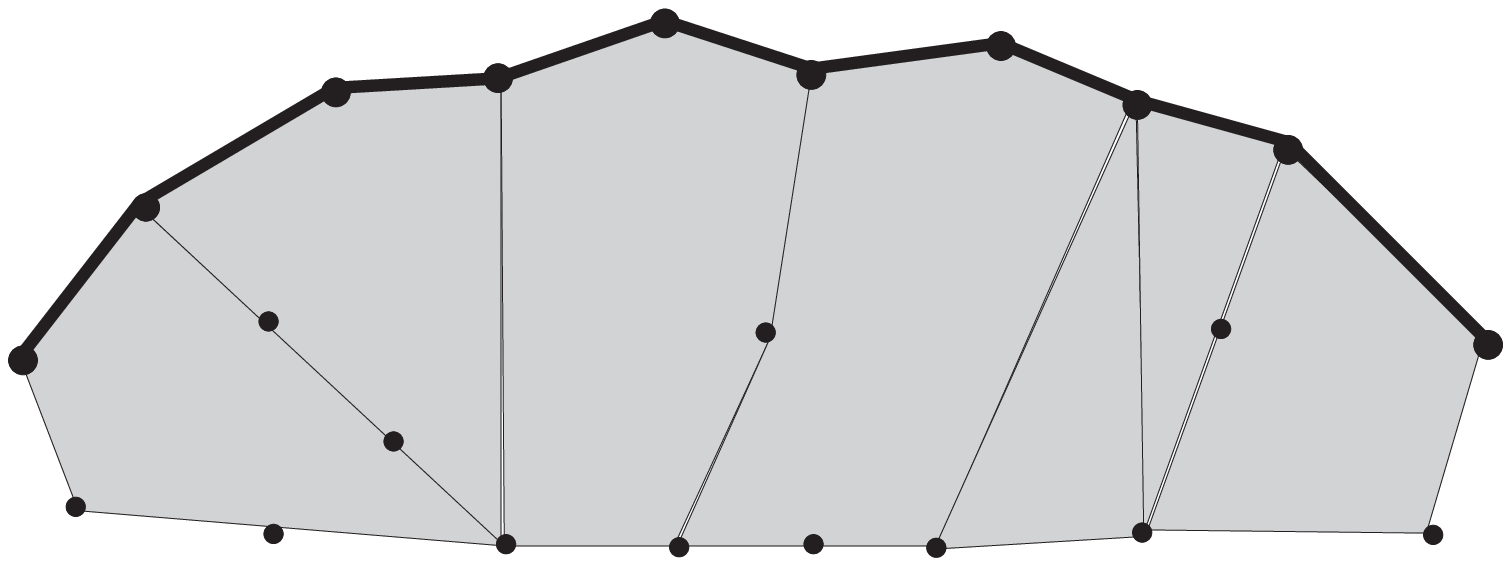}
\hspace{2cm}
\includegraphics[width=.25\textwidth]{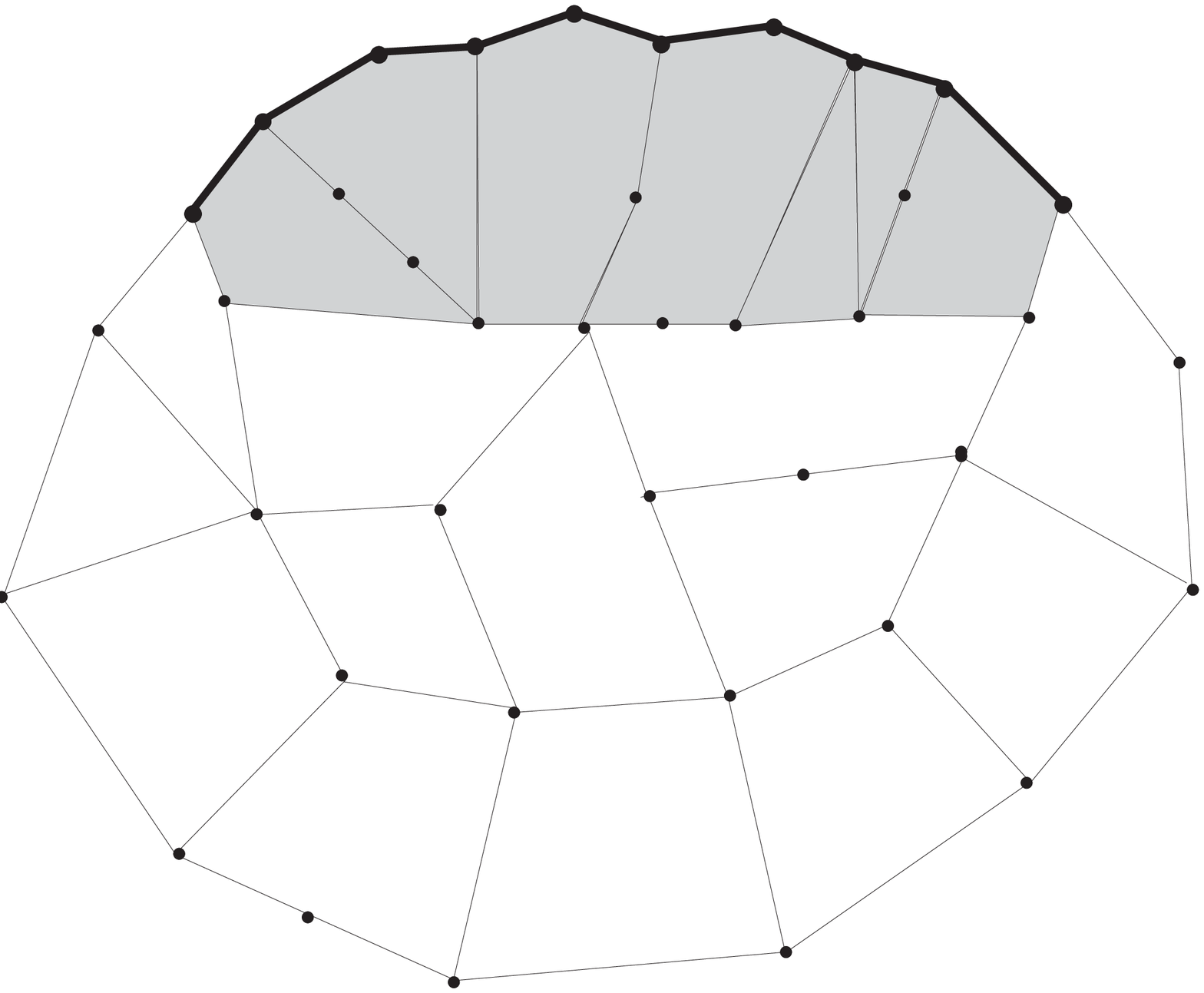}
\caption{\label{fig:typicalfan} On the left is a fan $F$ whose outer
path $Q$ is the bold path on its boundary.  The disc diagram on the
right contains the fan $F$ as a subcomplex.  Note that $Q$ is a
subpath of the boundary path of $D$.}
\end{figure}

\begin{defn}[Fan]
A {\em fan} $F$ is a $2$-complex homeomorphic to a closed disc, which
is the union of closed $2$-cells $\cup_{1\leq i\leq n} R_i$, with the
property that for each $i$, $F-R_i$ is the disjoint union of the
connected sets $\cup_{j < i}R_j$ and $\cup_{j >i}R_j$ (note that when
$i=1$ or $i=n$ one of these sets is empty.)
The {\em outer path} $Q$ of $F$ is a concatenation $Q=Q_1Q_2\dots Q_n$
where each $Q_i$ is a subpath of $\partial R_i$.  We refer the reader
to Figure~\ref{fig:typicalfan} for a picture of a typical fan.  The
unique path $S$ such that $QS^{-1}$ is the boundary cycle of $F$ will
be called the \emph{inner path} of $F$.

Given a map $F\rightarrow X$ there is a unique extension to a packed
map (Definition~\ref{def:packets}) where the $1$-skeleton of the
domain is unchanged.  We will denote this extended domain by $\widetilde F$
in analogy with $\widetilde R$.

We will only be interested in fans equipped with a map $F\rightarrow
X$ such that $\widetilde F \rightarrow X$ is a near-immersion, and we
will refer to such a mapped fan as a {\em fan in $X$}.  In this case,
we will also regard the outer path $Q\rightarrow F$ of $F$ as a path
$Q\rightarrow X$.

The disc diagram $D\rightarrow X$ {\em contains the fan $F\rightarrow
X$}, provided that $F\rightarrow X$ factors as $F\rightarrow
D\rightarrow X$, where the outer path $Q$ of $F$ maps to $\partial D$,
and the inner path $S$ of $F$ is an internal path in $D$.
\end{defn}

\begin{exmp}
The simplest fans are $i$-shells (Definition~\ref{def:i-shells}).  In
this case $\partial R$ is the concatenation $QS^{-1}$ where $Q$ is the
outer path of the $i$-shell, and $S$ is the concatenation of
$i$-pieces in $X$.
\end{exmp}

We will often be interested in a collection $\fantypes$ of fans in a
$2$-complex $X$ which satisfy additional properties.  The next three
definitions are technical conditions which will enable us to perform a
perimeter calculation for fan attachments parallel to the calculation
in Lemma~\ref{lem:attachments}.

\begin{defn}[Perimeter-reducing fan]\label{def:reducing-fan}
Let $X$ be a weighted $2$-complex.  The fan $F\rightarrow X$ is
\emph{perimeter-reducing} provided that the perimeter of $\widetilde
F\rightarrow X$ is less than the perimeter of its outer path
$Q\rightarrow X$.  In other words, $\perimeter(\widetilde F) <
\perimeter(Q)$. Similarly, it is \emph{weakly perimeter-reducing}
if $\perimeter(\widetilde F) \leq \perimeter(Q)$.
\end{defn}

\begin{defn}[Missing along outer path]
Let $Y\rightarrow X$ be a packed $1$-immersion,
let $F\rightarrow X$ be a fan, and let $Q\rightarrow Y$ be a lift of
the outer path of $F$ to $Y$.  We say $F$ is \emph{missing in $Y$
along $Q\rightarrow Y$} provided that the lift of $Q\rightarrow
\widetilde F\rightarrow X$ to the path $Q\rightarrow Y$ does not
extend to a lift of any $2$-cell $R$ of $\widetilde F$ that intersects
$Q$ in a nontrivial path.  Specifically, if $R$ is a $2$-cell of
$\widetilde F$ and $Q'=\partial R \cap Q$ is a nontrivial path, then
there should not exist a lift of $R \cup_{Q'} Q$ to $Y$ which extends
the path $Q \rightarrow Y$.  Equivalently $F\rightarrow X$ is missing
in $Y$ along $Q\rightarrow Y$ provided that for each $1$-cell $q$ in
$Q$, no side of a $2$-cell of $X$ at $x$ is present at both $y$ and
$f$, where $f$, $y$, and $x$ are the images of $q$ in $\widetilde F$,
$Y$, and $X$.
\end{defn}

\begin{defn}[Spread-out]\label{def:spread}
A fan $F\rightarrow X$ is \emph{spread-out} provided that the sides of
$2$-cells of $\widetilde F$ along $1$-cells in the outer path
$Q\rightarrow F$ project to distinct sides of $2$-cells along
$1$-cells in $X$.  This condition is certainly satisfied when the
outer path $Q\rightarrow \widetilde F$ projects to a path $Q
\rightarrow X$ which does not pass through any $1$-cell of $X$ more
than once. For instance, $F\rightarrow X$ is spread-out when
$\widetilde F\rightarrow X$ is an embedding, and it is spread-out when
$Q\rightarrow X$ is a (possibly closed) simple path.
\end{defn}

The following lemma calculates the perimeter of $Y \cup_Q \widetilde
F$ in terms of the perimeters of its constituents.

\begin{lem}[Fan attachment]\label{lem:fan-attachments}
Let $X$ be a weighted $2$-complex, let $\phi\colon Y \rightarrow X$ be
a packed $1$-immersion with $\perimeter(Y) < \infty$, and let
$Q\rightarrow Y$ be a lift of the outer path of a fan $F\rightarrow
X$.  If $F\rightarrow X$ is spread-out and $F$ is missing along
$Q\rightarrow Y$, then, letting $Y^+ = Y \cup_Q \widetilde F$, the
perimeter of the induced map $\phi^+: Y^+ \rightarrow X$ satisfies:
\begin{equation}\label{eq:fan}
\perimeter(Y^+) = \perimeter(Y) + \perimeter(\widetilde F) -
\perimeter(Q)
\end{equation}
\noindent
Thus, if $ F\rightarrow X$ is perimeter-reducing then $\perimeter(Y^+)
< \perimeter(Y)$ and if $F\rightarrow X$ is weakly perimeter-reducing then
$\perimeter(Y^+) \leq \perimeter(Y)$.
\end{lem}

\begin{proof}
The proof is similar to that of Lemma~\ref{lem:incomplete} where it is
obvious that a fan consisting of a single $2$-cell is spread-out.
Since the perimeter of $\phi\colon Y \rightarrow X$ is unaffected by
the addition or removal of redundant $2$-cells from $Y$
(Lemma~\ref{lem:redundant}), we may assume that $Y$ has no
redundancies.  By Lemma~\ref{lem:no-redundant->immersion} this means
that we may assume that $\phi$ is an immersion.

Next, we show that the map $Y^+ \rightarrow X$ is a near-immersion.
By the definition of a fan, distinct $1$-cells of $Q$ are sent to
distinct $1$-cells in $\widetilde F$ under the map $Q \rightarrow
\widetilde F$ and thus distinct $1$-cells of $Y$ are mapped to
distinct $1$-cells of $Y^+$ under the map $Y \rightarrow Y^+=Y \cup_Q
\widetilde F$.  This shows that the induced map $\sides_Y\rightarrow
\sides_{Y^+}$ is an injection.  On the other hand, the map $\widetilde
F \rightarrow X$ is a near-immersion by definition.  Combined with the
fact that $F\rightarrow X$ is spread-out, this shows that the induced
map $\sides_{\widetilde F}\rightarrow \sides_{Y^+}$ is also an
injection.  Thus, if $Y^+\rightarrow X$ fails to be a near-immersion,
it must fail along the path $Q \rightarrow Y^+$.  More precisely, to
show that $Y^+\rightarrow X$ is a near-immersion, it only remains to
be shown that a side of $X$ which lifts to a side of $Q$ in
$\widetilde F$ could not also lift to a side of $Q$ in $Y$.  This is
impossible because of the assumption that $F$ is missing along
$Q\rightarrow Y$.

If we assume for the moment that $Y$ is compact, then we can calculate
the perimeter of $Y^+$ using Equation~(\ref{eq:weight-stuff}) of
Lemma~\ref{lem:perimeter-rewritten}.  According to
Equation~(\ref{eq:weight-stuff}), the perimeter of $Y^+$ equals the
sum of the perimeters of its $1$-cells minus the weights of its
$2$-cells.  If we apply Equation~(\ref{eq:weight-stuff}) to
$\widetilde F$ and $Y$ separately then we would add the perimeters of
their $1$-skeletons and subtract the weights of their $2$-cells.  The
difference between these counts arises from the $1$-cells of $Q$ in
$\widetilde F$ which get identified with $1$-cells of $Y$ in the space
$Y^+$.  This proves Equation~(\ref{eq:fan}).  In the general case
where $\perimeter(Y)$ is finite but $Y$ is not compact, the proof
proceeds as in Lemma~\ref{lem:incomplete}, except that
Lemma~\ref{lem:perimeter-rewritten} is applied to the packed fan
$\widetilde F$ instead of the packet $\widetilde R$.
\end{proof}

Having established conditions under which we can control the change in
perimeter, it is now relatively easy to define a hypothesis and prove
a coherence theorem.

\begin{defn}[Fan reduction hypothesis]\label{def:fanRH}
A packed $1$-immersion $\phi\colon Y \rightarrow X$ admits a {\em fan
perimeter reduction} provided there is a perimeter-reducing spread-out
fan $F \rightarrow X$ and there exists a lift of its outer path to $Y$
such that $F$ is missing along $Q \rightarrow Y$.  A weighted
$2$-complex $X$ satisfies the \emph{fan reduction hypothesis} if each
packed $1$-immersion $\phi\colon Y \rightarrow X$ which is not
$\pi_1$-injective, admits a fan perimeter reduction.
\end{defn}

\begin{thm}[Fan coherence]\label{thm:fanCT}
Let $X$ be a weighted $2$-complex. If $X$ satisfies the fan reduction
hypothesis then $X$ satisfies the perimeter reduction hypothesis, and
thus $\pi_1X$ is coherent.
\end{thm}

\begin{proof}
Let $Y \rightarrow X$ be a $1$-immersion which is not
$\pi_1$-injective.  Since adding the $2$-cells necessary to make
$Y\rightarrow X$ a packed map does not increase perimeter, we may
assume it is packed without loss of generality.  By hypothesis, there
is a perimeter-reducing fan $F\rightarrow X$ which is spread-out and a
lift of its outer path to $Y$ such that $F$ is missing along
$Q\rightarrow Y$.  This can be used to create a complex $Y^+$ whose
perimeter is smaller, by Lemma~\ref{lem:fan-attachments}.  The fact
that $Y$ and $Y^+$ have the same $\pi_1$ image in $X$ is obvious.
Thus $X$ satisfies the perimeter reduction hypothesis.
\end{proof}

In most applications, we will only use a special case of
Theorem~\ref{thm:fanCT} which can be formulated in terms of disc
diagrams.

\begin{thm}[Diagram fan coherence]\label{thm:fan-diagramCT}
Let $X$ be a weighted $2$-complex and let $\fantypes$ be a collection
of perimeter-reducing spread-out fans in $X$.  If each nontrivial
 minimal area disc diagram $D\rightarrow X$ contains a spur or a fan
from $\fantypes$, then
$X$ satisfies the perimeter reduction hypothesis, and thus $\pi_1X$ is
coherent.
\end{thm}

\begin{proof}
By Theorem~\ref{thm:fanCT}, it is sufficient to show that $X$
satisfies the fan reduction hypothesis.  Let $Y$ be a compact
$2$-complex and let $Y \rightarrow X$ be a packed $1$-immersion which
is not $\pi_1$-injective.  There exists at least one essential
immersed closed path $P \rightarrow Y$ whose image in $X$ is a
null-homotopic closed immersed path.  We assume that $P \rightarrow Y$
has been chosen so that $\area(P)$ is as small as possible, and such
that the disc diagram $D\rightarrow X$ realizes this minimum area.  By
hypothesis, $D$ contains a fan $F\rightarrow D\rightarrow X$
which is perimeter-reducing and spread-out.  Let
$Q\rightarrow Y$ be the restriction of $P\rightarrow Y$ to the outer
path of $F$.  Since $P$ was chosen to have minimal area, the fan $F$
is missing along $Q \rightarrow Y$.  Consequently $X$ satisfies the
fan reduction hypothesis.  Indeed, if some some side of a $2$-cell $R$
in $\widetilde F$ is already present along some edge $q$ of the path
$P$, then $P$ is homotopic in $Y$ to a path $P'$ which travels around
the boundary of $D-\big(\text{Interior}(R)\cup q\big)$.  But
$\area(P')<\area(P)$, so we can find an immersed essential path in $Y$
which bounds a smaller area diagram in $X$, and this is impossible.
\end{proof}

We conclude this section with two further generalizations of results
from the previous sections.

\begin{thm}[Fan algorithm]\label{thm:fan algorithm}
If $X$ is a compact weighted $2$-complex which satisfies the fan
reduction hypothesis, then there is an algorithm which produces a
finite presentation for any subgroup of $\pi_1X$ given by a finite set
of generators.
\end{thm}

\begin{proof}
The proof is analogous to the proof of Theorem~\ref{thm:algorithm} and
we leave the details to the reader.
\end{proof}

\begin{thm}[Weak fan coherence]
Let $X$ be a compact weighted $2$-complex.  Let $\fantypes$ be a
collection of spread-out weakly perimeter-reducing fans, and suppose
that for each fan $F\in \fantypes$, we have $\perimeter(\packet{F}) <
\perimeter(\partial F)$.  If every nontrivial minimal area disc
diagram contains a spur or a fan in $\fantypes$ then $\pi_1X$ is
coherent.
\end{thm}

\begin{proof}
The proof is essentially a generalization of the proof of
Theorem~\ref{thm:pathCT} that uses fans instead of $2$-cells.  Let
$Y_1$ be a compact $\pi_1$-surjective packed subcomplex of a cover
$\widehat X$.  If $Y_1$ is not $\pi_1$-injective, then there is a minimal
area disc diagram $D\rightarrow \widehat X$ whose boundary cycle is an
essential immersed path in $Y_1$.
As in the proof of Theorem~\ref{thm:fan-diagramCT}, a minimal area
disc diagram for an essential immersed path in $Y$ yields a sequence
of weakly perimeter reducing spread-out fans that can be attached.

We claim that in the appropriate sense $\packet F$ is missing along
the $1$-cells in $\partial F$ that map to $1$-cells of $Y_i$ in $\widehat
X$.  Indeed, if some $2$-cell of $\packet{F}$ was already contained in
$Y_i$ then a corresponding $2$-cell $R$ of $F$ is contained in $Y_i$.
Let $\partial R$ be the concatenation $Q_1Q_2^{-1}$ where $Q_1$ is
the part of $\partial R$ that is the subpath of the outer path $Q$ of
$F$.
Now $P_i$ is homotopic in $D_i$ and $Y_i$ to a path $P_i'$ which is
identical to $P$ except that $Q_2$ is substituted for $Q_1$.  Since
$P_i'$ doesn't go around $R\subset D_i$, we see that $\area(P_i')\leq
\area(P_i)-1$ and therefore after removing spurs from a disc diagram
for $P_i'$ (and identifying some $1$-cells on the boundary), we obtain
an immersed path homotopic to $P_i$ in $Y_i$ whose area is strictly
less than the area of $P_i$ which is impossible.

We let $Y_{i+1}$ be the union of $Y_i$ with (the image of) $\packet{F}$
be a new compact subcomplex in $\widehat X$.  Now the outer path $Q$ of
$F$ extends to a path $\partial F\rightarrow \widehat X$, and the argument
breaks down according to whether $\partial F$ is contained in $Y_i$.
If $\partial F \subset Y_i$ then our hypothesis that
$\perimeter(\packet{F}) < \perimeter(\partial F)$ implies that
$\perimeter(Y_{i+1}) < \perimeter(Y_i)$.  If $\partial F$ is not a
path in $Y_i$, then $P$ is not null-homotopic in $Y_{i+1}$ since
$\pi_1(Y_i\cup\packet{F})= \pi_1(Y_i\cup F)$ and $Y_i\cup F$ collapses
onto the union of $Y_i$ and some nontrivial arcs.  Now our hypothesis
that $F$ is weakly perimeter reducing implies that
$\perimeter(Y_{i+1})\leq \perimeter(Y_i)$, and $Y_{i+1}$ contains the
essential immersed path $P_{i+1}$ with $\area(P_{i+1})<\area(P_i)$,
where $P_{i+1}$ is defined as follows: First remove the interiors of
$F$ and $Q$ from $D_i$ to obtain a diagram $D_i'$, and then fold
$\partial D_i'$ until the boundary is immersed.  Note that we cannot
obtain a sphere in this way, because otherwise $\partial F$ is the
same as $\partial D_i$ and so we could have used $F$ instead of $D_i$
to begin with, contradicting that $D_i$ is minimal area.

This process can only be repeated finitely many times without the perimeter
strictly decreasing and hence $\widehat X$ satisfies the perimeter reduction
hypothesis and so $\pi_1X$ is coherent.
\end{proof}

\section{Quasi-isometries and quasiconvexity}\label{sec:quasi}
In this section we review the interconnections between
quasi-isometries, quasiconvexity, and word-hyperbolicity.  Since these
results are well-known, we simply state the definitions and lemmas we
will need and refer the interested reader to \cite{ABC90},
\cite{Ep92}, and \cite{Sh91} for more detailed accounts.

\begin{defn}[Geodesic metric space]\label{def:geo}
Let $(X,d)$ and $(X',d')$ be metric spaces.  A map $\phi\colon X'
\rightarrow X$ which preserves distances is called an \emph{isometric
embedding} of $X'$ into $X$, and an isometric embedding of an interval
$[a,b]$ of the real line is called a \emph{geodesic} from $\phi(a)$ to
$\phi(b)$.  If any two points in $X$ can be connected by a geodesic,
then $X$ is  a \emph{geodesic metric space}.
\end{defn}

A fundamental example of a geodesic metric space is a connected graph
with the path metric.  Note that by a `graph' we mean a
$1$-dimensional CW-complex, so that loops and multiple edges are
allowed.

\begin{defn}[Path metric]\label{def:graphs}
The \emph{path metric} on a connected graph $\Gamma$ makes each
$1$-cell of $\Gamma$ locally isometric to the unit interval, and then
defines the distance between two arbitrary points of $\Gamma$ to be
the length of the shortest path between them.  It is easy to see that
such a minimal path always exists, and that it will be a
geodesic. Thus connected graphs are geodesic metric spaces using the
path metric.
\end{defn}

\begin{defn}[Cayley graph]\label{def:cayley}
Let $X$ be a connected $2$-complex, and let $\widetilde X$ be its
universal cover.  Since the $1$-skeleton $\widetilde X^{(1)}$ is a
connected graph it is a geodesic metric space with the path
metric. If $X$ has a unique $0$-cell, then $X$ is the standard
$2$-complex of some group presentation $G = \langle A | \script{R}
\rangle$, and the graph $\widetilde X^{(1)}$ is the \emph{Cayley graph
of the presentation}.  Alternatively, the Cayley graph, often denoted
$\Gamma(G,A)$, can be defined as follows.  Begin with a $0$-cell set
corresponding to the elements of $G$ and an edge set labeled by the
elements of $G \times A$.  Then attach the edges to the $0$-cells so
that the edge labeled $(g,a)$ begins at the $0$-cell $g$ and ends at
the $0$-cell $ga$.  We endow $\Gamma(G,A)$ with the path metric.
Since $G$ can be identified with the $0$-skeleton of $\Gamma(G,A)$, we
can metrize $G$ by giving it the subspace metric. Since this metric on
$G$ depends on the generating set $A$, we will denote the resulting
metric space by $G_A$.
\end{defn}

Although distinct generating sets for $G$ will produce distinct
metrics using this procedure, all of the metrics on a finitely
generated group will be roughly equivalent.  We will now make this
precise.

\begin{defn}[Quasi-isometry]\label{def:qgeo}
Let $(X,d)$ and $(X',d')$ be metric spaces and let $\phi\colon X'
\rightarrow X$ be a map between them.  If there exist constants $K >
0$ and $\epsilon \geq 0$ such that for all $x,y \in X'$,
\[
K d(x,y) + \epsilon > d(\phi(x),\phi(y)) > \frac{1}{K} d(x,y) -
\epsilon
\]

\noindent then $\phi$ is a \emph{$(K,\epsilon)$-quasi-isometric
embedding} of $X'$ into $X$.  The special case of a
$(K,\epsilon)$-quasi-isometric embedding of an interval of the real
line into $X$ is  a
\emph{$(K,\epsilon)$-quasigeodesic}.  If every point in $X$ is within
a uniformly bounded distance of a point in the image of $\phi$, then
$\phi$ is a \emph{$(K,\epsilon)$-quasi-isometry} between $X'$
and $X$.  A map will be called a \emph{quasi-isometry} if it is a
$(K,\epsilon)$-quasi-isometry for some choice of $K$ and $\epsilon$,
and the spaces involved will be said to be quasi-isometric.  The
notion of quasi-isometry is an equivalence relation on spaces in the
following sense.  If there is a quasi-isometry from $X$ to $Y$ then
there also exists a quasi-isometry from $Y$ to $X$, and if $\rho: X
\rightarrow Y$ and $\phi\colon Y \rightarrow Z$ are quasi-isometries,
then the composition $\phi \circ \rho$ is also a quasi-isometry.
\end{defn}

If a finitely generated group $G$ acts in a reasonably nice way on a
reasonably nice space, then the group, using the metric derived from
its Cayley graph, will be quasi-isometric to the space it acts on. The
following theorem contains the precise statement of this fact.

\begin{lem}[Theorem~3.3.6 of \cite{Ep92}]\label{lem:nice-action}
Let $X$ be a locally compact, connected, geodesic metric space.  Let
$G$ be a finitely generated group which acts on $X$ properly
discontinuously and cocompactly by isometries.  Then for any point
$x\in X$, and for any finite set of generators $A\subset G$, the map
$G \rightarrow X$ defined by $g\rightarrow gx$ is a quasi-isometry,
where we give $G$ its Cayley graph metric relative to the generating
set $A$.
\end{lem}

As a corollary, we see that changing generating sets induces a
quasi-isometry.

\begin{cor}\label{cor:qi}
If $A$ and $B$ are finite generating sets for a group $G$, then the
metric spaces $G_A$, $G_B$, $\Gamma(G,A)$, and $\Gamma(G,B)$ are all
quasi-isometric. Furthermore, the quasi-isometry $G_A\rightarrow G_B$
is induced by the identity map $G\rightarrow G$.
\end{cor}

A second fundamental notion is that of a quasiconvex subspace of a
metric space.

\begin{defn}[Quasiconvexity]\label{def:qcon}
A subspace $Y$ of a geodesic metric space $X$ is
\emph{$K$-quasiconvex} if there is a $K$-neighborhood of $Y$ which
contains all of the geodesics of $X$ that begin and end in $Y$.  A
subspace is $\emph{quasiconvex}$ provided that it is
$K$-quasiconvex for some $K$.  The notion of quasiconvexity can be
extended to groups and subgroups via Cayley graphs.  Specifically, a
subgroup $H$ of a group $G$ generated by $A$ is
\emph{quasiconvex} if the $0$-cells corresponding to $H$ form a
quasiconvex subspace of $\Gamma(G,A)$. The group $G$ generated by $A$
is \emph{locally quasiconvex} if every finitely generated
subgroup is quasiconvex.
\end{defn}

We record the following two properties of quasiconvex subgroups.  See
\cite{Sh91} and the references therein for details.

\begin{lem}[Proposition~1 of \cite{Sh91}]\label{lem:short}
If $H$ is a quasiconvex subgroup of a group $G$ generated by finite
set $A$, then $H$ itself is finitely generated.
\end{lem}

\begin{lem}\label{lem:qc->qie}
Let $G$ be a group with finite generating set $A$ and let $H$ be a
subgroup of $G$ with finite generating set $B$.  If $H$ is a
quasiconvex subspace of $\Gamma(G,A)$, then $H_B$ is
quasi-isometrically embedded in $\Gamma(G,A)$.
\end{lem}

Although the various metrics which have been defined for a group $G$
are all equivalent up to quasi-isometry (Corollary~\ref{cor:qi}), the
generating set $A$ does need to be specified in
Definition~\ref{def:qcon}.  This is because the notion of
quasiconvexity is not well-behaved under quasi-isometries.  In
particular, the group $\Z\times \Z$ shows that the converse of
Lemma~\ref{lem:qc->qie} is false.  The dependence of
quasiconvexity on generating sets and the distinction between
quasiconvex subgroups and quasi-isometrically embedded subgroups
disappears once we restrict our attention to word-hyperbolic groups.

\begin{defn}[Hyperbolic spaces and groups]\label{def:hyperbolic}
Let $x$, $y$, and $z$ be points in a geodesic metric space $X$ and let
$\Delta$ be a triangle of geodesics connecting $x$ to $y$, $y$ to $z$
and $x$ to $z$.  This geodesic triangle is \emph{$\delta$-thin}
if each of the sides is contained in a $\delta$-neighborhood of the
union of the other two.  If there is a uniform $\delta$ such that
every geodesic triangle in $X$ is $\delta$-thin, then $X$ is a
{\em $\delta$-hyperbolic space}.  A group $G$ generated by a finite set
$A$ is \emph{word-hyperbolic} if its Cayley graph $\Gamma(G,A)$
is $\delta$-hyperbolic.
\end{defn}

One of the key properties of $\delta$-hyperbolic spaces is that
geodesics and quasigeodesics stay uniformly close in the following
sense:

\begin{lem}[Proposition 3.3 of \cite{ABC90}]\label{lem:close}
Let $x$ and $y$ be points in the $\delta$-hyperbolic metric space
$X$. Then there are integers $L(\lambda,\epsilon)$ and
$M(\lambda,\epsilon)$ such that if $\alpha$ is a
$(\lambda,\epsilon)$-quasigeodesic between the points $x,y$ and
$\gamma$ is a geodesic $[xy]$, then $\gamma$ is contained in an
$L$-neighborhood of $\alpha$ and $\alpha$ is contained in an
$M$-neighborhood of $\gamma$.
\end{lem}

It is easy to deduce from Lemma~\ref{lem:close} that the property of
being $\delta$-hyperbolic for some $\delta$ is preserved by
quasi-isometries between geodesic metric spaces, even though the
specific value of $\delta$ may have to be changed.  Combined with
Corollary~\ref{cor:qi}, this shows that the property of a group being
word-hyperbolic is independent of the choice of a finite generating
set.

\begin{cor}\label{cor:preserves quasi}
If $X$ and $X'$ are geodesic metric spaces, $X$ is
$\delta$-hyperbolic, and $\phi\colon X \rightarrow X'$ is a
quasi-isometry, then a subspace $Y$ is quasiconvex in $X$ if and only
if $\phi(Y)$ is quasiconvex in $X'$.
\end{cor}

As a consequence, the quasiconvexity of a subgroup in a
word-hyperbolic group does not depend on the generating set.

\begin{cor}\label{cor:independent}
Let $H$ be a subgroup of the word-hyperbolic group $G$ and let $A$ and
$B$ be finite generating sets for $G$.  The subgroup $H$ will be
quasiconvex in $\Gamma(G,A)$ if and only if $H$ is quasiconvex in
$\Gamma(G,B)$.  In particular, when $G$ is word-hyperbolic, the
quasiconvexity of a subgroup is independent of the choice of finite
generating set for $G$.
\end{cor}

Thus, for word-hyperbolic groups there is the following partial
converse to Lemma~\ref{lem:qc->qie}.

\begin{lem}\label{lem:qie->qc}
Let $G$ be a word-hyperbolic group with finite generating set $A$ and
let $H$ be a subgroup of $G$ with finite generating set $B$.  If $H_B$
is quasi-isometrically embedded in $G_A$, then $H$ is a quasiconvex
subgroup of $G$.
\end{lem}

\section{Fan quasiconvexity theorems}\label{sec:lqc-main}
In this section we prove our main technical results about local
quasiconvexity.  Since the reader has already seen arguments utilizing
perimeter-reducing fans in Section~\ref{sec:fanCT}, we will treat only the
general fan case.  The reader not yet completely comfortable with the
language of fans should keep in mind the special case of
a fan consisting of a single $2$-cell $R$ with
outer path $Q$ in $\partial R$.

\begin{defn}[Straightening]\label{def:qrh}
Let $X$ be a weighted $2$-complex, let $\fantypes$ be a collection of
fans in $X$, and let $K$ and $\epsilon$ be constants.  A path $P
\rightarrow X$ can be \emph{$(K,\epsilon)$-straightened} if there
exists a sequence of paths $\{P=P_1, P_2, \ldots, P_t\}$ such that for
each $i$, $P_{i+1}$ is obtained from $P_i$ by either removing a
backtrack, or by pushing across a fan in $\fantypes$.  In addition,
the final path $P_t$ must satisfy the following condition:
Consider the lift of $P_t$ to $\widetilde{X}$ and let $d$ denote the
length of a geodesic in $\widetilde{X}^{(1)}$ with the same endpoints.
There must exist a path $P'\rightarrow \widetilde{X}$ with the same
endpoints as $P_t\rightarrow \widetilde{X}$ such that $P'$ lies in a
$K$-neighborhood of $P_t$ and such that
\begin{equation}\label{eq:length-condition}
K\cdot d +\epsilon > |P'|.
\end{equation}
If every fan in $\fantypes$ is spread-out and perimeter-reducing, and
if for some fixed choice of $K$ and $\epsilon$, every path $P
\rightarrow X$ can be $(K,\epsilon)$-straightened, we say that $X$
satisfies the \emph{straightening hypothesis}.
\end{defn}

The following is our main technical result about the
straightening hypothesis.

\begin{thm}[Subgroups quasi-isometrically embed]\label{thm:lqc}
Let $X$ be a compact weighted $2$-complex.  If $X$ satisfies the
straightening hypothesis, then every finitely generated subgroup of
$\pi_1X$ embeds by a quasi-isometry.  Furthermore, if $\pi_1X$ is
word-hyperbolic then it is locally quasiconvex.
\end{thm}

\begin{proof}
Let $G = \pi_1X$, let $H$ be a subgroup of $G$ which is generated by a
finite set $B$ of elements, and let $\widehat X$ be the based covering
space of $X$ corresponding to the inclusion $H \subset G$.  Next, let
$C$ be a wedge of finitely many circles, one for each generator in
$B$, and let $\phi\colon C \rightarrow X$ be a map which sends each
circle to a based path in $X^{(1)}$ representing its corresponding
generator.  This map lifts to a map $C \rightarrow \widehat X$, and we
let $Y_0$ denote the image of $C$ in $\widehat X$.  Since the weights
on the sides of the $2$-cells are nonnegative integers, the perimeter
of any compact subcomplex of $\widehat X$ is finite and nonnegative.
In particular, there exists some compact connected subcomplex $Y
\subset \widehat X$ which contains $Y_0$ and which does not admit any
fan perimeter reductions.  For instance, we can choose $Y$ to be of
minimal perimeter among all compact connected subcomplexes of
$\widehat X$ containing $Y_0$.  If $Y \rightarrow \widehat X$ were to
admit a fan perimeter reduction, then Lemma~\ref{lem:fan-attachments}
and Lemma~\ref{lem:surj} would allow us to create a slightly larger
subcomplex which had a strictly smaller perimeter, contradicting the
way $Y$ was chosen.

Let $Z$ be a $K$-neighborhood of $Y$ in $\widehat X$, let $\widetilde
X$ be the based universal cover of $X$, let $\widetilde Y$ denote the
based component of the preimage of $Y$ in $\widetilde X$, and let
$\widetilde Z$ denote the based connected component of the preimage of
$Z$ in $\widetilde X$.  Since $Z \supset Y \supset Y_0$, we see that
$Z$ contains a set of paths which generate $\pi_1\widehat X = H$ and
thus the action of $H \subset G$ on $\widetilde X$ stabilizes
$\widetilde Z$.  In particular, the preimages of the basepoint of
$\widehat X$ in $\widetilde X$ are contained in $\widetilde Z$ and
these $0$-cells are in one-to-one correspondence with the elements of
$H$.  Using this correspondence, we will consider $H$ as a subspace of
$\widetilde X$.  Let $H_{\widetilde{Z}}$ be the metric on $H$ defined
by the $1$-skeleton of $\widetilde{Z}$.  Specifically, define
$d(h,h')$ to be the length of the shortest path in $\widetilde Z$
between the appropriate $0$-cells of $H$.

Since $Y$ does not admit any fan perimeter reductions and $Z$ is the
$K$-neighborhood of $Y$, the straightening hypothesis allows us to
conclude that every pair of points in $H \subset \widetilde X$ is
connected by a path in $\widetilde Z$ which satisfies
Equation~(\ref{eq:length-condition}).  In particular, given a path $P$
in $\widetilde Y$ connecting a pair of points in $H$, we can follow
the sequence of alterations to obtain paths $P_1, \ldots, P_t$ without
leaving the subcomplex $\widetilde Y$, and since the path $P'$ lies in
a $K$-neighborhood of $P_t$, the path $P'$ does not leave the
subcomplex $\widetilde Z^{(1)}$.  Finally, since this is true for all
pairs of points in $H \subset \widetilde X$, this shows that the map
$H_{\widetilde Z} \rightarrow \widetilde X^{(1)}$ is a
$(K,\epsilon)$-quasi-isometric embedding.

Since $X$ is compact, $G$ has some finite generating set $A$.  More
specifically, if we select a maximal spanning tree for $X^{(1)}$, then
a generator corresponding to each $1$-cell not in the spanning tree is
sufficient.  Consider the following diagram of maps between metric
spaces where the metric spaces $G_A$ and $\Gamma(G,A)$ are the metric
on the group and the metric on its Cayley graph.
\[
\begin{array}{ccccc}
H_B & \rightarrow & G_A & \rightarrow & \Gamma(G,A)\\
\downarrow & & \downarrow &&\\
H_{\widetilde{Z}} & \rightarrow & \widetilde X^{(1)} &&
\end{array}
\]
We have shown that the bottom map is a quasi-isometric embedding.
Since $G$ acts properly discontinuously and cocompactly on $\widetilde
X^{(1)}$, by Lemma~\ref{lem:nice-action} the map $G_A \rightarrow
\widetilde X^{(1)}$ is a quasi-isometry.

Next, as remarked above, the action of $H \subset G$ on $\widetilde X$
stabilizes $\widetilde Z$.  The action of $H$ on $\widetilde Z^{(1)}$
is clearly properly discontinuous since it is a restriction of $G$ on
$\widetilde X$, and it is cocompact since the quotient of $\widetilde
Z$ by $H$ is the compact space $Z$.  Thus, by
Lemma~\ref{lem:nice-action} the map $H_B \rightarrow \widetilde
H_{\widetilde{Z}}$ is also a quasi-isometry.  Combining these three
maps we see that the map $H_B \rightarrow G_A$ is a quasi-isometric
embedding.  Since by Corollary~\ref{cor:qi}, $G_A \rightarrow
\Gamma(G,A)$ is also a quasi-isometry, the map from $H_B$ to
$\Gamma(G,A)$ is a quasi-isometric embedding as well.  Finally, if $G$
is word-hyperbolic then it follows from Lemma~\ref{lem:qie->qc} that
$H$ is quasiconvex.
\end{proof}

We are unable to answer the following problem about the relationship between
the straightening hypothesis and word-hyperbolicity.
However, one can show the answer is affirmative if one adds to
Definition~\ref{def:qrh}
the requirement that $P_t$ lie in a $K$-neighborhood of $P'$.
\begin{question}[Straightening and Hyperbolicity]\label{prob:qrh->hyp}
Suppose the compact weighted $2$-complex $X$ satisfies the
straightening hypothesis with respect to some finite collection
$\fantypes$ of fans in $X$. Is $\pi_1X$ word-hyperbolic?
\end{question}

Two conditions which immediately imply the straightening hypothesis
are a decrease in length and a decrease in area.  More explicitly, if
for every immersed path $P\rightarrow X$ which does not lift to a
$(K,\epsilon)$-quasigeodesic, the path $P$ can be pushed across a
perimeter-reducing fan to obtain a new path of strictly smaller
length, then $X$ will satisfy the straightening hypothesis. Indeed the
sequence of reductions cannot continue indefinitely because the length
decreases each time, and thus they terminate at a path $P_t = P'$
which lifts to a quasigeodesic. The most important condition
that implies the straightening hypothesis will involve the following notion:

\begin{defn}[$J$-thin]\label{def:J-thin}
A disc diagram $D$ with boundary cycle $PQ^{-1}$ is called {\em
$J$-thin} for some $J\in\N$, if every $0$-cell in $P$ is contained in
a $J$-neighborhood of $Q$ and vice-versa.
\end{defn}

\begin{thm}[Diagrammatic local quasiconvexity criterion]
\label{thm:fan-diagramLQC}
Let $\fantypes$ be a finite collection of perimeter-reducing spread-out fans
in the compact weighted $2$-complex $X$,
and let $J\in \N$.
Suppose that for every minimal area disc diagram $D\rightarrow
X$ with boundary cycle $PQ^{-1}$,
 either $D$ is $J$-thin or $D$
contains a spur or fan in $\fantypes$ whose outer path is a
subpath of either $P$ or $Q$.
Then $\pi_1X$ is a locally quasiconvex word-hyperbolic group.
\end{thm}

\begin{proof}
We will first give an argument that shows that $X$ satisfies
the straightening hypothesis. We will then apply
a special case of this argument to see that $\pi_1X$ is word-hyperbolic.
The result will then follow from Theorem~\ref{thm:lqc}.

Let $P \rightarrow X$ be a path, let $P \rightarrow \widetilde X$ be a
lift of this path to the universal cover $\widetilde X$, and let $Q
\rightarrow \widetilde X$ be a geodesic in $\widetilde X$ with the
same endpoints.  Since the path $PQ^{-1}$ is a closed null-homotopic
path in $\widetilde X$, its projection to $X$ has the same
properties. Let $D\rightarrow X$ be a minimal area disc diagram with
boundary cycle $PQ^{-1}$.

Let $D$ be oriented so that the path $P$ proceeds from left to right
across the top of the diagram and the path $Q$ proceeds from left to
right across the bottom. Using the diagram $D$, we will now construct
an explicit sequence of paths $P_i$ demonstrating the straightening
hypothesis.  Along the way we will need to define a sequence of paths
$Q_i$, a sequence of diagrams $D_i$, and a sequence of diagrams $E_i$
as well.  The idea will be to systematically remove portions of $D$
from the top and bottom.  At each stage of this process the paths
along the top and bottom will be $P_i$ and $Q_i$, the diagram between
$P_i$ and $Q_i$ will be $D_i$, and the diagram bounded by $Q_iQ^{-1}$
will be $E_i$.  At the end we will reach a diagram $D_t$ with boundary
paths $P_t \rightarrow D_t$ and $Q_t \rightarrow D_t$ such that
$P_tQ_t^{-1}$ is the boundary cycle of $D_t$, and $D_t$ is $J$-thin.
At each stage $D_i$ and $E_i$ will be subdiagrams of $D$.  The reader
is referred to Figure~\ref{fig:lqcdiagram} for an illustration of the
diagram $D$ as well as some of the relevant paths and subdiagrams
appearing in the final situation.

\begin{figure}\centering
\includegraphics[width=.7\textwidth]{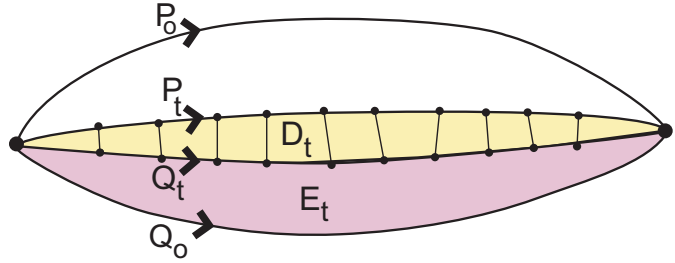}
\caption{
The final subdiagram $D_t$ of $D$ illustrated above is $1$-thin.
\label{fig:lqcdiagram}}
\end{figure}

To begin the process, let $P_0 = P$, $Q_0 = Q$, $D_0 = D$, and let
$E_0$ be the diagram without $2$-cells consisting of the path
$Q$. This is a diagram since $Q$ is a geodesic in $\widetilde X$ and
thus simple in $D$.  Since Definition~\ref{def:qrh} permits the
removal of backtracks, we may assume that our initial path
$P\rightarrow D$ is immersed.  Now for each $i$, since both $P_i
\rightarrow D_i$ and $Q_i \rightarrow D_i$ are immersed paths, by
hypothesis the diagram $D_i$ is either $J$-thin, or $D_i$ contains a
fan $F \in \fantypes$ whose boundary path is a subpath of either $P_i$
or $Q_i$.  Removing the interior of the fan and the interior of its
boundary path from $D_i$ defines a new connected and simply-connected
diagram $D_{i+1}$ where the new boundary paths across the top and
bottom of $D_{i+1}$ are $P_{i+1}$ and $Q_{i+1}$. Notice that by
hypothesis, the path $P_{i+1}$ is either identical to $P_i$ or it is
obtained from $P_i$ by a perimeter-reducing push across a fan, and
likewise, $Q_{i+1}$ is either identical to $Q_i$ or it is obtained
from $Q_i$ by a perimeter-reducing push across a fan.  In either case,
the fact that the original path is an immersion implies that the new
path is an immersion as well.  Since the total number of cells in the
diagrams $D_i$ is decreasing, this process must eventually terminate
at a diagram $D_t$ which is $J$-thin.

It only remains to show that the final path $P_t$ lies close to a path
$P'$ which satisfies the length condition mentioned in
Definition~\ref{def:qrh}.
Let $L\geq 1$ be a bound on the perimeter of a $1$-cell in $X$ (which
is compact), let $f$ be a bound on the lengths of inner paths of fans
in $\fantypes$, and let $K$ be the larger of $J$ and $fL+1$.
Let $P' = Q_t$ and note that since $D_t$ is $J$-thin, $P'$ lies in a
$K$-neighborhood of $P_t$.

We will now show that $\size{Q_t} \leq K\size{Q}$. Let
$\packet{D}_i\to X$ and $\packet{E}_i\to X$ be the packed versions of
$D_i\to X$ and $E_i\to X$.  Now since $Q_{i+1}$ is obtained from $Q_i$
by a perimeter-reducing push across a fan, we have
$\perimeter(\packet{E}_{i+1}) < \perimeter(\packet{E}_i)$, and
consequently, $t \leq \perimeter(\packet{E}_0) = \perimeter(Q) \leq L
\size{Q}$.  Finally, since $\size{Q_{i+1}} \leq f+\size{Q_i}$ for each
$i$, we have $\size{Q_t}\leq ft+ \size{Q} \leq (fL+1)\size{Q} \leq
K\size{Q}$.

We will prove that $\pi_1X$ is word-hyperbolic by showing that
$\pi_1X$ has a linear isoperimetric function \cite{ABC90}.  We will
apply the above argument except that we will exchange the roles of $P$
and $Q$, so that $Q$ is an arbitrary immersed null-homotopic path and
$P$ is the trivial path.  First observe that since no outer path of a
fan can be a subpath of $P$, we see that $P=P_1=P_t$.  We will now use
the fact that $\area(Q)\leq \area(Q_t)+\area(E_t)$ to estimate
$\area(Q)$.

Let $m$ be the maximal number of $2$-cells occurring in a fan in $\fantypes$.
Then $E_t$ was obtained from $E_1$ by adding at most $\perimeter(E_0)$ fans,
and so $\area(E_t)\leq mt \leq mL\size{Q}$.
Since the diagram $D_t$ is $J$-thin, the path $Q_t$ is contained in a
$J$-neighborhood of the path $P_t=1$. Letting $C$ denote the ball of
radius~$J$ in $\widetilde X$, Lemma~\ref{lem:cell-bound} (established
below) implies that $\area(Q_t)\leq M\size{Q_t}$ where $M$ is a
constant that depends only on $X$ and $J$.
But $\size{Q_t}\leq K\size{Q}$ so $\area(Q)\leq \area(E_t)+\area(Q_t)
\leq mL\size{Q}+ MK\size{Q}$, and the isoperimetric function is linear
as claimed.
\end{proof}

\begin{lem}\label{lem:cell-bound}
If $\widetilde X$ is a simply-connected $2$-complex and $C$ is a
compact, connected subspace of $\widetilde X$, then there is a
constant $M$, depending only on $C$, such that for every closed path
$Q\rightarrow C$, we have $\area(Q)\leq M \size{P}$.
\end{lem}

\begin{proof}
Let $S_1, \ldots, S_j$ be the finitely many simple closed nontrivial
paths in $C$, and for $1\leq i\leq j$ let $M_i=\area(S_i) =
\area(D_i)$ where $D_i\rightarrow X$ is a minimal area disc diagram
with boundary cycle $S_i$.  Such a $D_i$ exists for each $S_i$ since
$X$ is simply-connected.  Let $M$ be the maximum value of
$\frac{M_i}{\size{S_i}}$.  Intuitively, $M$ measures the maximum
number of $2$-cells needed per $1$-cell in a simple closed path.  The
necessary inequality is now easy to establish by breaking any closed
path $Q\rightarrow C$ into simple closed paths and backtracks,
creating disc diagrams for each individually, and then reassembling
them into a disc diagram $D$ for the original path.  Since the maximum
number of $2$-cells needed per $1$-cell is bounded by $M$ for each
individual portion, $\area(Q)\leq\area(D)\leq M\size{P}$ as claimed.
\end{proof}

We close this section with a problem analogous to
Problem~\ref{prob:coherence}:

\begin{question}\label{question:lqc}
Let $X$ be a compact weighted $2$-complex such that $\pi_1X$ is
word-hyperbolic. Let $K$ and $\epsilon$ be fixed constants, and
suppose that for every path $P \rightarrow X$, there exists a
$(K,\epsilon)$-quasigeodesic $Q$ and a disc diagram $\Delta$ with
boundary cycle $PQ^{-1}$, such that $\perimeter(\Delta) <
\perimeter(P)$.  Is $\pi_1X$ locally quasiconvex?
\end{question}

\section{Small cancellation II}\label{sec:sc-lqc}
In this section we apply the local quasiconvexity result to various
small cancellation groups.

\begin{thm}[Local quasiconvexity]\label{thm:lqc-sc-1}
Let $X$ be a weighted $2$-complex which satisfies $C(6)$-$T(3)$
$[C(4)$-$T(4)]$.  Suppose
$\perimeter(S) < n\weight(R)$
 for each $2$-cell $R \rightarrow X$
and path $S\rightarrow \partial R$
which is the concatenation
at most three~$[$two$]$ consecutive pieces in the
boundary of $R$.
Then $\pi_1X$ is coherent.
\end{thm}

\begin{proof}
This follows immediately from Corollary~\ref{cor:fan-classification}
and Theorem~\ref{thm:fan-diagramLQC}.
\end{proof}

\begin{exmp}[Surface groups]\label{surfacegps}
Let $X$ be the standard $2$-complex of the presentation $\langle a_1,
\dots, a_g \mid a_1^2 a_2^2 \ldots a_g^2 \rangle$. Then $X$ is the
usual cell structure for the nonorientable surface of genus $g$.
Clearly for any $g\geq 2$, $X$ satisfies $C(2g)$-$T(4)$, and the
pieces are of length~$1$.  Using the unit perimeter we see that the
weight of any piece ($1$-cell) is $2$ and the weight of the $2$-cell
is $2g$.  Thus by Theorem~\ref{thm:coh-sc-1}, $\pi_1X$ is coherent for
$g \geq 2$ and by Theorem~\ref{thm:lqc-sc-1}, $\pi_1X$ is locally
quasiconvex for $g>2$.

A similar result holds if we let $X$ be the standard $2$-complex of
the presentation
\[\langle a_1, b_1, \dots, a_g, b_g\mid
[a_1,b_1][a_2,b_2] \ldots [a_g,b_g] \rangle\] \noindent so that $X$ is
an orientable surface of genus $g$.  For $g \geq 1$ the $2$-complex
$X$ satisfies $C(4g)$-$T(4)$, the pieces are of length $1$, the weight
of each piece is $2$, and the weight of the $2$-cell is $4g$.  Thus by
Theorem~\ref{thm:coh-sc-1}, $\pi_1X$ is coherent for $g \geq 1$ and by
Theorem~\ref{thm:lqc-sc-1}, $\pi_1X$ is locally quasiconvex for $g
>1$.
\end{exmp}

The fact that these methods can be used to prove the coherence and
local quasiconvexity of surface groups is to be expected since the
boundary of a $2$-manifold was one of the original motivations for the
notion of perimeter introduced in Section~\ref{sec:perimeter}.  Here
is a more novel application of Theorem~\ref{thm:coh-sc-1}.

\begin{thm}\label{thm:few-occurrences}
Let $G = \langle a_1,\ldots \mid R_1,\ldots \rangle$ be a
presentation that satisfies $C'(1/n)$.  If each $a_i$
occurs at most $n/3$ times among the $R_j$, then $G$ is coherent and
locally quasiconvex.
\end{thm}

\begin{proof}
In case $n<5$ this is obvious, because any generator appears at most
$5/3$ times, and consequently at most once, which implies that the
group is free.  On the other hand, when $n\geq 5$, the complex $X$
satisfies $C(6)$.  Use the unit weighting, and let $S$ be a subpath of
the $2$-cell $R_j$ consisting of three consecutive pieces.  The small
cancellation assumption implies that $|S| < 3 \cdot
\frac{1}{n}|\partial R_j|$.  On the other hand, the bound on the
number of occurrences of each generator shows that $\perimeter(e) \leq
\frac{n}{3}$ for each $1$-cell $e$.  Thus $\perimeter(S) < 3
\frac{n}{3} \frac{1}{n} |\partial R_j| = |\partial R_j|$.  By
Theorem~\ref{thm:coh-sc-1} the group is coherent and by
Theorem~\ref{thm:lqc-sc-1} it is locally quasiconvex.
\end{proof}

In Section~\ref{sec:1relator} we showed that certain one-relator
groups with torsion are coherent.  Since it was shown in \cite{Pr83}
that $\langle a_1,\dotsc \mid W^n\rangle$ satisfies the $C(2n)$ small
cancellation condition, we can apply Theorem~\ref{thm:lqc-sc-1} to
obtain a local-quasiconvexity result as well.

\begin{thm}\label{thm:lqc-one-rel}
Let $G = \langle a_1,\ldots | W^n \rangle$ be a one-relator group
with $n \geq 3|W|$.  Then $G$ is locally quasiconvex.
\end{thm}

\begin{proof}
Let $X$ be the $2$-complex of this presentation with the unit
weighting.  Since the weight of the unique $2$-cell is $n |W|$, the
weight of its packet is $n^2 |W|$. On the other hand, by the spelling
theorem (Theorem~\ref{thm:spelling}), the length of a piece is less
than $|W|$ and we can assume that each generator appears in $W$ fewer
than $|W|$ times, for otherwise the group is virtually free and the
theorem is obvious.  Consequently, the perimeter of each $1$-cell is
at most $n(|W|-1)$ and since $n\leq |W|$, this is bounded by
$(n-1)|W|$.  So the perimeter of a piece is strictly less than
$(n-1)|W|\cdot |W|$. Thus for $n \geq 3|W|$, $G$ is locally
quasiconvex by Theorem~\ref{thm:lqc-sc-1}.
\end{proof}

A more detailed examination of the local quasiconvexity of one-relator
groups with torsion has been carried out by Hruska and Wise in
\cite{HrWi01}.  As an application of Theorem~\ref{thm:lqc} they are
able to prove the following:

\begin{thm}
Let $G = \langle a_1,\dotsc | W^n \rangle$ be a one-relator group
with $n \geq |W|$.  Then $G$ is locally quasiconvex.
\end{thm}

Our next application is to finitely presented small cancellation
groups with torsion.  As in the one-relator case, small
cancellation groups with sufficient torsion will always be
coherent and locally quasiconvex.  In the proof we will need the
following classical lemma about words in the free group (see
\cite{Gu93}).

\begin{lem}\label{lem:x+y}
Let $X$ and $Y$ be cyclically reduced words in the free group
which are not proper powers.  If $X$ and $Y$ are not cyclic
conjugates, and a word $U$ is both a subword of a power of $X$
and a subword of a power of $Y$, then $|U| \leq |X| + |Y|$.
\end{lem}

\begin{thm}[Power theorem]\label{thm:powers}
Let $\langle a_1, \dotsc \mid W_1, \dotsc \rangle$ be a finite
presentation, where each $W_i$ is a cyclically reduced word which is
not a proper power.  If $W_i$ is not freely conjugate to $W_j^{\pm1}$
for $i \neq j$, then there exists a number $N$ such
that for all choices of integers $n_i \geq N$ the group $G = \langle
a_1, \dotsc \mid W_1^{n_1}, \dotsc \rangle$ is coherent.
Specifically, the number
\begin{equation}\label{eq:powers}
  N = 6 \cdot \frac{|W_{\text{max}}|}{|W_{\text{min}}|} \sum |W_i|
\end{equation}

\noindent
has this property, where $W_{\text{max}}$ and $W_{\text{min}}$ denote
longest and shortest words among the $W_i$, respectively.
Moreover, if $n_i >N$ for all $i$, then $G$ is locally quasiconvex.
\end{thm}

\begin{proof}
We will assume that $\sum |W_i| \geq 2$.  Essentially, the only
case that this assumption eliminates is the presentation $\langle
a_1,\dotsc \mid a_1 \rangle$ and the presentation $\langle
a_1,\dotsc \mid \ \rangle$, and the theorem is trivial in these
cases.
Let $N$ be the number satisfying Equation~(\ref{eq:powers}) and
choose $n_i \geq N$ for all $i$. We will first show that $\langle
a_1,\dotsc \mid W_1^{n_1}, \dotsc \rangle$ satisfies the $C(6)$
condition.

By Lemma~\ref{lem:x+y} the length of the longest piece between the
$2$-cell labeled $W_i^{n_i}$ and the $2$-cell labeled $W_j^{n_j}$
is bounded above by $|W_i| + |W_j| \leq 2 \cdot |W_{\text{max}}|$.
Notice that $2 \cdot |W_{\text{max}}|$ is thus a uniform bound on
the size of a piece which is independent of the size of the
chosen $n_i$. Since the length of the $i$-th $2$-cell is $n_i
\cdot |W_i| \geq n_i \cdot |W_{\text{min}}|$, and since by
assumption,
\[
n_i \cdot |W_{\text{min}}| \geq N \cdot |W_{\text{min}}| = 6
\cdot |W_{\text{max}}| \sum |W_i| \geq 6 \cdot 2 \cdot
|W_{\text{max}}|
\]

\noindent we can conclude that the $C(6)$ condition is satisfied.

Next we will choose a weighting on the sides of the $2$-cells of
the standard $2$-complex for $\langle a_1,\dotsc \mid W_1^{n_1},
\dotsc \rangle$, and verify that the weight criterion of
Theorem~\ref{thm:coh-sc-1} is satisfied.  Let $n = \prod n_i$ and
assign a weight of $n/n_i$ to each of the sides of the $2$-cell
$R_i$ corresponding to the relator $W_i^{n_i}$.  If we let $S$
denote a path in $\partial R_i$ consisting of at most three
consecutive pieces, then we must show that $n_i\weight(R_i) \geq
\perimeter(S)$. This follows from the following string of
inequalities:
\begin{displaymath}
  n_i \cdot \weight(R_i) \geq n_i \cdot n \cdot |W_{\text{min}}|
  \geq 6 \cdot |W_{\text{max}}| \cdot n \sum |W_i| \geq \perimeter(S)
\end{displaymath}

\noindent The first inequality is true since $\weight(R_i) =
\frac{n}{n_i} \cdot |W_i^{n_i}|= n \cdot |W_i|$, which is clearly
greater than or equal to $n \cdot |W_{\text{min}}|$.  The middle
inequality uses the restrictions assumed on $n_i$ and $N$.
Finally, the last inequality is a combination of two
observations: (1) The perimeter of a single $1$-cell will be at
most the sum of the weights of all of the sides in the complex,
so that $\weight(e) \leq \sum \frac{n}{n_i} |W_i^{n_i}| = n \sum
|W_i|$.  And (2) since a piece has length at most $2 \cdot
|W_{\text{max}}|$, the perimeter of three consecutive pieces is
bounded above by $6 |W_{\text{max}}| \cdot n \sum |W_i|$.
Together these show the final inequality.  The weight criterion of
Theorem~\ref{thm:coh-sc-1} is thus satisfied and the group $G'$ is
coherent.  Similarly, if the inequalities are strict, then by
Theorem~\ref{thm:lqc-sc-1}, it is locally quasiconvex.
\end{proof}

Although there is some overlap between the groups studied in this
article and those studied by Feighn and Handel in
\cite{FeHa99}, the methods and the results are distinct.
One indication of this is that all of the groups in
\cite{FeHa99} are indicable (i.e. admit a homomorphism
onto $\Z$), whereas Theorem~\ref{thm:powers} can be used to
construct coherent groups which are perfect.

\begin{cor}\label{cor:perfect}
There exist perfect groups which satisfy the perimeter reduction
hypothesis and are thus coherent.
\end{cor}

\begin{proof}
The following example illustrates the idea: Consider the
presentation
$$\langle a, b \mid a, aaba^{-1}b^{-1}, b,
bbab^{-1}a^{-1} \rangle.$$
 If $N$ is chosen to satisfy
Equation~(\ref{eq:powers}), then the following group is coherent by Theorem~\ref{thm:powers}.
$$G = \langle a,b \mid
a^N, {(aaba^{-1}b^{-1})}^{N+1}, b^N, {(bbab^{-1}a^{-1})}^{N+1}
\rangle$$
But the following presentation for the abelianization of $G$ shows that it is trivial:
$$\langle a,b\mid [a,b],
a^N, a^{N+1}, b^N, b^{N+1}\rangle$$
\end{proof}

\section{$3$-manifold groups}\label{sec:manifold}
In this section, we use the theorems about coherence and local
quasiconvexity in small cancellation groups to show that a large
family of $3$-manifold groups are coherent and locally
quasiconvex. We begin with a theorem about branched covers of
$2$-complexes.

\begin{thm}\label{thm:branched 2-complex}
Let $X$ be a compact $2$-complex, and suppose that no $2$-cell of
$X$ is attached by a proper power. Then there exists a constant
$d$ depending only on $X$ such that for every branched cover
$\widehat X \rightarrow X$ where (1) the branching is over the
barycenters of $2$-cells of $X$, and (2) all of the branching
degrees are at least $d$, the fundamental group $\pi_1\widehat X$
is coherent.  Similarly, there is another constant $d$ such that
the compact branched covers satisfying these conditions have a
locally quasiconvex fundamental group.
\end{thm}

\begin{proof}
The proof is similar to the proof of Theorem~\ref{thm:powers},
which can actually be deduced from it.  First notice that it is
sufficient to consider the case where no two $2$-cells have the
same attaching maps.  It follows from Lemma~\ref{lem:x+y} that
there is a bound depending only on $X$ on the length of a piece
in the boundary of a $2$-cell of $\widehat X$.  Thus the lengths
of the pieces remain bounded as the branching degrees increase,
whereas the length of the boundary of each $2$-cell grows
linearly with $d$.  It follows that for large $d$ the
presentation satisfies small cancellation conditions.

Using the unit perimeter, the maximum sum of the weight of three
consecutive pieces remains constant, whereas the weight of each
$2$-cell grows linearly with $d$.  Thus, for sufficiently
large~$d$ the weight criteria for coherence
(Theorem~\ref{thm:coh-sc-1}) and for local quasiconvexity
(Theorem~\ref{thm:lqc-sc-1}) will be satisfied.
\end{proof}

\begin{rem}\label{rem:numerous branched}
Observe that finite branched covers with high branching degrees
correspond to certain finite index subgroups of the fundamental group
of the space we obtain when we remove the barycenter of each
$2$-cell. Since the fundamental group of this space is free, and thus
residually finite (\cite{LySch77}), these types of covers are
numerous.
\end{rem}

\begin{thm}\label{thm:branched manifold}
Let $M$ be a compact $3$-manifold with a combinatorial cell
structure. There exists $d$ depending on $M$ such that the
following holds: Let $B\rightarrow M\setminus M^{(0)}$ be a
branched cover with at least $d$~fold branching along each
$1$-cell of $M$.  Then $\pi_1B$ is coherent or even locally
quasiconvex.
\end{thm}

\begin{proof}
Let $M_1$ denote the underlying manifold of $M$ equipped with the
cell structure dual to $M$.  The branched covers $B\rightarrow
M\setminus M^{(0)}$ branched along $1$-cells of $M$ correspond to
branched covers $B' \rightarrow M_1^{(2)}$ branched along centers
of $2$-cells of $M_1^{(2)}$.  Note that $B'$ is a subspace of
$B$.  Furthermore, the branching degree at each $1$-cell of $M$
is the same as the branching degree at the center of the
corresponding $2$-cell of $M_1^{(2)}$.

The obvious strong deformation retraction of $M\setminus M^{(0)} $
onto $M_1^{(2)}$ induces a strong deformation retraction of $B$
onto $B'$, so that $\pi_1B \cong \pi_1B'$.  It is therefore
sufficient to prove the analogous result for branched covers
$B'\rightarrow M_1^{(2)}$ where the branching occurs over the
centers of $2$-cells of $M_1^{(2)}$, and this is exactly what was
proved in Theorem~\ref{thm:branched 2-complex}.
\end{proof}

We illustrate Theorem~\ref{thm:branched manifold} with the
following example.

\begin{exmp}\label{exmp:branched cube}
Let $M$ denote the usual cell division for the $3$-torus $T^3 =
S^1\times S^1\times S^1$.  We will show that any branched cover of
$M\setminus M^{(0)}$ along the $1$-cells of $M$ has a coherent
fundamental group provided that the degree~$d$ of branching is
$\geq 3$.

First observe that since $M_1$ is obviously isomorphic to $M$, it
is easy to see that $M_1^{(2)}$ contains exactly three $2$-cells,
each of which is a square.  Each piece of $M_1^{(2)}$ has
length~$1$ and perimeter~$4$.
If each branching degree is $\geq d$ then the presentation
satisfies $C(4d)$.  Since the perimeter of three consecutive
pieces is $3 \cdot 4$ and the weight of each $2$-cell at least
$4d$, the fundamental group will be coherent when $d \geq 3$ (by
Theorem~\ref{thm:coh-sc-1}) and locally quasiconvex when $d > 3$
(by Theorem~\ref{thm:lqc-sc-1}).
\end{exmp}

As mentioned in the introduction, Scott and Shalen proved that all
$3$-manifold groups are coherent, so the coherence assertion in
Theorem~\ref{thm:branched manifold} is certainly not new.
Nevertheless, it is interesting to be able to recover this
special case from the different point of view of this paper.

The local quasiconvexity assertion is a bit trickier to obtain using
prior results.  However, it seems likely that the branched covers of
Theorem~\ref{thm:branched manifold} can be constructed from hyperbolic
$3$-manifolds with non-empty boundary and no cusps, by gluing along
annuli.  Thus the local quasiconvexity appears to follow from the
following theorem of Thurston's \cite{Th77} which we quote from
\cite[Proposition~7.1]{Mo84}, together with the theorem due to Swarup
\cite{Sw93} that in the case that $N$ has no cusps, geometrically
finite subgroups are quasiconvex.

\begin{thm}[Thurston]\label{thm:atoroidal}
Let $N$ be a geometrically finite hyperbolic manifold such that
$\partial\text{Core}(N)$ is nonempty. Then every covering space
$N'$ of $N$ with a finitely generated fundamental group is also
geometrically finite.
\end{thm}

We believe that the branched covers with branching degree $\geq
2$ of the $3$-torus in Example~\ref{exmp:branched cube} are
atoroidal hyperbolic $3$-manifolds with boundary, and hence the
local quasiconvexity does follow from Theorem~\ref{thm:atoroidal}.

\section{Related properties}\label{sec:other}
There are a number of properties of groups which are closely related
to coherence and local quasiconvexity.  In this final section we
examine briefly how our techniques can be used to produce results
about three of these related topics: Howson's property, finitely
generated intersections with Magnus subgroups, and the generalized
word problem.

\subsection{Finitely-generated intersections}
A group is said to satisfy the \emph{finitely generated intersection
property} (or f.g.i.p.) if the intersection of any two finitely
generated subgroups is also finitely generated.  In 1954 Howson proved
that free groups have the f.g.i.p., and as a result this property is
sometimes referred to as Howson's property.  As was shown in
\cite{Sh91}, every quasiconvex subgroup is finitely generated and the
intersection of any two quasiconvex subgroups is again quasiconvex.
Combining these two facts, one sees that every locally quasiconvex
group satisfies the f.g.i.p.  In particular, all of the groups we have
shown to be locally quasiconvex, also have Howson's property.

In this subsection we show how the algorithm of
Section~\ref{sec:algorithms} can be used to explicitly construct the
finitely generated intersection of two finitely generated subgroups
using the perimeter techniques we have already introduced.
The method is a 2-dimensional generalization of that described by
Stallings for graphs \cite{St83}. The
construction of the finitely generated intersection will proceed in
two steps.  The first step will be to reduce this property to a
property of the fiber product of the spaces corresponding to these
subgroups.  The second step will be to show how this property can be
achieved using perimeter reductions.

\begin{defn}[Fiber products]
Let $X$ be a complex, let $A\rightarrow X$ and $B\rightarrow X$ be
maps, and let $D \subset X\times X$ be the {\em diagonal} of $X \times X$
so $D= \{ (x,x) \mid x\in X \}$.  The
\emph{fiber-product} $A\otimes B \rightarrow X$ is defined to be the
subspace of $A\times B$ which is the preimage of $D$ in $A\times B
\rightarrow X\times X$. Identifying $D$ with $X$, there is a natural
map $A\otimes B \rightarrow X$, and in fact, the following diagram
commutes:
\[
\begin{array}{ccc}
A  \otimes  B   &   \rightarrow &    B   \\
\downarrow        &   \searrow    &   \downarrow   \\
A               &   \rightarrow  &    X
\end{array}
\]
Notice that if both of the maps $A\rightarrow X$ and $B\rightarrow X$
have compact domain, then so does $A\otimes B \rightarrow X$, and if
both maps are immersions, then so is their fiber product.  Notice also
that if $\phi$ is the map $A \rightarrow X$ and $B$ is a subcomplex of
$X$, then $A\otimes B$ is $\phi^{-1}(B)$.
Note that the map $A\otimes B\rightarrow X$ induces a cell structure
on $A\otimes B$ such that the map $A\otimes B\rightarrow X$ is
combinatorial.  Furthermore, if $A$, $B$, and $X$ are based spaces and
the maps from $A$ and $B$ preserve basepoints, then $X\times X$ and
$A\otimes B$ have natural basepoints.  The \emph{based component} $C$
of $A\otimes B$ is the component of $A\otimes B$ containing this
basepoint. In particular we have $\pi_1(A\otimes B)=\pi_1C$.
\end{defn}

\begin{lem}\label{lem:fiberproduct}
Let $A\rightarrow X$ and $B\rightarrow X$ be based maps and let $a$,
$b$, and $x$ be their basepoints.  If for every element of $\pi_1A
\cap \pi_1B$ there is a closed path $P \rightarrow X$ based at $x$
which lifts to closed paths $P \rightarrow A$ and $P \rightarrow B$
based at $a$ and $b$, then the image of $\pi_1 (A \otimes B)$ in
$\pi_1X$ is the intersection of the images of $\pi_1A$ and
$\pi_1B$. In particular, when $A$ and $B$ are compact, the fiber
product $A\otimes B$ will also be compact, and $\pi_1(A\otimes B)$
is finitely generated.
\end{lem}

\begin{proof}
Since the based component of $A\otimes B$ factors through $A$ and
$B$, the image of its fundamental group must be contained in the
intersection of the images of their fundamental groups.  On the other
hand, by hypothesis, each element in this intersection has a
representative which lifts both to $A$ and to $B$, and thus to
$A\otimes B$ as well.  In particular, $\pi_1(A\otimes B)$ also maps
onto this intersection.
\end{proof}

While we have already proven in Theorem~\ref{thm:lqc-sc-1}
that the following groups
are locally quasiconvex and hence have the finitely generated intersection
property, the following theorem gives an explicit and relatively
efficient method of computing this intersection.

\begin{thm}[f.g.i.p using $i$-shells]\label{thm:fgip}
Let $X$ be a weighted $C(4)$-$T(4)$-complex $[C(6)$-$T(3)]$.  If
every $i$-shell with $i\leq 2$ $[i\leq 3]$ is perimeter-reducing then
$\pi_1X$ has the finitely generated intersection property and the
intersection of two finitely generated subgroups of $\pi_1X$ can be
constructed explicitly.
\end{thm}

\begin{proof}
First note that by Theorem~\ref{thm:fan-classification}, $X$ satisifies
the $2$-cell reduction hypothesis.  Let $H$ and $K$ be finitely
generated subgroups of $\pi_1X$ and let $A_1\rightarrow X$ and
$B_1\rightarrow X$ be compact complexes chosen so that the images of
their fundamental groups are $H$ and $K$.  We will now show by
construction that $H \cap K$ is finitely generated.  Let $A_2
\rightarrow X$ and $B_2\rightarrow X$ be obtained from $A_1$ and $B_1$
by running the perimeter reduction algorithm described in
Section~\ref{sec:algorithms}.  When the algorithm stops, no further
folds or $2$-cell perimeter reductions can be performed on $A_2
\rightarrow X$ or $B_2 \rightarrow X$.  Thus, by the remark at the end
of Definition~\ref{def:prh}, $\pi_1A_2$ and $\pi_1B_2$ can now be
viewed as subgroups of $\pi_1X$.  Next, let $A_3\rightarrow X$ and
$B_3\rightarrow X$ be complexes obtained from $A_2$ and $B_2$ by
attaching to every vertex $v$, a copy of each $2$-cell in $X$ whose
boundary cycle contains the image of $v$ in $X$.  These copies of
$2$-cells are attached only at the vertex $v$, and the maps into $X$
are extended in the obvious way.  Finally, let $A_4\rightarrow X$ and
$B_4\rightarrow X$ be obtained by rerunning the perimeter reduction
algorithm on $A_3$ and $B_3$.  Note that each of these steps adds only
a finite number of $2$-cells, and consequently, since $A_1$ and $B_1$
are compact, so are $A_4$ and $B_4$.  It should also be clear that
$\pi_1A_2 = \pi_1A_3$ and $\pi_1B_2 = \pi_1B_3$, and since running the
perimeter reduction algorithm does not change the image of the
fundamental group in $X$, $\pi_1A_4 = \pi_1A_2 = H$ and $\pi_1B_4 =
\pi_1B_2 = K$.  It only remains to show that $A_4 \rightarrow X$ and
$B_4 \rightarrow X$ satisfy the conditions of
Lemma~\ref{lem:fiberproduct}.

For each element of $H \cap K$ we can choose closed paths $P
\rightarrow A_2$ and $Q \rightarrow B_2$ whose images in $X$ represent
this element.  Since the concatenation
 $P \rightarrow X$ followed by the inverse of $Q
\rightarrow X$ is null-homotopic in $X$, there is a minimal area
disc diagram $D \rightarrow X$ whose boundary cycle is $PQ^{-1}$.  Let
$P\rightarrow A_2$ and $Q\rightarrow A_2$ be chosen (among
combinatorial paths homotopic to them) so that the corresponding disc
diagram $D$ is of minimal area.

Observe that neither $P \rightarrow D$ nor $Q\rightarrow D$ contains a
boundary arc which is the complement of at most $3$ pieces in $D$ ($2$
pieces in the $C(4)$-$T(4)$ case).  To see this, suppose that such a
subpath existed in $P$.  By hypothesis, there would then be a
perimeter reduction on $P \rightarrow X$ which would yield a
corresponding reduction of $P \rightarrow A_2$ (since no perimeter
reductions exists for $A_2 \rightarrow X$). But this would yield a new
path $P$ and a new disc diagram $D$ with fewer $2$-cells. We therefore
conclude by Corollary~\ref{cor:fan-classification}, that $D$ is
$K$-thin.  It is now easy to see that $P \rightarrow X$ actually lifts
to $B_4$ as well as $A_4$.  Thus $A_4 \rightarrow X$ and $B_4
\rightarrow X$ satisfy Lemma~\ref{lem:fiberproduct} and the proof is
complete.
\end{proof}

The fact that Theorem~\ref{thm:fgip} uses the $2$-cell reduction
hypothesis rather than the path reduction hypothesis is crucial.  The
examples below satisfy the path reduction hypothesis but fail to have
the finitely generated intersection property.

\begin{exmp}\label{exmp:fgip-fails}
Let $X$ be the standard $2$-complex of the presentation $F_2 \times \Z
= \langle a,b,t \mid [a,t], [b,t] \rangle$.  If we assign a weight of~$0$
 to every side incident at the edge labeled $t$ and a weight of $1$
to every other side of the complex, the result is a $C(4)$-$T(4)$
complex which satisfies the path reduction hypothesis.  However, the subgroup
$\langle a,b \rangle \cap \langle at,bt\rangle$ is not
finitely generated.

Similarly, if we let $X$ be the standard $2$-complex of the
presentation $\langle a,b,t\mid t^{-1}a^2tb^3\rangle$, and we assign a
weight of $1$ to all of the sides incident at the edge labeled $t$ and
a weight of $0$ to all of the other sides, the result is a
$C(6)$-$T(3)$-complex which satisfies the path reduction hypothesis.
Since this group has a free factor which is commensurable with
$F_2\times Z$, it also fails to have Howson's property.
\end{exmp}

\begin{rem}[f.g.i.p using fans]
We note that the construction above works
to explicitly compute the intersection between two
subgroups of $\pi_1X$ if $X$ satisfies the
more general criterion of Theorem~\ref{thm:fan-diagramLQC}.
The main difference is that one attaches fans in $\fantypes$
instead of $i$-shells, and one adds a diameter~$J$ ``neighborhood''
in passing from $A_2$ and $B_2$ to $A_3$ and $B_3$.
\end{rem}

\subsection{Finitely-generated intersections with Magnus subgroups}
In this subsection we show how the intersections between specific
subgroups can sometimes be shown to be finitely generated under a
weaker set of assumptions.  For example, even if we only assume the
path reduction hypothesis, we can sometimes prove that certain
subgroups $H\subset G$ have the property that $H\cap K$ is finitely
generated for every finitely generated subgroup $K\subset G$.

\begin{thm}\label{thm:fgip-zero}
Let $X$ be a compact based $2$-complex with non-negative integer
weights assigned to the sides of its $2$-cells, and let $M$ be a based
subgraph of $X^{(1)}$ with $\perimeter(M)=0$.  If the weight of each
$2$-cell in $X$ is positive and $X$ satisfies the path reduction
hypothesis, then for any finitely generated subgroup $H$ in $\pi_1X$,
the intersection $\pi_1M\cap H$ is also finitely generated.
\end{thm}

\begin{proof}
Let $\widehat X$ be the based covering space of $X$ corresponding to
the subgroup $H$.  Let $K\subset X$ be a based $\pi_1$-surjective
subcomplex.  According to Theorem~\ref{thm:GCT}, there exists a
compact connected subcomplex $Y \subset \widehat X$ such that
$K\subset Y$, and such that $Y\rightarrow X$ is of minimal perimeter
among all such complexes and hence such that $\pi_1Y\rightarrow
\pi_1\widehat X$ is an isomorphism.  Let $\phi$ denote the map from
$Y$ to $X$, and let $C$ denote the connected component of
$\phi^{-1}(M)$ which contains the basepoint of $Y$.  Note that
$\phi^{-1}(M)$ is the fiber product $Y \otimes M$ of the maps
$Y\rightarrow X$ and $M\hookrightarrow X$.  We will show that
$H=\pi_1(Y\otimes M)=\pi_1(C)$ where $C$ is the based component of
the fiber product.

It is obvious that $\pi_1(Y\otimes M) \subset ( \pi_1Y \cap \pi_1M)$
and so it is sufficient to show that $\pi_1Y\cap \pi_1M \subset
\pi_1(Y\otimes M)$.  Suppose that $P \rightarrow M \subset X$ is a
closed immersed path which represents an element of $\pi_1\widehat X$,
then $P$ lifts to a closed path $\widehat P$ in $\widehat X$.  Let
$Y'=Y\cup\widehat P$ and note that $Y'$ is compact and connected.  Since
$\widehat P$ is a closed path, $Y'$ is formed from $Y$ by adding finitely
many arcs whose endpoints lie in $Y$.  In particular, $Y'-Y$ consists
of a set of $0$-cells and perimeter zero $1$-cells, and so
$\perimeter(Y')=\perimeter(Y)$.  If $Y'\neq Y$ then $\pi_1Y$ is a
proper free factor of $\pi_1Y'$ and hence the map $Y'\rightarrow
\widehat X$ is not $\pi_1$-injective.  Part~A of Theorem~\ref{thm:GCT}
provides a complex $Y_t$ with $Y\subset Y'\subset Y_t$ and
$\perimeter(Y_t) < \perimeter(Y') = \perimeter(Y)$. This contradicts
the minimality of $\perimeter(Y)$.

In conclusion $Y'= Y$ and $P\rightarrow X$ lifts to $Y$.  By
Lemma~\ref{lem:fiberproduct} it now follows that $\pi_1 M \cap H =
\pi_1(C)$, and by the compactness of $C$ this intersection is finitely
generated.
\end{proof}

We note that when $X$ satisfies the strict hypothesis, generators for
an intersection can be computed.  To see that the hypothesis that
$\perimeter(M)=0$ cannot be dropped, let $M$ denote the $\langle
a,b\rangle$ subcomplex of Example~\ref{exmp:fgip-fails}.  The
following application is an example from small-cancellation theory:

\begin{exmp}
In the group $G =\langle a,b,c,d \mid abcd dacb badc \rangle$, the
subgroup generated by $a$ and $b$ has a finitely generated
intersection with every other finitely generated subgroup of $G$.  To
apply Theorem~\ref{thm:fgip-zero} we let $X$ be the standard
$2$-complex of the presentation, and let $M$ be the graph consisting
of the unique $0$-cell of $X$ together with the edges corresponding to
the letters $a$ and $b$.  Note that $X$ is a $C(12)-T(4)$-complex and
that all its pieces are of length~$1$.  Also, if we assign a weight of $0$
to all of the sides present at the edges labeled $a$ and $b$, and a
weight of $1$ to all of the others, then the unique $2$-cell has
weight~$6$, and the perimeter of each piece is at most~$3$.  Thus $X$
satisfies the path reduction hypothesis. \end{exmp}

\begin{cor}\label{cor:Magnus-intersections}
A Magnus subgroup of $G=\langle A | W^n \rangle$ has finitely
generated intersection with any f.g. subgroup of $\langle A | W^n
\rangle$ provided that $n \geq |W|$.
\end{cor}

\begin{proof}
Recall that a Magnus subgroup of $G$ is any subgroup which is
generated by a proper subset $B \subset A$.  Let $X$ be the standard
$2$-complex of $G$ and let $M$ be the graph formed using the unique
$0$-cell and the $1$-cells corresponding to the generators in $B$. We
assign a weight of $0$ to every side incident at an edge in $M$ and a
weight of $1$ to all of the other sides.  To show that $X$, with this
weighting, satisfies the path reduction hypothesis, by
Theorem~\ref{thm:wt-one-rel} we only need to show that the weight of the
packet of the $2$-cell is at most the perimeter of a subword of (a
cyclic conjugate of) $W$.  In order to apply Theorem~\ref{thm:wt-one-rel}, the
weight of $R$ must be strictly positive.  This corresponds to the
existence of a letter in $W$ which is not in $M$.  When this is not
the case, the subgroup generated by $B$ is a free factor of $G$ and
the result follows immediately from the theory of free products.
Thus, we may assume that the weight of $W$ is indeed positive, and
that Theorem~\ref{thm:wt-one-rel} applies.

Let $k$ be the number of times that elements outside of $M$ occur in
$W$.  The weight of the unique $2$-cell is then $nk$, and the weight
of its packet is $n^2 k$.  On the other hand, since the perimeter of a
single edge in $X$ is at most $nk$, the perimeter of a subword of $W$
is at most $nk|W|$.  Thus, the path reduction hypothesis is
satisfied whenever $n \geq |W|$.
\end{proof}

\subsection{Generalized word problem}
Our final related property concerns a generalization of the word
problem for a group.

\begin{defn}[Generalized word problem]
A subgroup of $G$ generated by elements $V_1, \ldots, V_r$ is said to
have \emph{solvable membership problem} provided it is decidable
whether an element $U \in G$ lies in $\langle V_1, \ldots, V_r
\rangle$.  If the membership problem is solvable for every finitely
generated subgroup of $G$, then $G$ is said to have a \emph{solvable
generalized word problem}.  The name alludes to the fact that it
includes the question of membership in the trivial subgroup (the word
problem for $G$) as a particular case.
\end{defn}

We will need the following lemma.

\begin{lem}\label{lem:closing-paths}
Let $X$ be a $2$-complex with non-negative integer weights assigned to
the sides of its $2$-cells.  If $X$ satisfies the (weak) path
reduction hypothesis and $Y \rightarrow X$ is a map which does not
admit a (weak) perimeter reduction or a fold, then every path $P
\rightarrow Y$ whose image is a closed and null-homotopic path in $X$
will be a closed path in $Y$ as well.
\end{lem}

\begin{proof}
Let $P_i \rightarrow X$ ($i=1,\ldots t$) be a sequence of folds and
perimeter reductions which starts with the path $P \rightarrow Y
\rightarrow X$ and ends with the trivial path.  Since $Y \rightarrow X$
does not admit any folds or (weak) perimeter reductions, all of the
alterations to the path $P \rightarrow X$ can be mimicked in the path
$P \rightarrow Y$.  In particular, since the final path $P_t$ is the
trivial path and since the endpoints of $P_i$ are the same throughout this
process, the original path $P \rightarrow Y$ must have been closed.
\end{proof}

\begin{thm}\label{thm:generalizedWP}
Let $X$ be a $2$-complex with non-negative integer weights assigned to
the sides of its $2$-cells.  If $X$ satisfies the path reduction
hypothesis, then $\pi_1X$ has a solvable generalized word problem.
\end{thm}

\begin{proof}
Let $V_1,\ldots, V_r, U$ be a set of closed paths in $X$ with a common
basepoint.  To decide whether the element of $\pi_1X$ represented by
$U$ is in the subgroup $H$ generated by the elements corresponding to
the $V_i$ we proceed as follows.  Let $Y_1$ be the wedge of $r$ closed
paths and a single open path attached to the others at only one of its
endpoints.  We define the map $Y_1
\rightarrow X$ so that it agrees with the maps
$V_i\rightarrow X$ and $U\rightarrow X$.
Let $p$ denote the basepoint of $Y_1$, and let $q$ be the other
endpoint of the open path in $Y_1$ (see
Figure~\ref{fig:membershipwedge}).  Let $Y_t$ be the final complex
produced by running the perimeter reduction algorithm on the map $Y_1
\rightarrow X$.  Recall that the process also constructs a
$\pi_1$-surjective map $Y_1 \rightarrow Y_t$ such that the composition
$Y_1 \rightarrow Y_t \rightarrow X$ is the original map $Y_1
\rightarrow X$.  We claim that vertices $p$ and $q$ are identified
under the map $Y_1 \rightarrow Y_t$ if and only if the element
corresponding to the path $U \rightarrow X$ lies in the subgroup $H$.
Both directions of this implication need to be verified.

\begin{figure}\centering
\includegraphics[width=.5\textwidth]{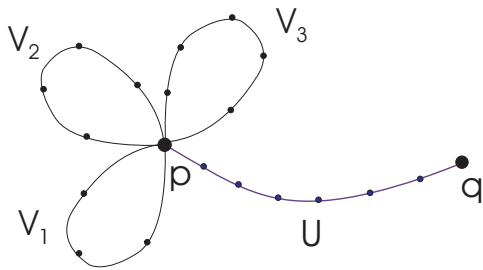}
\caption{The complex $Y_1$ described in
Theorem~\ref{thm:generalizedWP}\label{fig:membershipwedge}}
\end{figure}

First, observe that if $p$ and $q$ are identified in $Y_t$, then $U
\rightarrow Y_1$ is sent to a closed path $U \rightarrow Y_t$.  Since
by the remark at the end of Definition~\ref{def:prh}, $\pi_1Y_t$ can
be considered a subgroup of $\pi_1X$ whose image is the same as the
image of $\pi_1Y_1$ (which is $H$), the element represented by $U$
must lie in the subgroup $H$.

Conversely, suppose that the element represented by $U$ is indeed in
the subgroup $H$.  Then there exists a closed path $V \rightarrow Y_1$
whose image in $X$ is homotopic to $U$.  In particular, the path
$VU^{-1}$ is a path in $Y_1$ whose image in $X$ is null-homotopic.  By
Lemma~\ref{lem:closing-paths}, the image of this path in $Y_t$ must be
closed.
\end{proof}

  We conclude this section with the following corollary.

\begin{cor}\label{cor:generalizedWP}
If $G$ is a group which satisfies the hypotheses of
Theorem~\ref{thm:wt-one-rel} or Theorem~\ref{thm:coh-sc-1}, then each
of its finitely generated subgroups has a decidable membership
problem.
\end{cor}

\begin{proof}
The proofs of these theorems actually show that these groups satisfy
the path reduction hypothesis in addition to the perimeter reduction
hypothesis.  Thus, Theorem~\ref{thm:generalizedWP} can be applied.
\end{proof}


\begin{thebibliography}{10}

\bibitem{ABC90}
J.~M. et~al. Alonso.
\newblock Notes on word hyperbolic groups.
\newblock In {\em Group theory from a geometrical viewpoint (Trieste, 1990)},
  pages 3--63. World Sci. Publishing, River Edge, NJ, 1991.
\newblock Edited by H. Short.

\bibitem{Ba74}
Gilbert Baumslag.
\newblock Some problems on one-relator groups.
\newblock In {\em Proceedings of the Second International Conference on the
  Theory of Groups (Australian Nat. Univ., Canberra, 1973)}, pages 75--81.
  Lecture Notes in Math., Vol. 372. Springer, Berlin, 1974.

\bibitem{BeBr97}
Mladen Bestvina and Noel Brady.
\newblock Morse theory and finiteness properties of groups.
\newblock {\em Invent. Math.}, 129(3):445--470, 1997.

\bibitem{CoLeRi90}
Thomas~H. Cormen, Charles~E. Leiserson, and Ronald~L. Rivest.
\newblock {\em Introduction to algorithms}.
\newblock MIT Press, Cambridge, MA, 1990.

\bibitem{Ep92}
David B.~A. Epstein, James~W. Cannon, Derek~F. Holt, Silvio V.~F. Levy,
  Michael~S. Paterson, and William~P. Thurston.
\newblock {\em Word processing in groups}.
\newblock Jones and Bartlett Publishers, Boston, MA, 1992.

\bibitem{FeHa99}
Mark Feighn and Michael Handel.
\newblock Mapping tori of free group automorphisms are coherent.
\newblock {\em Ann. of Math. (2)}, 149(3):1061--1077, 1999.

\bibitem{Gu93}
Victor~S. Guba.
\newblock The word problem for the relatively free semigroup satisfying ${T}\sp
  m={T}\sp {m+n}$ with $m\geq 3$.
\newblock {\em Internat. J. Algebra Comput.}, 3(3):335--347, 1993.

\bibitem{HrWi01}
G.~Christopher Hruska and Daniel~T. Wise.
\newblock Towers, ladders and the {B}. {B}. {N}ewman spelling theorem.
\newblock {\em J. Aust. Math. Soc.}, 71(1):53--69, 2001.

\bibitem{LySch77}
Roger~C. Lyndon and Paul~E. Schupp.
\newblock {\em Combinatorial group theory}.
\newblock Springer-Verlag, Berlin-New York, 1977.
\newblock Ergebnisse der Mathematik und ihrer Grenzgebiete, Band 89.

\bibitem{McWi-windmills}
Jonathan~P. McCammond and Daniel~T. Wise.
\newblock Coherence tests for one-relator groups.
\newblock In preparation.

\bibitem{McWi-lqc}
Jonathan~P. McCammond and Daniel~T. Wise.
\newblock Locally quasiconvex small cancellation groups.
\newblock In preparation.

\bibitem{McWi-fanladder}
Jonathan~P. McCammond and Daniel~T. Wise.
\newblock Fans and ladders in small cancellation theory.
\newblock {\em Proc. London Math. Soc. (3)}, 84(3):599--644, 2002.

\bibitem{Mo84}
John~W. Morgan.
\newblock On {T}hurston's uniformization theorem for three-dimensional
  manifolds.
\newblock In {\em The Smith conjecture (New York, 1979)}, pages 37--125.
  Academic Press, Orlando, FL, 1984.

\bibitem{Ne68}
B.~B. Newman.
\newblock Some results on one-relator groups.
\newblock {\em Bull. Amer. Math. Soc.}, 74:568--571, 1968.

\bibitem{Pa-alg}
Oliver Payne.
\newblock Private communication.
\newblock 2001.

\bibitem{Pr83}
Stephen~J. Pride.
\newblock Small cancellation conditions satisfied by one-relator groups.
\newblock {\em Math. Z.}, 184(2):283--286, 1983.

\bibitem{Ri82}
E.~Rips.
\newblock Subgroups of small cancellation groups.
\newblock {\em Bull. London Math. Soc.}, 14(1):45--47, 1982.

\bibitem{Sch-cox}
Paul~E. Schupp.
\newblock Coxeter groups, perimeter reduction and subgroup separability.
\newblock Preprint 2001.

\bibitem{Sc73}
G.~Peter Scott.
\newblock Finitely generated $3$-manifold groups are finitely presented.
\newblock {\em J. London Math. Soc. (2)}, 6:437--440, 1973.

\bibitem{Sh91}
Hamish Short.
\newblock Quasiconvexity and a theorem of {H}owson's.
\newblock In {\em Group theory from a geometrical viewpoint (Trieste, 1990)},
  pages 168--176. World Sci. Publishing, River Edge, NJ, 1991.

\bibitem{St83}
John~R. Stallings.
\newblock Topology of finite graphs.
\newblock {\em Invent. Math.}, 71(3):551--565, 1983.

\bibitem{Sw93}
G.~A. Swarup.
\newblock Geometric finiteness and rationality.
\newblock {\em J. Pure Appl. Algebra}, 86(3):327--333, 1993.

\bibitem{Th77}
William~P. Thurston.
\newblock Geometry and topology of $3$-manifolds.
\newblock Lecture Notes, {P}rinceton University, 1977.

\bibitem{vK33}
E.~R. van Kampen.
\newblock On some lemmas in the theory of groups.
\newblock {\em Amer. J. Math.}, 55:268--273, 1933.

\bibitem{Wall79}
C.~T.~C. Wall.
\newblock List of problems.
\newblock In C.~T.~C. Wall, editor, {\em Homological group theory (Proc.
  Sympos., Durham, 1977)}, pages 369--394. Cambridge Univ. Press, Cambridge,
  1979.

\bibitem{Wh78}
George~W. Whitehead.
\newblock {\em Elements of homotopy theory}.
\newblock Springer-Verlag, New York-Berlin, 1978.

\bibitem{Wi-nocore}
Daniel~T. Wise.
\newblock A covering space with no compact core.
\newblock To appear in Geometriae Dedicata.

\bibitem{Wi98}
Daniel~T. Wise.
\newblock Incoherent negatively curved groups.
\newblock {\em Proc. Amer. Math. Soc.}, 126(4):957--964, 1998.

\end{thebibliography}
\end{document}